\documentclass{article}
\usepackage{amssymb,amsmath}
\usepackage[mathscr]{euscript}

\usepackage[cp1251]{inputenc}
\usepackage[english,russian]{babel}

\usepackage{hyperref}

\newcounter{sec}
\newcounter{punct}[sec]

\def\punct{\refstepcounter{punct}{\arabic{sec}.\arabic{punct}.  }}
\def\COUNTERS{\addtocounter{sec}{1}
              \setcounter{punct}{0}
          }

\begin{document}

\def\br{}

\def\sm{\smallskip}

\begin{center}

{\huge\bf The
Kolmogorov reform

\sm

 of mathematical education,
 
 \medskip
 
  1970-1980}

\medskip

\huge\sc

Yury A. Neretin

\end{center}

{\small 
This reform was a part of a transnational movement as 'the New Math' in the United States
or the activity of the Lichnerowicz commission in France. The Russian version of these events involved the whole secondary and high school (4-10 grades).
The leaders of the reformation were Alexey Markushevich and Andrey Kolmogorov.
Their team consisted of best educationalists and best specialists in supplementary education.
 The project offered unqualified support of the Academy of Sciences and of the Ministry of education. However,
 the program of the reform was too optimistic, apparently
 some strategical ideas were unrealizable (and certainly there were no sufficient preliminary experiments for their verification).
 Starting September 1972 Soviet school was depressed, 
 enormous efforts of participants of the reformation  and simplification of the new school curriculum could not improve the situation.
 The reform was stopped by Lev Pontryagin in autumn 1980. The counterreformation required quick
 not-prepared steps and  was not well-ordered.
 In 1981-1985 the heavy crisis was overcome but the mathematical education did not return to the pre-reform
 level.}

\bigskip

{\bf 1. Landscape.}
The%
\footnote{English introduction (p.1--\pageref{eng-fin}).
	 Russian text of the paper, 
	(p.\pageref{rus}--\pageref{rus-end}).
The bibliography (p.\pageref{bib-begin}--\pageref{bib-end}) is  Russian
with few texts in English and Russian texts that have English translations.%
} fifties and the sixties were a time of a florescence  of the
Soviet mathematical education, in particular of education
of schoolchildren. Besides regular school lessons there were 
many other forms of educational activity:

\sm

---  'Circles' (group studies of schoolchildren with somebody elder, a university student, a professor,
a teacher,
an educationalist, etc.).

\sm

--- 'Olympiads' (since 1960 they involved all regions of the Soviet union).

\sm

--- 'Facultatives' (non-obligatory series of additional lessons in schools).

\sm

---  'Correspondence school' since mid-sixties
 (when the University proposed a series of problems
and corresponds with schoolchildren in a country).

\sm 

There were published lot of popular books of a good level (apparently this
editorial program was supervised
by Alexey Markushevich), such books were inexpensive and easily available
in bookshops.

All these forms were possible in other sciences and other branches of
human activity; mathematicians were the most successful among other sciences.

In the beginning of the  sixties first mathematical schools appeared in Moscow.
In 1963 in Moscow, Leningrad, Kiev, and Novosibirsk there were created 
boarding physical-mathematical schools for children from province. In Moscow this was the famous
Kolmogorov school%
\footnote{Колмогоровский интернат, or Восемнадцатый интернат, or Сунец.}.

Nikolay Konstantinov organized an informal structure (Konstantinov's system)
 consisted mainly of students
and postgraduate students,
this structure was involved to work in mathematical schools, circles,
olympiads.

In 1970 Isaak Kikoin and Andrey Kolmogorov initiated a  physical and ma\-the\-ma\-ti\-cal
 journal 'Quant' for schoolchildren. This journal also had a correspondence division.
 
This florescence was outsides regular schools lessons...

\sm

In the fifties school mathematics%
\footnote{During a long period Soviet school of general education  was divided into
a primary school, 1-4 grades (children entered to school being 7 years old),
a secondary school (5-7 grades), and high school (8-10 grades). Majority
of teenagers  leaved schools of general education after 7 grades [[неполное среднее образование]]
and 
went to professional schools (further ways of education were not closed for them). The level of education in 8-10 grades [[старшие классы средней школы]]
 was quite high and this school left very positive memory. After the high
 school, alumni finished education [[полное среднее образование]],
 or entered to high professional schools [[техникумы]], or to universities [[институты and университеты]].
 In years of the Kolmogorov  reform 1970-80 the system was changed. A primary school occupied 1-3 grades,
 a secondary school 4-8 grades, high school 9-10 grades. Half of teenagers
 leaved school after 8 grades and entered to professional schools [[профессионально-технические училища]],
 obligatory 3 years.    
}
 included 'Arithmetics', 'Geometry', 'Algebra', and 'Trigonometry'. 'Geometry' was based on classical approach
 of Euclid type ('Planimetry' and 'Stereometry'). 'Algebra' and 'Trigonometry' included po\-ly\-no\-mials and elementary functions. Traditionally, in the Soviet Union school mathematics (mainly, 'Geometry') was the main learning tool for logical thinking,
 structures of courses  strongly corresponded to achievement of this goal. An additional aim
 of 'Geometry' was
 spatial imagination development. For these reasons school mathematics was not only a natural
 science but a humanitarian discipline of the great importance.      

\sm 

{\bf 2. Reform. Version 1959.} At end of fifties and beginning of sixties, a group
consisting 
of several bright figures of school supplementary education and prominent educationalists
de\-ci\-ded that school curriculum is too archaic and started to develop a
modernization program. 

\sm 

The leader of this group was Alexei Markushevich (1908-1979).
He graduated Middle-Asian University at Tashkent under supervision
of Vsevolod Romanovski (1930), then he learned in the graduate school of Moscow
State University under supervision of Mikhail Lavrentyev. He 
 is best known for his seminal book 'Theory of functions of
complex variables' [[<<Теория аналитических функций>>]], 1950, which is still
 used by mathematicians. Later he published a text-book on complex analysis,
 which was widely used in Russian universities, 'Theory of Abelian functions', and a book on history
 of complex analysis.
 
 However, after 1950 mathematical pedagogics becomes his main domain of activity.
 In 1950 he was chosen a vice-president of the Academy of Pedagogical Sciences
 of the Russian Federation [[АПН РСФСР]] and occupied this position in  1950--1958 и 1964--1975
 (in 1967 this Academy was transformed to the Academy of Pedagogical Sciences of the Soviet Union).
 In 1958-1964 he was  a vice-minister of Enlightenment%
 \footnote{The ministry of school education [[Министерство просвещения]].}
  of the Russian Federation. In 1959-1966 he was a deputy of the Supreme Council
  of the Russian Federation.
  
  He was an influential person in popular book publishing and initiated
  (or participated in initiation of) several
  successful editorial projects (in popular mathematics and not only in mathematics).
   A good selection in such  business
  is very important,  I think that many people that were young  at that time
  could read good books due to Markushevich's efforts.
  
  He himself wrote several short good popular mathematical books, which were translated 
  to many languages...

\sm 

Reformers discussed possibilities to include to obligatory school curriculum
derivatives, integrals, geometric transformations%
\footnote{More precisely, to enlarge volume of this topic.}, vector algebra, 
combinatorics, probability (and statistics), complex numbers,
elements of analytic geometry, elements of programming, the set-theoretical approach,  mathematical
logics (and even more, for instance, groups, rings, fields, or a rigor exposition of geometry%
\footnote{Expositions of geometry of Euclid type were not rigorous...}).

School mathematical curriculums of the fifties contained some archaic elements, cutting
of such elements and some recombinations  gave a
 base for a modernization. On the other hand, mathematics had to develop
 logical thincking and creativity. There existed a danger to break non-obvious 
 psychological constructions created during a long time
 by  Russian teachers and educationalists (as Andrey Kiselev%
 \footnote{Kiselev (1852-1938) graduated from PhysMath of the Imperial Moscow  University
 	and was
 	  a teacher of a school ('gymnasium'[[гимназия]]) in Orel.
 	His textbooks 'Arithmetics', 'Algebra', 'Geometry' 
 	appeared respectively in 1884, 1888, 1892. He himself revised
 	them infinite number of times. These textbooks were widely used
 	until beginning of the thirties. In 1938-40 due actions of mathematical
 	organizations  (as the Moscow mathematical society, the Steklov Institute,
 	 the Academy of Sciences of the Soviet union) modernized versions
 	 of Kiselev's book's were return as basic school texbooks
 	 and were used until 1956.  The last high school gradation 
that used (modified) 'Stereometry' of Kiselev was in 1976.}, Nikolay Rybkin, Nil Glagolev,
 Sergey  Novoselov, Aleksandr Barsukov, and many more)
 using reasonings together with  trial and cut methods...
 
 \sm  
 
 The first version of a new curriculum was proposed in
 1959.  New elements were: derivatives, complex numbers, vectors,
 and more geometric transformations%
 \footnote{More precisely, the purpose was to enlarge 
 a volume of derivatives, transformations, and vectors in the school curriculum}. We know results of experiments
 of the next three decades and can say that this program was more-or-less realistic
 and this was a maximum of possible%
 \footnote{There were other possible maximums
 	and it is non-obvious that precisely this maximum was desirable.
 Also a necessity to achieve a maximum is debatable.}.
 
 Algebraic part of the project was solved in the textbook 'Algebra
 and ele\-men\-tary functions, 9-10 grades' by teachers Evgeny Kochetkov and Elena Kochetkova,
 which was introduced to schools in 1966. This simple textbook included ele\-men\-tary functions,
 derivatives, complex numbers, and induction. Near 1966 well-known geometer Alexey
 Pogorelov created a new geometric course including vectors and geometric transformations.
 This way was rejected by Kolmogorov who intended to write
 a planimetric textbook himself. Pogorelov's book was
 published in 1969-70 as a teacher manual and was
  used as an (unpolished) texbook in 1982 to stop the disaster
 produced by the Kolmogorov reform. It remains one of the most
 popular geometric textbooks in Russian schools up to now. 
 
In 1963 Markushevich (who was a vice-minister at that time)
introduced to schools a new texbook 'Geometry, 9 grade' by Vladimir Boltyansky
and Isaak Yaglom,  it was completely devoted to geometric transformations
and vector algebra. In the next year  minister Evgeny Afanasenko called off
this book as unsuitable for schools of general education. This failure was not analyzed
by reformers (this textbook was quite peaceful comparatively texbooks
of the Kolmogorov project and the failure could be used as a good base for further efforts).
In this year moods have changed and the Academy of Pedagogical sciences created a revolutionary 
program.

\sm 

{\bf 3. Revolution in mathematical education.} 
In 1959 Vladimir Boltyansky, Naum Vilenkin, and Isaak Yaglom
published the paper [БВЯ1959] with a project of a future school curriculum.
They proposed a simultaneous global  trans\-fo\-rma\-tion of all mathematical
subjects in 4-10 grades with changing of a logical structure
of each topic. A promised prize was a possibility to introduce new topics to a curriculum.
Initially, this project was published as a work of abstract art.
However, we know the future. In particular, we know that precisely authors of the paper
represented professional mathematics in the  Markushevich and Kolmogorov team 
(the rest were educationalists), see [Колм1965-1]. 

 A highly experienced educationalist (or highly skilled teacher)  can estimate results of local transformations 
 of a curriculum. Results of a global transformation can be verified only by an experiment.
 For a project of [БВЯ1959] parallel experiments in different grades are impossible,
 because we can involve pupils of $n$-th grades only if they were subjects
 of experiments in 4, 5, \dots, $n-1$ grades. Therefore in the ideally successive  case
 for a verification of the project we need at least 7 years of experiments and at least 1
 year for an independent evaluation of its results.
 
 The authors of Curriculum-1959 took interest in [БВЯ1959]-project...
 At the end of 1964 Markushevich [Мар1964] announced a sketch of a new program
 developed by the Academy of Pedagogical Sciences, detailed were not published.
 However, this was a program based on global transformation of a curriculum, and it looks similarly
 to the Kolmogorov program 1967-8.
 
 \sm 
 
 Kolmogorov jointed to reformers in the same 1964, in the natural way he became a lieder
 of the project. Versions of the program were published in 1966 and 1967. The authors of the version-1967 were Vladimir Boltyansky, Andrey Kolmogorov, Yury Makarychev, Aleksey Markushevich,
 Galina Maslova, Konstantin Neshkov, Aleksey Semushin, Antonin Fetisov, 
 Aleksandr Sheresevsky, Isaak Yaglom.

 There was lot of public admiration for new ideas and 
  few  publications with a professional analysis%
  \footnote{
 The professional  journal 	'Mathematics in school' was controlled by reformers.
 There were published three collections of opinions  in
  [О-прог-1967-3], [О-прог-1967-4], [О-прог-1967-1]. Few negative
letters are contained  in [О-прог-1967-4]. However, several texts in that collection starting from
words 'This is a great project' contained very unpleasant observations.}.

Mathematical branch of the Academy of sciences of the Soviet Union
sup\-por\-ted the program. In 1968 the new curriculum was approved be the Ministry
of Enlightenment of the Soviet Union. In September 1970 fourth grades
of
schools of general education of USSR began work according the new curriculum...

\sm  
  
{\bf 4. Few  comments on Curriculum-1968.}  

A) Clearly, {\sc the new curriculum 
was overcommitted}. Apparently, the authors thought
that it will be possible to omit some topics when strictly necessary.

\sm 

B) {\sc The set theoretical ideology}. The word 'set' was new in school,
but this word is neutral, usage of terms 'union' and 'intersection'
also is semi-neutral. 
The set theoretical ideology was a more serious thing. I'll say few examples.

\sm 

{\footnotesize 
1. The authors of 'Algebra-6' [МММ] define binary relations, explain them
for finite sets, after this they define functions in terms of binary relations,
and  use functions of real variables. 

\sm

2. {\sc Definition} ('Geometry 9-10' [КСЯ]).
A {\it vector} (or a parallel translation), defined by a pair $(A,B)$
of non-coinciding points is a transformation of space that send each point
$M$ to the point $M_1$ such that the ray $MM_1$ is codirectional
to the ray $AB$ and the distance $|MM_1|$ equals distance $AB$. 

\sm 

3. From 'Geometry, 6' [КСНЧ], first edition, pp. 11-13:

{\sc Definition.} A point $X$ is {\it lying between} points $A$ and $B$,
if it differs from both of them and 
$$
m(A,X)+m(X,B)=m(A,B).
$$

{\sc Definition.} The set consisting 
of two non-coinciding points  $A$ and $B$
and all point between them is called a {\it  segment} $AB$.

\sm

{\sc Definition.} A {\it ray} $OA$ is the set of points
consisting of all points of the segment   $OA$
and all ponts  $X$, for which  $B$ lyes between  $O$ and $X$.

\sm

... The {\it line} $AB$ consists of 

1) all points of the segment  $AB$;

2) all points $X$, for which  $B$ lyes between $A$ and $X$;

3) all points $Y$, for which  $A$ lyes between $B$  $Y$.

}

\sm 

Such kinds of statements and arguments were spread over all courses
of 6-10 grades...  However, I think that  the fatal mistake of reformers was another.

\sm

C) {\sc The transformation of geometry}. Kolmogorov wrote
the textbook 'Geometry. 6-8 grades' (planimetry)
together with eminent educationalists Ale\-ksandr Semenovich, Fedor Nagibin,
Rostislav Cherkasov. The main claim was to base 
 elementary geometry on geometric transformations and vector algebra.
 This idea was many times and many countries%
 \footnote{It seems that in Soviet Union where were 
unsuccessful attempts to realize it. I did not investigated
these topic.} repeated starting 1900. From a point of view
 of abstract logic this possibility is doubtless. But  both topics are 
 difficult for teaching and themselves require certain base. There were no
 experiments confirming that this is possible.
 
 \sm 
 
 {\bf Experiments.} The Kolmogorov curriculum was approved in 1968
 when a package of textbooks was not written.
 In this 1968  experiment with 'Mathematics, 4 degree' started. 
 After the first year the curriculum of 4 degree was seriously cut.
After the second year it was cut second time...
 
 Experiments with systematic courses of geometry and algebra (6 degree) were started
 in September, 1970. In the same 
  moment there was started a total introduction
  of the Kolmogorov curriculum to schools of general education
 of the whole country (starting 4 degree). Bridges were burn.
 
 \sm
 
 The so-called 'Kolmogorov program' was a strange victory of a theory over a practice.
 However, such victory happened not only in Russia but also
 in US and France%
 \footnote{From a French school textbook (quotation from [Маш2006]):
 {\sc Definition.} {\it 
 	We say that an orthogonal endomorphism of $\phi$ of $E_3$ is a vector rotation to express the fact that the subspace
 	of variables invariant under $\phi$ has dimension 1 or 3.
}}.
 
 \sm 
 
 {\bf Reality.}  Two years 1970-1971 and 1972-1972 were seemingly calm.
 In Sep\-tem\-ber, 1972, Soviet school was shocked. 'Geometry, 6'
 of Kolmogorov and 'Algebra, 6' edited by Markushevich came to school.
 Pupils understood nothing...
 
 The following 8 years were a time of enormous efforts.
 Reformers claimed that teachers are not qualified enough and started
 a wide program of continuous education of teachers. Simultaneously,
 they simplified textbooks and tried to improve them.
 Situation did not improved%
 \footnote{An example. All grades of 10-year school had to pass final exams 
 	in core subjects. Normally, these exams were very friendly, however they
 	were an important part of the education machine. Starting 1977 final
 	geometric exams
were canceled since the textbook 'Geometry-9-10' [КСЯ] 	edited
by Zalman Skopets
 was completely improper.}.
  An introduction
 of each new party of textbooks
 (7, 8, 9, 10 degrees) produced new difficulties...

 Since teachers were denounced, nobody heard them. However there existed
 millions of angry parents. In 1977 first grades, which learned new school
 mathematics, came to Universities, and university professorship
 was not happy.

 22 December 1977 the Council of Ministers and the Central committee of
 the Communist party published a resolution [ЦК1977] containing the following words:
  
  \sm 
  
 {\it  In certain cases school curriculums and texbooks are
  overloaded with un\-ne\-ces\-sa\-ry information and secondary materials, this
  is an obstacle for
   a development of student's  skills of independent creative work.}

\sm

This was a usual cloudy style of such resolutions of that time. More
in\-for\-ma\-ti\-on was distributed through closed channels. Here main of 'certain cases'
 was obvious without additional detailing.
 
 In fact, the Goverment proposed to mathematicians to simplify curriculum and textbooks...
 
 \sm 
 
 The main problem of reformers was the texbook on planimentry
 by Kol\-mo\-go\-rov, Semenovich, Cherkasov. This texbook was based on 
 ideological mistakes and can not be simplified, a simplification
 required a complete rewriting of the textbook and a counterreformation
 by hands of reformers.
 
 \sm 
 
 {\bf The end.} In December 1978 the mathematical branch
 of the Academy of sciences of the Soviet Union voted against the
 Kolmogorov curriculum. Their resolution was not published...
  Positions of reformers seemed to be strong. 
 
  Several academicians
 (Alexei Pogorelov, Andrey Tikhonov, Aleksandr Ale\-ksan\-drov)  started to prepare new text-books.
  
  Lev Pontryagin found a way (see [Пон2008]) to 'Communist', the main journal of the Communist party.
  His pitiless paper [Пон1980] 'On mathematics and quality of its teaching' was published
  in September 1980.  The reform was finished.
  
  \sm 
  

 \bigskip 
 
$$\times\qquad\times\qquad\times\qquad\times\qquad\times\qquad\times$$

\bigskip 

The text after the content is  Russian [[После оглавления идет основной текст статьи
по-русски]].

\bigskip

{\sc
	
\begin{center}
	
	Content
	
	\end{center}	
	
1. Two definitions (instead of a preface)

2. Two myths

3. A reform or a revolution?

4. Dieudonne and Lichnerowicz

5. The process has begun

6. The successful experiment

7. Clash with reality

8. Overkeel

9. Dizzy with Success (instead of an afterword)

10. Addendum. The team of Markushevich and Kolmogorov

11. Bibliography
\label{eng-fin}

}

\vspace{33pt}

\label{rus}
Попытаюсь описать картину, которая складывается из моих личных наблюдений над процессом введения колмогоровских учебников, 
обсуждений в математической тусовке второй половины 70х
годов, поисков  в старых изданиях и в современной литературе.
Как получилось, что реформа 1968-1980гг, возглавленная самим великим Колмогоровым, закончилась таким звонким фиаско?
И как вообще эта реформа случилась?

\sm 

{\sc Историография} по этой теме небогата (если не брать в рассмотрение многочисленные  статьи,
сводящиеся к моралистическим или эмоциональным оценкам).
А.~М.~Абрамов, ученик Колмогорова по педагогической линии и известный деятель образования 
90х-00х годов, последовательно отстаивал идею правильности и разумности Реформы,
см. [Абр1988], [Абр2003], [Абр2010]. С ультраконсервативных позиций критикует реформу И.~П.~Костенко
[Кост]. Историю реформы обсуждал Ю.~М.~Колягин в книге <<Русское математическое образование>> [Коля2001].
Он сам был деятельным участником событий, сторонником Реформы во время ее подготовки
и начала ее проведения, автором нескольких реформистских пособий в 1969-1975гг. В
начале 1978г.
он предстает в качестве одного из лидеров контр-реформистской комиссии при Минпросе РСФСР.
 Вместе с О.~А.~Саввиной
он опубликовал сборник документов <<{\it Бунт российского министерства и отделения математики АН СССР. (Материалы по реформе школьного математического образования 1960-1970-х гг.)}>> [КС2012].
Самый   интересный документ там -- {\it Стенограмма Общего собрания Отделения математики АН СССР,  5 декабря 1978}, ниже цитируется как [Стен]. C точки зрения содержательного анализа  реформы важна статья А.~Л.~Вернера [Вер2012] 

Основные исторические источники, использованные в данном исследовании --
это математико-педагогические  публикации 50-80х годов и школьные учебники тех лет.
В архивах я не работал, опубликованных архивных источников, кроме [КС2012],
мне, к сожалению, неизвестно.


 
%

\bigskip                                                       
                                                          
{\sc Оглавление
	
	\sm

1. Пара дефиниций (вместо предисловия) 

2. Два мифа

3. Реформа или революция?

4.  Дьедонне и Лихнерович

5. Процесс пошел

6. Успешный эксперимент

7. Столкновение с действительностью

8. Оверкиль

9. Головокружение от успехов (вместо послесловия)

10. Дополнение. Соратники Маркушевича и Колмогорова.

11. Список литературы%
\footnote{ В списке литературы почти все издания открываются по ссылке.}}

\section{Пара дефиниций (вместо предисловия)}

\COUNTERS

Я учился по этим учебникам, когда они  были экспериментальными и еще не пошли в массовую школу.
Об этом расскажу ниже, а пока определение прямой,
которое меня удивило тогда и не перестают казаться удивительными и сейчас. Из учебника геометрии за 6-ой класс, в самом начале систематического курса планиметрии.

\sm 

{\sc Определение.} Точка $X$ лежит {\it между} точками $A$ и $B$, если она отлична от каждой из этих точек,
и для нее соблюдается равенство
$$
m(A,X)+m(X,B)=m(A,B).
$$

{\sc Определение.} Множество, состоящее из двух несовпадающих точек $A$ и $B$ и всех
точек, лежащих между ними, называется {\it отрезком} $AB$.

\sm

{\sc Определение.} {\it Лучом} $OA$ называется множество точек, которое состоит из всех точек
отрезка $OA$ и всех точек $X$, для которых $B$ лежит между $O$ и $X$.

\sm

И, наконец, прямая,

... Таким образом, {\it прямая} $AB$ состоит из

1) всех точек отрезка $AB$;

2) всех точек $X$, для которых $B$ лежит между $A$ и $X$;

3) всех точек $Y$, для которых $A$ лежит между $B$ и $Y$.

\sm

Это из первого издания учебника 
А~Н.~Колмогоров,
Ф.~Ф.~Нагибин, А.~Ф.~Семенович и Р.~С.~Черкасов  <<{\it Геометрия, 6}>>, страницы 11-13.
В дальнейшем учебник менялся, и эта премудрость несколько смягчилась
Я еще запомнил, что в учебнике (это уже за 7 класс) очень сложно объяснялось понятие вектора.
Но лучше процитируем определение вектора из учебника  
В.~М.~Клопский, З.~А.~Скопец,  М.~И.~Ягодовский  {\it Геометрия}, 9-10 классы. Под редакцией З.~А.~Скопеца,

\sm 

{\sc Определение вектора.}
{\it <<Вектором (параллельным переносом)>>}, определяемым парой $(A,B)$ несовпадающих точек
пространства, называется преобразование пространства, 
при котором каждая точка $M$ отображается на такую точку $M_1$, что луч $MM_1$ сонаправлен с лучом $AB$ и расстояние $|MM_1|$ равно расстоянию 
$|AB|$.

\sm

Это определение было ва всех изданиях учебника и дожило до конца Реформы.
Именно с него  Понтрягин начал свою убийственную статью [Пон1980]:

\sm 

{\it В этом сплетении слов разобраться нелегко, а главное — оно бесполезно, поскольку не может быть применено ни в физике,
ни в механике, ни в других науках.

Что же это? Насмешка? Или неосознанная нелепость? Нет, замена в учебниках многих сравнительно простых,
наглядных формулировок на громоздкие, нарочито усложненные, оказывается, вызвана стремлением... усовершенствовать (!) преподавание математики.

Если бы приведенный мною пример был только досадным исключением, то ошибку, по-видимому, легко можно было бы устранить. 
Но, на мой взгляд, в подобное состояние, к сожалению, пришла вся система школьного математического образования...}

\sm 

Эта пара определений являются  своего рода вершинами, но такие пики не вырастают на ровном месте...

Не вызывает сомнений, что Колмогоров лично работал над учебником
<<Геометрия, 6-8 классы>>. Его соавторы  Ф.~Ф.~Нагибин, А.~Ф.~Семенович и Р.~С.~Черкасов
были известными и заслуженными педагогами. Написанная ими книга не могла не иметь многих достоинств.
Тем не менее эта книга не годилась  именно для того, для чего она была написана - для того, чтобы быть
начальным учебником геометрии. 


Впрочем, совсем уж сильно удивляться этому не стоит. Боюсь, что слишком многие, учившиеся
в институтах и университетах в последние лет 40,
сталкивались в качестве слушателей с вузовскими  курсами математики, 
(иногда будучи не лишенными многих достоинств),
однако не рассчитаны ни на какого живого человека. 
Они в своем роде совсем просты: если отучиться на мехмате МГУ, то в курсе всё проще простого, а если не отучиться - 
то понять его нельзя, потому как там опущены
некоторые необходимые для понимания куски и потому что туда напихано много ненужных украшательств.
 Комплект школьных учебников эпохи Колмогорова был одним из ранних образцов <<курсов ни для кого>>.
О том, как и почему  это получилось, пойдет речь
в настоящей работе.

\section{Две мифа}

\COUNTERS

Прокомментирую две сентенции, часто повторяющиеся при обсуждении предмета и искажающие современное  восприятие этой истории.

\sm 

{\bf\punct Легенда о началах анализа.}
Есть распространенное мнение%
\footnote{Которое делится на черную балладу о зловредности начал анализ для школы, и на оптимистическую трагедию: несмотря 
ни на что реформа победила!}, что элементы анализа были введены в школу
благодаря реформе Колмогорова 1970-1980. Это не так.
В 1966-1975гг. основным учебником алгебры для 9-10 классов  
был Е.~С.~Кочетков, Е.~С.~Кочеткова <<Алгебра и элементарные функции>>. Беру вторую часть учебника (10 класс), издание 1968 года
[КК1968]. Объем
учебника - 284 страницы. Первая глава посвящена  тригонометрии - 64 страницы.
Дальше - показательная функция и логарифм - 60 стр., функции и пределы - 66 стр, производная - 50 стр., комплексные числа - 
34 стр, метод математической индукции - 19 стр.

Можно спорить о деталях, о том, что показательную функцию и производную хотелось бы пораньше иметь на  физике и т.п. Но
и производная, и комплексные числа в школе были до колмогоровской реформы%
\footnote{М. Ю. Колягин [Коля2001] пишет, что производная появилась в школе в 1954г.
Я отдельно этот вопрос не исследовал.}. В результате колмогоровской
реформы и ее отката комплексные числа из программы выпали (и индукция вместе с ними),
эта тема, как мне кажется, довольно безобидная, была в итоге заменена на достаточно сомнительный <<интеграл>>. 
Вот и весь итог... Что касается  интеграла, то я с 1983-2001гг работал на факультете Прикладной математики
МИЭМ (в те времена весьма приличном), было заметно, что студенты были неплохо знакомы с производной до обучения в институте, а следов их знакомства
с интегралом что-то видно не было.

Начиная с этого места, многократно цитируется  стенограмма [Стен1978]
Общего собрания Отделения математики АН СССР 5 декабря 1978 года, изданная Колягиным и Саввиной в [КС2012]
(это собрание было посвящено катастрофическому положению с учебниками математики):

\sm

{\sc ШАБУНИН} (представитель физтеха). {\it Элементы высшей математики в программе средней школы
появились значительно раньше - мы имели 20 лет назад эти элементы в меньшем объеме,
а в 60е годы они были, может быть, чуть меньше.} 

\sm

{\bf \punct Легенда о Киселёве.} Учебник Киселёва сейчас является символом дореформенной школы,
Золотого Века, когда геометрия
была понятной и приятной.

\sm

{\sc Арнольд}, 2009: {\it  Я бы рекомендовал в преподавании в школе вернуться к Киселёву.}

\sm 

У Киселёва было три наиболее известных учебника, <<Арифметика>> (1884), <<Алгебра>> (1888), <<Геометрия>> (1892), 
они оказались исключительно удачными и живучими. История с этими учебниками была
такая. Они были популярны до  революции и продолжали использоваться в 20е годы. На 1929-30гг.
пришли тяжелые потрясения системы образования, и тогда учебники Киселева (кроме <<Алгебры>>)
из обихода выпали. Новые учебники, по мнению математической общественности, 
были неудовлетворительны. После нескольких лет писания писем в Наркомпрос
в конце 1936-начале 1937 года последовало несколько выступлений математических организаций.
В качестве образца приведу отрывок из резолюции 
 Московского математического общества [ВММО1937] от 12 апреля 1937 года (см. также [Нер2016]):

 \sm
 
 {\it 
По вине невежественного руководства со стороны Управления средней школы Наркомпроса, в частности по вине А.~И.~Абиндера, 
учебная литература по математике находится в настоящее время на чрезвычайно низком уровне.
Управление средней школы Наркомпроса, получая в течение ряда лет и со стороны научных организаций
и со стороны учительства сигналы о безграмотности стабильного учебника геометрии Гурвица и Гангнуса, 
никакой подготовительной работы для замены 
этой безусловно вредной книги не вело. Книги Гурвица и Гангнуса должны быть изъяты и ни в каком случае не переиздаваемы.
Временно стабильным учебником геометрии должен быть объявлен курс Киселева под редакцией Н.~А.~Глаголева.
 Учитывая, что составление оригинального учебника геометрии потребует времени и что курс Киселева на ближайшее время
 может удовлетворить потребности преподавания (книга написана безусловно грамотно 
 и усовершенствовалась на протяжении 40 изданий), собрание считает необходимым объявление конкурса
 с длительным сроком (не менее трех лет) на 1) стабильный учебник по геометрии, 
 2) курс элементарной геометрии для учителей.
}

 \sm
 
 В итоге  с 1938-40 видоизмененный 
<<Киселев>> вторично был  пущен в дело в качестве стабильных учебников.
<<Арифметика>>  была 
переработана аж самим А.~Я.~Хинчиным, <<Алгебра>> -- А.~Н.~Барсуковым, а <<Геометрия>> --  Н.~А.~Глаголевым
(с характером правки Хинчина можно ознакомиться по его статье [Хин1941]).

\sm

Следующее преобразования последовали в 1955.
Учебником арифметики в 5-6 классах стал  <<Шевченко>>. В средних
классах стабильным учебником алгебры  стал <<Барсуков>>, а геометрии --  <<Никитин>>.
Стереометрия осталась киселевской, с алгеброй старших классов я не разобрался,
<<Тригонометрия>> Рыбкина сменилась на Новосёлова.

\sm 

C 1966г. в старшие классы  
в качестве основного 
  учебника алгебры пошел Кочетков-Кочеткова (упомянутый чуть выше), последний выпуск по нему состоялся в
  1976г.

\sm

От Киселева к моменту Колмогоровской реформы оставалась лишь <<Стереометрия>>, и то отчасти реформированная.
В легенде правда то, что учебники Киселева продержались долго благодаря их замечательным достоинствам,
и правда  то, что предреформенные учебники были простыми.

\sm

Да, кстати, информация  для ультраконсерваторов, призывающих <<долой производную>> и <<назад к Киселёву>>
(см. многочисленные статьи И.~П.~Костенко [Кост]
и их клоны). В 1920х годах издавался еще один учебник Киселёва [Кис1925] (в каталоге РНБ присутствует 8-е издание 1930 года)
для <<трудовой школы>>,
 <<Элементы алгебры и  анализа>>:

 \sm
 
{\it  настоящая книга содержит в себе элементы анализа бесконечно-малых, с его
применениями к вопросам элементарной геометрии и начальной механики, и краткие
сведения по аналитической геометрии, без которых элементарный курс математики был
бы неполным;

в конце книги помещены некоторые дополнения к обычному курсу алгебры,
напр., теория соединений, бином Ньютона,.. и другие.}

\sm 

Комплексные числа там тоже были...


Так что реформу 1959 года, о которой пойдет речь в следующем
разделе, можно при желании охарактеризовать словами <<Назад, к Киселеву!>>.

И еще для ультрапочвенников и противников <<растленного влияния>>:  Киселёв Андрей Петрович 
(см. его биографию в [Анд1967], [АА2002]) кончал московский физмат, 
 при создании своих учебников использовал многочисленные французские, немецкие и английские книги%
 \footnote{Cм. длинные  списки 
в  предисловиях
к  [Кис1892], [Кис1925].}.

\section{Реформа или революция?}

\COUNTERS

{\bf \punct Реформа сверху.}
Идея преобразования школьной математики в середине 60х годов  исходила из двух источников.  Один из них был довольно умеренным,
но он способствовал продвижению идеи реформы.

Тогда совершался переход ко всеобщему среднему образованию.
Раньше была школа-семилетка, а для меньшей части -- школа-десятилетка%
\footnote{На короткое время при Хрущеве 1958-1965, семилетка была заменена
на восьмилетку, а десятилетка на 
одинадцатилетку.} (и как раз эта советская
старшая школа оставила  весьма положительные
воспоминания у тех, кто в ней учился). Потом перешли от семилетки прешли к восьмилетке, а потом уже ко всеобщему среднему
(которое, впрочем, делилось на общеобразовательную школу и профессионально-техническое образование).
Переход ко всеобщему среднему ставил вопрос о новых учебниках. Минпрос  связывал с эти определенные надежды.
Цитирую стенограмму [Стен1978] 1978 года:

\sm

{\sc КОРОТКОВ В.М.} (зам. министра просвещения СССР):...{\it  Отмечался ненужный концентризм построения курса,
т.е., своеобразный повтор курса... 
Есть педагогический концентризм,
а есть концентризм, связанный с изменением ступеней школы. У нас были и 7-летние и 10-летние школы, и та, 
и другая должны были дать соответствующую подготовку.
 В связи с переходом на единую общеобразовательную школу встала проблема ненужного концентризма...}
 
 \sm
 
 В принципе, за счет устранения <<концентризма>> можно было добиться определенного выигрыша%
 \footnote{Кстати, в <<Геометрии>> Никитина 1956-1971 за средние классы были сведения по Стереометрии,
 думаю, 
что это упрощало дальнейшее восприятие школьниками учебника Киселева.
Начало <<Стереометрии>> Киселева, по-моему, не очень простое,
и никакой моральной подготовки к этому в предшествовавшей ей до 1956 года Планиметрии Киселева не было.}
 (но выигрыша небольшого),
 с другой - переход на всеобщее среднее образование должен был
 породить проблемы, которые в государственных,  образовательно-министерских  и в АПН-овских кругах, быть может,
 и не вполне просчитывались%
 \footnote{Если смотреть из сегодняшнего дня, то возникает предположение, что этот переход мог
 быть вынужденным и вызванным проблемами с молодежью, не пошедшей в
 старшую школу.}.
 
 \sm
 
 {\bf \punct Реформа снизу.}
 Избавление от <<концентризма>> не требовало революции: переход от 11-летней школы к 10-летней тоже сопровождался  перекройкой учебников,
 ее мирно провели. Мирно можно было  провести и переход к всеобщему среднему образованию.

Однако тогда на уроках физики тогда говорили об электронах и ядрах,
на химии - о строении молекул. В то же время алгебра не выходила за рамки 17 века, а геометрия за рамки третьего века до нашей эры
(правда наше трехмерное пространство с тех пор тоже сильно не изменилось, как и наши представления о нём).
Идея модернизации курса математики более сложна, чем в прочих предметах, потому что в тех присутствует  описательная часть, а математика - логическая 
конструкция, из под которой трудно изымать  нижние  элементы 
(и, кроме того, школьная математика несет на себе большую гуманитарную нагрузку). 
В последние десятилетия профессиональная математика столкнулась с проблемой модернизации
 университетских программ и университетских учебников (куда ушли времена <<Курса высшей математики>> В.~И.~Смирнова?).
 Тогда в начале 1960х математики столкнулись с  этой головоломкой на уровне школ. К сожалению, они вовремя не поняли, 
 сколь эта головоломка страшна%
 \footnote{См., впрочем, [Луз1920]}.
 
 В 50х - начале 60х годов
многие известные математики и многие известные педагоги считали, что программу средней школы следует реформировать,
устранить из нее архаичные и искусственные
элементы и вместо этого ввести новые разделы.

Процитирую А.~Я.~Хинчина [Хин1961], хотя написано это было раньше обсуждаемого времени:

\sm

{\it 
Как-то   мне   пришлось    спросить    несколько   опытных   учителей    пятых   классов   о   том, 
какой   примерно   процент  учащихся   действительно   
научается   решать   арифметические    задачи,   не   являющиеся    простыми    
вычислительными    примерами,    т.    е.    такие,   где    способ   решения,    как    
бы   прост   он   ни   был,   должен   быть   найден   самим  учащим...
  Конечно,   решив   целый    ряд    совершенно   однотипных   задач,   ученик   без   
труда    решит    задачу    в    точности    того    же    типа   (этим    объясняется    
отсутствие   сплошных   провалов   на  экзаменах  и  контрольных   работах);   
но  добиться,   чтобы    ученик    самостоятельно    нашел    решение    задачи    
нового,    хотя    бы    и    очень    простого    типа,  —  это,    по   единодушному   
мнению    учителей,   есть    дело,   у дающееся     только   в   самых     исключительных    случаях...

Сопоставим   с   этим    другой    факт:   хорошо   известно,   что   большая   
часть   наших    ученых   математиков,   как   правило,   становится    в   тупик   
перед  задачами  элементарной   арифметики.   Я  лично  охотно   признаюсь,   
что    всякий    раз,    когда    ученик   пятого   класса   среди   моих    знакомых    
просил   меня    помочь   решить    арифметическую    задачу,   дело   это   для   
меня   оказывалось    весьма   тяжелым,   а   подчас   я   терпел   и   полную   неудачу.  
Я,   как  и  большинство   моих  товарищей,   легко  решал,   конечно,   
предложенную   задачу   естественным   алгебраическим   путем   (т.   е. 
составлением    уравнения    или    системы    уравнений);    но   ведь   надо   было   
во  что   бы  то   ни  стало  обойтись  без  алгебраического  анализа!    Обычно    
если   мне   в  конце    концов    удавалось    найти    такое   решение,   оно   в
же   оставляло   меня    неудовлетворенным:     моя   научная   совесть   неумолимо    подсказывала     мне,  
что    тут    остается    какой-то   туман,   не   всё   
ясно.

    В   результате,    как   правило,   и   ученика   мое  решение   не   удовлетворяло,  
    и   он   явно   лишь    из    вежливости    принимал   его.    Иногда    в    
таких   случаях    я    потом    пытался    узнать,   как    же   объяснил   решение   
задачи   учитель?   Должен   признаться,   что   и   в  рассуждениях   учителя   
для   меня   почти   всегда   оставался   тяжкий   элемент   ненатуральности    и    
искусственности...

.... Кстати,  хорошо  известно  и  многократно  отмечалось,   что,  как   правило,   
ни  оканчивающие   школу,   ни  студенты   педвузов,   ни  начинающие    учителя    
(ни,    прибавим    от    себя,   научные   работники)   не   умеют   решать   
арифметических   задач,   да   и  вряд   ли  на  всем  свете  кто-нибудь   умеет   
решать   их,   кроме   учителей   пятых   классов.}

\sm

Статья Хинчина%
\footnote{
С легкой руки И.~П.~Костенко [Кост] по интернету гуляют сентенции типа <<программа реформы была в основном сформулирована Хинчиным
в конце 30х годов>>.
Есть посмертный сборник педагогических работ Хинчина [Хин1963], Александр Яковлевич был настроен на введение в школу элементов анализа
(которые -  напомню - были и у  Киселёва), но усмотреть сходство довоенных взглядов Хинчина  с колмогоровской реформой
я лично не могу, разве что и там, и там речь шла о реформе. Но реформы бывают разные.}
была найдена при посмертном разборе его бумаг, Б.~В.~Гнеденко предположительно датировал ее 1938-39 годами.
Речь идет об изощренных задачах типа <<из одной трубы в бассейн вливается, из двух выливается>>, которые современным людям,
как правило, 
и не знакомы. Я лично помню рассказы старших людей об ужасе, который вызывали
эти многовопросные задачи.
Плюс от них все же был, после них школьники  встречали появление икса и игрека как избавление.
Не знаю, когда в точности это ушло,
уверенно могу сказать, что в конце сороковых годах оно еще было, и что,  будучи школьником,
я видел вбивающий в смертную тоску пожелтевший сборник 
подобных задач%
\footnote{При написании настоящей работы, мне не удалось такие книги найти.
Кстати, ни учебник Киселева, ни унаследовавший ему учебник Шевченко такого не требовали.}
(были ли они при мне, тоже не помню, потому что с третьего класса учился по экспериментальным учебникам).

В связи с освященностью подобных задач древней традицией у читателя может возникнуть вопрос, 
не было ли в этом высшей сермяжной правды?
Вряд ли.  Понятно, что текстовые задачи полезны, и полезно
в определенных дозах их решать без привлечения уравнения. 
Стоит, однако, заметить, что в различных сегментах системы образования, развивающихся на своей собственной основе без внешних коррекций, с легкостью возникают
странные
опухолевые структуры. Ныне здравствующие могут вспомнить математику вступительных экзаменов, высокоразвитую индустрию раскрытия неопределенностей,
или такие особые разделы математики как задачи С3, С4, С5, С6 современного ЕГЭ. Развивается все это само собой для удобства тех или иных
образовательных структур и групп участников образовательного процесса.  Раскрывать тучи неопределенностей или безуспешно искать $\delta$ по 
$\epsilon$ 
на самом деле весьма комфортно и избавляет препода от более утомительных родов деятельности. Скорее всего,  передозировка текстовых арифметических задач
развилась сама собой в пограничной зоне между младшей и средней школой и развилась по правилам младшей школы.

\sm

При чтении педагогико-математической литературы конца 50х-начала 60х годов все время 
наталкиваешься на идеи, исходящие от разных авторов, о включении
  в школьную программу производной, интеграла, 
 геометрических преобразований, векторной алгебры,
комбинаторики, теории вероятностей (и статистики), комплексных чисел, элементов аналитической
геометрии, теоретико-множественного подхода
и математической логики.
В принципе почти любой отдельный элемент списка  впихнуть
было можно,
и можно было сделать это относительно разумным образом (кроме статистики, интеграла и логики).
Но  разными другими мечтами было необходимо пожертвовать. И решить это надо было очень жестко...

Понятно, что при попытке модернизации математических программ основной вопрос состоит не в том,
что туда можно ввести, а в том, что для этого надо сократить.
Так или иначе, различные резервы для сокращения, которые могли бы дать
возможность для  модернизации школьной программы, 
не приводящей (или почти не приводящей) к ее
фактическому усложнению, в 50е годы были. Но были они не так уж велики.

\sm

{\bf \punct Программа-1959 и попытки ее воплощения в жизнь.}
В том же 1959 году, 
был разработан проект новой школьной  программы по математике
(последовательные варианты в  [Прог1959], [Прог1960], [Прог1961]). 
Сюжет этот  не изучен (и, как будто, почти выпал из истории). По моему, в программе присутствовали нереализуемые намерения, и неудачные детали,
но в целом  она, скорее, была положительным явлением (и осталась бы таковым, если бы итоги экспериментов 
первой половины
60х были бы адекватно оценены).
В школу предлагалось ввести пределы, производную, комплексные числа, геометрические преобразования, векторы%
\footnote{Более точно,предполагалось учеличить дозы геометрических преобразований и векторов.
Если говорить не о деталях, а самой идее, то она была благой и достижимой.}.
Был объявлен  конкурс учебников (на рецензируемых рукописях не было фамилий авторов,
а вместо этого стояли <<девизы>>).

В рамках Проекта-1959 многое удалось сделать. Был отменен предмет <<Тригонометрия>>,
за счет этого был сделан новый предмет для старших классов - <<Алгебра и элементарные функции>>.
Был запущен в оборот учебник Кочеткова и Кочетковой
-- простой учебник, о котором уже говорилось --
где многие желания удалось мирно претворить в жизнь. Этот учебник, кстати, занял первое место в конкурсе
(см.[ГП1965]) 
в номинации <<алгебра для старших классов>>. Председателем комиссии по алгебре был зав. кафедры алгебры
мехмата МГУ А.~Г.~Курош, а учебник выходил под редакцией О.~Н.~Головина, работавшего на этой кафедре.
(По другим номинациям первых мест не было).

\sm

С геометрией получилось хуже. В 1963г. в массовую школу (тогда была одиннадцатилетка)
пошел учебник

\sm

Болтянский В.~Г. , Яглом И.~М. {\it Геометрия. Учебное пособие для 9 класса средней школы.}
- М.: Учпедгиз. 1963 (тираж  2 300 000), 1964 (тираж 2 000 000), 
издана также в Киеве, 1963, тираж 230 000.

\sm

Следов его предварительной обкатки
(с помощью  экспериментальных учебников с малыми тиражами) по библиотечным каталогам [Кат-РНБ], [Кат-МГУ]
не видно... Трудно усомниться 
в том, что учебник
провел лично А.~И.~Маркушевич, который в 1958-1963гг был замом министра просвещения РСФСР.

\sm

Посвящен учебник  геометрическим преобразованиям с добавлением векторной алгебры -- это предполагалось
изучать весь 9 класс (46 часов). 
Что можно сказать об этой книге? Абстрактно, эта книга  хорошая,
она  была  полезной для увлекающихся школьников,
для учителей, ее можно было использовать как дополнительный учебник в мат.школах.
Но из этого не следует, что она годилась как базовый школьный учебник.
Ни эта книга, ни ее
клоны больше не переиздавались. Следов ее воздействия в ближайших по времени  учебниках для средней школы 
не видно. В официозной реформационной статье [БМ1975] есть такая фраза (без прямого упоминания данного учебника):

\sm

{\it Понятие интеграла и его применения совсем не вошли в программу} [эта  фраза о <<Кочеткове>>]. {\it
Аналогично обстояло дело с геометрическими преобразованиями. Их изучение в 
IX классе оказалось некоторой не совсем оправданной <<надстройкой>> над курсом планиметрии.
Идея геометрических преобразований очень мало использовалась
в самой геометрии.}

\sm

А вот что говорит современный антиреформистский автор Ю.~М.~Колягин%
\footnote{Тогда, как он и сам честно подчеркивал, он участвовал в реформации.}
[Коля2001]:

\sm

{\it учебник в 1958 году} [Ю.~Н.: очевидная опечатка] {\it был сразу внедрен в массовую школу 
(через год он был отменен Министерством просвещения как
непригодный для массовой школы).}

\sm

Введение геометрических преобразований и векторов в умеренных дозах не представляло неразрешимой  проблемы, 
такой курс был написан  А.~В.~Погореловым в середине 1960х. Но тогда уже в игру вступил Колмогоров,
которых хотел писать <<Геометрию>> сам. Кроме того (как мы обсудим чуть ниже) и настроения в целом были совсем иные.
В итоге <<Геометрия>> Погорелова была издана в 1969-1970гг
как <<учебное 
 пособие>>, [Пог1969]. Именно эти учебники, несколько видоизмененные и еще недоработанные,
 были брошены на тушение пожара в 1982 году.
Погорелов сделал существенно новый учебник геометрии, но, по-моему, введение <<преобразований>> и векторов можно было бы  проделать
и в рамках традиционного подхода, использовав имевшиеся небольшие резервы%
\footnote{Строго говоря, это было не введение, а увеличение доз,
	потому что и векторы (1956-1966), и геометрические преобразования в программе
	были. Так как обучение в 4-5 классах переориентировалось в пользу геометрии,
	это давало выигрыш, который можно было бы разыграть, модифицируя уже существовавшие
	геометрические учебники, то есть <<Никитина>> или <<Киселева>>.}.

Все могло бы кончиться благополучно...
Но, как мы знаем, этого не случилось.

\sm


{\bf\punct Маркушевич.}
За Программой-1959, очевидно, стоял Алексей Иванович Маркушевич (1908-1979), который с 1958 года
был заместителем министра просвещения РСФСР. Остановимся подробнее на этой яркой и неоднозначной фигуре
(о нем см. статьи   [ГЛШ1958], [Шаб1968], [АЛМШ1978], [АГКЛШ1980], [Шаб1968],  [Чер1988],  [Тих2009]).

Учился в Ташкентском университете у В.~И.~Романовского, а потом в аспирантуре в Москве у М.~А.~Лаврентьева.
  В 34-35 заведовал кафедрой в Калининском педе. С 1935 года работал на мехмате МГУ. 
  В 1950 опубликовал  знаменитую в профессиональных кругах 
  книгу  (то ли учебник, то ли монография, то ли трактат) <<{\it Теория аналитических функций}>>
  (расширенное второе издание в 1967-68 году;
  английские издания 1965-1967, 1970),
  она и до сих пор  используется математиками
  (и в России, и за рубежом;  там подробно излагаются разные тонкие красивые вещи,
  вроде большой теоремы Пикара,
  лучей Жюлиа или соответствия граничных элементов
  по Каратеодори). Больше  активно математикой
  не занимался. Позже  издал  учебник по функциям комплексного переменного, 
  пользовавшийся популярностью во многих вузах,
  а также <<{\it Классическую теорию абелевых функций}>> (1979, английское издание 1992).
  
  Основной областью приложения сил Маркушевича после 1950 года
  была педагогика.
  
  В 1945 он становится членом-корреспондентом Академии педагогических наук (АПН) РСФСР (создана в 1943 году),
а в 1950 -- действительным членом  и вице-президентом.
Когда АПН РСФСР была преобразована в АПН СССР, он становится ее вице-президентом.
Пост вице-президента  он занимал в 1950--1958  и 1964—1975 годах, в промежутке  был
заместителем министра просвещения РСФСР%
\footnote{Министерство просвещения СССР было организовано в 1966г.}
в 1958—1964.

В 1959-1966гг. -  депутат Верховного совета РСФСР.

 В 1934—1937, 1943—1947 годах  заведовал редакцией математики в Издательстве технико-теоретической литературы.
Выступал инициатором издания серий книг <<{\it Библиотека учителя}>>,
<<{\it Популярные лекции по математике}>>. Участвовал в их издании.
Один из инициаторов   и один из редакторов первых двух изданий  <<{\it Детской энциклопедии}>> в 12 томах, 
1958-1962, 1964-1969,
и главный редактор третьего издания 1971-1978. Раздел по математике (<<{\it Числа и фигуры}>>) 
был очень завлекательным и понятно написанным.

 Разумеется, высокое качество популярной литературы по математике 1950-70гг во многом было заслугой
 авторов книг. Но это было в не меньшей степени и заслугой издателей, и думаю, что коснулось всех
 молодых людей, с интересом учившихся
 математике в те годы.

Кстати, он и сам написал несколько популярных тонких книжек для школьников и студентов:
 <<{\it Возвратные последовательности}>> (1950),
 <<{\it Замечательные кривые}>> (1952),
<<{\it Площади и логарифмы}>> (1952),
<<{\it Комплексные числа и конформные отображения}>> (1954),
<<{\it Целые функции}>> (1965),
<<{\it Замечательные синусы: введение в теорию эллиптических функций}>> (1965).

Эти книжечки  переводились на английский, испанский (в этом случае издательством <<Мир>>),
немецкий, польский, румынский, чешский, болгарский, турецкий  языки...

%


Коль скоро Маркушевич был главным педагогом-математиком Советского Союза, то вопрос о школьной программе
был в его прямой компетенции.

\sm

{\bf \punct Революционная программа.}
В 1957-1961гг. были изданы 6 сборников   
\href{http://www.mathnet.ru/php/archive.phtml?jrnid=mp&wshow=contents&option_lang=rus}{\it Математика, ее преподавание, приложения и история»},
(это довольно толстые книжки объемом 270-400стр.). Редакторами были Я.~С.~Дубнов, А.~А.~Ляпунов, А.~И.~Маркушевич, 
Там, в частности есть статья, с которой, по-моему, и начинается история колмогоровской Реформы, или,
может правильнее говорить, Революции:

\sm

	В. Г. Болтянский, Н. Я. Виленкин, И. М. Яглом,
\href{http://www.mathnet.ru/php/archive.phtml?wshow=paper&jrnid=mp&paperid=529&option_lang=rus}
{\it О содержании курса математики в средней школе}, Математика, ее преподавание, приложения и история, 4, 1959, 131–143

\sm

Напомню, что Наум Яковлевич Виленкин и Владимир Григорьевич Болтянский - замечательные математики.
Все три автора статьи много занимались популяризацией математики среди школьников,
в особенности Яглом и Болтянский.

\sm

В статье содержится  проект новой школьной программы. В качестве критики и обоснования там
говорилось много разумного (и полезного для дела). Но дальше предлагалось революционное  решение, причем решение глобальное,
охватывающее всю среднюю и старшую школу, арифметику, алгебру, геометрию и тригонометрию,
и требующее перестройки не только школьных программ,
но и всего стиля обучения. По сравнению с Программой 1959г. в качестве дополнительного выигрыша 
в школьную программу 
предлагалось ввести интеграл и метод координат (элементы аналитической геометрии).
А, например, комплексные числа (вместе с их тригонометрической формой) должны были быть в 8 классе.


Ни догмой, ни руководством к действию эта программа тогда еще не была, или, по крайней мере, таковой не объявлялась.
Из редакционного введения к сборнику:

\sm

{\it 
Эти   статьи} [имеется в виду еще статья В.~И.~Левина [Лев1959]]{\it 
не  следует   рассматривать   как   проекты   методических   и   
программных   документов,   как   не  следует   и  критиковать   их  за   то,   что   
в  них  не сказано,  не  учтено   и  не  подсчитано.  Это,   скорее,   «размышления  вслух»  научных   работников-математиков   о  том,
что  их  не  удовлетворяет    ныне   в  школьной   математике   и  какой  они   хотели   бы  ее   видеть. 
}

В том же сборнике был опубликован  критический комментарий к проектам 
за подписью И.~Н.~Бронштейна и А.~М.~Лопшица [БЛ1959].

Процитируем статью 
B. Г. Ашкинузе, В. И. Левина
и А. Д. Семушина  [АЛС1960] которые, по-видимому, входили в число авторов Программы-1959.
Их отношение к Революции,
 было выражено словами:

\sm

{\it 
Успешному   решению   задачи   во   многом   может   помочь  активное  участие   научной    математической    общественности.  
Насколько     интересные     
результаты   может   дать   такое   участие,   показывает   хотя  бы   программа,   
составленная   В.   Г.  Болтянским,   Н.   Я.  Виленкиным   и   И.  М.  Ягломом.  
В   ней   отражено   много   новых   прогрессивных  идей,  постепенное   освоение   
которых   нашей   школой   несомненно   внесло  бы   свежую  струю  в  обучение  
математике.   Однако   в   целом    эта   программа    настолько   далека   от   сложившейся 
системы    преподавания    математики,     что    осуществление     ее     
в    массовом    порядке    в   ближайшем     будущем    не    представляется    возможным,  
даже   если   бы   существовали   соответствующие   учебники.
}

\sm

Заметьте, авторы сомневаются не в возможностях школьников, не в проблемах человеческой психологии,
не в объективных сложностях, с которыми приходится сталкиваться, развивая у людей логическое мышление.
Они сомневаются (конечно, вполне справедливо) в способности современной системы преподавания  провести
программу в ближайшем будущем.

\sm

Вернемся к программе Болтянского--Виленкина--Яглома. Авторы предлагали крупный выигрыш за счет перестройки
и перетряски всей школьной математической программы за 4-10 классы (<<Арифметика>>, <<Алгебра>>,
<<Геометрия>>, <<Тригонометрия>>). Надо иметь в виду, что массовая математическая педагогика -- знание
экспериментальное. Хороший специалист может просчитать эффект от локальных возмущений того или иного курса.
Результат глобальной перестройки просчитать не может никто. Остается только надеяться на эксперимент.
Но тут возникает новая проблема. Так как перестройка глобальна, то учебники $n+1$ класса можно экспериментировать
лишь после того, как проэкспериментированы все учебники 4, 5, \dots, $n$-ого класса. Поэтому 
минимальное время эксперимента - 7 лет, после чего еще год нужен для анализа итогов. То есть при идеальных
результатах эксперимента запуск учебников в общую школу возможен лишь на 9-10 год. 

\sm

{\bf\punct Настроения 1965г.} Неудача с учебником Болтянского и Яглома случилась в 1964г. В тот год случились и иные
события, речь о которых пойдет ниже, а мы перенесемся в следующий, 1965г.

В качестве памятника тогдашним настроениям приведу несколько цитат из программной статьи 
Б.~В.~Гнеденко [Гнед1965]
в последнем номере журнале <<Математика в школе>> (жирный шрифт - мой, Ю.Н., по ходу дела я вставляю краткие комментарии,
но основные объяснения идут после приводимых отрывков),

\sm

{\it 
Прежде всего отметим, что ни в одном
школьном предмете нет и не может быть такого положения, чтобы он излагался почти
в точности так же, как пятьдесят, сто и
двести лет назад; чтобы в них ничего не сообщалось о результатах, тенденциях 
и возможностях современной науки. Ни в одном,
если не считать математики. По школьным
учебникам математики нельзя даже установить, что над миром 
пронеслась настоящая
научная революция и что математика вошла во все области знания и практической
деятельности в качестве мощного орудия
анализа и исследования. Из учебника невозможно узнать, что сама математика уже
не та, какой она была не только во времена
Евклида и Птоломея, но даже во времена К.Гаусса, М.В.Остроградского 
или же Б. Римана и П.Л.Чебышева. Современное содержание курса математики дает поразительно мало возможностей 
увязать его с удивительными успехами науки и с волнующими задачами наших
дней. Но даже то, что возможно, не делается. Нет спора, классическая математика
школьной программы является основой всей
современной математики и ее применений,
но красотами первых камней фундамента
трудно увлечься, если не увидеть хотя бы
контуров здания, которое из них может
быть выстроено...

\sm

{\bf Меня мучительно интересует вопрос о правильности
традиционной схемы урока}: повторение
пройденного на базе опроса учащихся (как
правило, наиболее слабых), объяснение нового материала,
беглый опрос или решение
задач у доски с целью выяснения степени усвоения 
вновь объясненных сведений.
Не 
слишком ли  много
тратится при этом
времени напрасно? Все ли это время используют
учащиеся для активной работы над
предметом урока? Приучает ли такая система учащегося к ответственности 
за порученное дело и не лишает ли она его какой-то
доли самостоятельности?

\sm

{\bf Сейчас все чаще и настойчивее раздаются
голоса о перегрузке учащихся, о невозможности расширения программы из-за недостатка времени.
Об этом много пишут и газеты общего направления, и специальные
педагогические журналы.} Мне кажется, что
во всех разговорах о перегрузке учащихся
имеется значительная доля недоразумений.
Я сам склонен говорить не о перегрузке,
а о непозволительном расточительстве времени, о крайне нерациональном его использовании. 
{\bf Речь может и должна идти о недопустимости перегрузки учащихся пассивным сидением,
когда не развиваются способности,
 притупляется внимание, теряется интерес к занятиям. Я убежден, что в школе производительно
используется
несравненно
меньше половины рабочего времени каждого учащегося....}

\sm

Не преуменьшаем ли мы развитие современного школьника,
 который еще
с дошкольных лет
 привык считать и читать
графики? {\bf Не следует ли пересмотреть
всерьез стиль школьного преподавания} в сторону
большей самостоятельности учащихся и большей напряженности их труда?
Учебник должен привлекать детей, заставлять их думать,
а не отпугивать и не усыплять...

\sm

... Вот почему важно
продумать курс математики так, чтобы его
изучение было интересно; содержание было современно, будило мысль и развивало
способности, а также открывало пути как в
научную, так и в практическую деятельность...

\sm

Мы подошли теперь к основной задаче настоящей статьи и, пожалуй, к основной задаче современной школы - каким должен
быть курс математики в школе? Я не буду
при этом касаться сравнительно малочисленных школ специально математического
или физико-математического направления,
а буду говорить о школах общего назначения, охватывающих подавляющее большинство детей школьного возраста...

\sm

....{\bf в наших
программах имеются и элементы аналитической геометрии, и основы теории пределов,
и даже элементы дифференциального
исчисления. Это так, но в ней нет элементов
интегрального исчисления, которые  могли
бы работать уже в самой средней школе}}  {\it и
многое упростить в изложении курса геометрии, а также физики.

\sm

{\bf Теория вероятностей в своих элементарных
понятиях очень проста и доступна}} [Ю.Н.: она однако ж опирается на элементарную комбинаторику
- предмет не столь уж простой для преподавания], {\it ее
легко увязать, с одной стороны, с традиционным материалом алгебры, а с другой с
тем, что сообщается в физике и химии, а
также в биологии. Кроме того, концепция
случайного так прочно вошла в современную практику, что без нее нельзя обходиться
ни в физике, ни в технике, ни в экономике,
ни в науках о живой природе. Необходимость
введения элементов теории вероятностей
в программы массовой средней школы
ощущается многими педагогами у нас, ощущается
она и в других странах. Недаром
в США, Японии и Югославии уже сделаны
попытки введения ее в программы средней
школы и созданы пробные учебники} [Ю.Н.: Берегитесь, а то нас опередят] {\it

\sm

{\bf Теория множеств} - одно из изящнейших
созданий математики XIX в.- легла в основу
всей современной математики, а также
многих концепций физики и техники. Немногие сведения, которые не займут большого
времени, могут раздвинуть интеллектуальные
горизонты учащихся....

\sm

{\bf Математическая логика} за последние двадцать
лет совершила буквально триумфальный
марш, и не только внутри математики,
но и в технической практике....

\sm

Я убежден, что те {\bf первичные сведения о программировании}, о которых идет речь, о кажут
огромную услугу в последствии многим
воспитанникам школы} [Ю.Н.:это замечание было справедливым, но каждая новая добавка усиливала дальнейшую нерезиновость программы,
кстати тогда надо было и обеспечивать доступ к ЭВМ].....

\sm

{\it
Я не касаюсь содержания курса геометрии, в последнее время у нас много и интересно об этом говорилось. В конечном счете
{\bf идеи, выдвинутые в свое время Ф. Клейном}} [Ю.Н.: См.[Бор1914]]. Именно эта <<коренная перестройка>> с опорой на геометрические
преобразования окажется основным
провальным элементом проекта], {\it сохранили свое значение и теперь. Переход
к этим идеям потребует {\bf коренной перестройки
сложившегося курса геометрии},
и этот переход следует совершить....

\sm

{\bf 
От педагога требуется,
чтобы он сумел увлечь своим предметом, и
работа учащегося проходила бы не под постоянным давлением принуждения, а под
знаком увлеченности, личной заинтересованности. К сожалению, этого добиваются далеко не все учителя и по отношению далеко
не ко всем ученикам....}}

\sm

Это уже конец 1965 года, статья Гнеденко - прокламация в поддержку начинающейся реформы Колмогорова,
но почти все сказанное
многократно  было бы сказано до того, я  цитирую именно ее лишь потому, что она удобна для цитирования.

Новым является то, что повышенные ставки теперь связаны не только с глобальной перетряской 
всей программы (об опасности чего уже говорилось выше), но и с переходом на новый стиль
обучения в школе и на всеобщее увлечение школьников математикой.
Такие надежды совсем удивительны:  даже в отобранных классах и хороших вузах всеобщая увлеченность не 
достигается, а в мат.классах что-то не наблюдалось особо творческого отношения к другим предметам, поэтому естественно 
было и не ожидать повального увлечения математикой в массовой школе.

\sm

Кстати, даже из приведенных цитат видно, что голоса против Реформы (о ней тогда уже говорили и писали) были.
С какими-то неизвестными нам авторами спорят и Колмогоров с Ягломом в [КЯ1965].

\br


\section{Дьедонне и Лихнерович}

\COUNTERS

Колмогоровская реформа была частью международного движения, математики изрядно начудили в ряде стран,
во Франции реформу разрабатывала комиссия под руководством  Лихнеровича (одного из Бурбаков).
Исследование тогдашних педагогических взглядов на Западе - отдельная работа,  автор за нее не брался,
читатель может найти много интересного в
[Фей1965], [Кла1973],  [Маш2006], [Фил2014]. Влияние интернациональных реформаторов 
на советских математиков и педагогов не вызывает сомнений (с другой стороны, это было время Sputnik crisis на Западе, мир наш взаимосвязан, быть может и наши реформаторы поспособствовали
заграничным движениям).

\sm

{\bf\punct Бурбаки.}
Когда мы со сверстниками, будучи студентами мехмата, обсуждали реформу учебников, то   замечали, что реформа близка по духу к книге 
[Дье1975]

\sm 

Дьедонне
{\it Линейная алгебра и элементарная геометрия}. 

\sm 

Она была издана по-русски в 1975 году  (перевод под редакцией  И.~М.~Яглома), французские издания 1964, 1968, 1968. 
Книга содержала жесткую  высокохудожественную критику традиционной школьной геометрии,  отчасти справедливую (а вполне несправедливую 
критику никто слушать не будет), взамен предлагалась программа сомнительной осуществимости. 

Дьедонне был членом влиятельной группы крупных французских математиков, выступавших под псевдонимом <<Никола Бурбаки>>.
В собственно математике они известны многотомным трактатом <<Элементы математики>> - труд, обладающий многими достоинствами
и в целом полезный (но понятность изложения не входит в число его свойств, и очевидно, что ее достижение не входило в цели
авторов, отдельные представители Бурбаков, например, Серр или Диксмье, писали свои собственные книги достаточно понятно). Эта группа также 
имела взгляды на школьное математическое образование.
Процитируем книгу
[Фил2014] Phillips, Christopher J. {\it The New Math: A Political History}, 2014

\sm 

{\it
Упор Бурбаки  на упразднение устаревших методов также повлиял на реформаторов
школьной
программы.
В речи 1959 года один из самых резких членов группы Бурбаки Жан Дьедонне
провозгласил <<Евклид должен уйти>>. Для Дьедонне евклидова геометрия 
была <<мертвым грузом>> с незначительным числом интересных вопросов 
и малой применимостью (re\-le\-van\-ce).
Полезным кусочкам можно благополучно обучить за несколько часов,
остальное имеет <<такое же отношение к тому, чем математики (чистые и прикладные)
занимаются сейчас, как магические квадраты и шахматные задачи>>. 
Математики
могли бы найти лучшие способы для обучения строгим логическим и аналитическим методам.
Дьедонне теоретизировал, что только вера и традиции удерживают Евклида в 
учебных программах. Хотя нападение концентрировалось на евклидовой геометрии,
оно в итоге относилось ко всем устаревшим предметам.

Речь Дьедонне вызвала острые споры. Он выступал на Royaumont Con\-fe\-ren\-ce --
на математическом форуме  по проблемам школьных
учебных программ.
Даже среди реформаторов мало кто сочувствовал почти полному изгнанию евклидовой геометрии
из школьных программ.
Однако провокационные речи Дьедонне оказали значительное влияние на американских
слушателей} [Albert] {\it Tucker..., Howard Fehr,...,} [Edward] {\it Begle, Robert Rourke, and Marshall
Stone.}

\sm

Упомянутые американские математики (среди которых такие знаменитости как Маршалл Стоун и Альберт Таккер) запустили реформу
New Math в Соединенных штатах.

\sm 

Во Франции была запущена реформа, 
возглавленная комиссией Лихнеровича.
Вот образцы сентенций из французского школьного учебника математики, Hachette, 
цитируется по книге Mashaal {\it Bourbaki: Secret society of mathematicians},
[Маш2006]

\sm

{\sc Определение.} {\it 
Подмножество $I$ упорядоченного множества $E$ называется интервалом $E$,
если оно удовлетворяет импликации} 
$$
(
x\in I 
\text{\,\,\,and\,\,\, } 
y\in I
\text{\,\,\,and\,\,\, }
x\le z\le y 
)
\Rightarrow z\in I
$$ 

{\sc Теорема.} {\it 
Тождественное отображение является единственным автоморфизмом поля вещественных чисел.

\sm 

Сложение и умножение функций определяет на множестве $L(V)$ эндоморфизмов векторного пространства
структуру унитарного кольца} [кольца с единицей] 

\sm 

 [Дальше привожу английский перевод с французского, потому что перевести на русский это уже невозможно. Но смысл примерно понятен] 
 
{\sc Definition.} {\it 
We say that an orthogonal endomorphism of $\phi$ of $E_3$ is a vector rotation to express the fact that the subspace
of variables invariant under $\phi$ has dimension 1 or 3.
} [Пытаюсь перевести: Движение пространства, оставляющее 0 на месте, является  вращением, если подпространство его
неподвижных точек имеет размерность 1 или 3].

\sm 

 В связи с французской реформой Понтрягин [Понт1980] цитировал слова, сказанные 
великим математиком Жаном Лерэ
в 1976 году:

\sm 

{\it 
Развитие понятия множества в последнее время значительно расширило область применения и силу математических методов, 
но значит ли это, что преподавание математики юношам и девушкам должно быть основано на этом понятии,
то есть проходить по схеме, принятой в прекрасном трактате Н. Бурбаки? Ответ может быть только отрицательным... 
Можно ли строить курс математики для юношества логически на теории множеств, то есть выразить сущность
этой теории на простом и доступном языке? 
Во Франции это пытались сделать с самонадеянностью, основанной на непонимании, что не могло не привести к катастрофе...

Торжество методики, основанной на повторении многословных определений, имеет самые серьезные социальные последствия. 
С одной стороны, это отваживает от научного образования способных юношей, которые лишены привилегии иметь взрослого руководителя,
способного объяснить им, что они правы, не понимая того, что им преподают, с другой стороны, это привлекает к занятиям
как раз наименее способных и думающих учеников, которые учат наизусть и повторяют, не понимая смысла...

Извращенная ситуация, в которой оказалось преподавание математических дисциплин во Франции, в большей степени,
чем в англо-саксонских странах, возникла из вполне законного стремления к прогрессу. Наши самые искренние 
и цельные реформаторы не сумели отстранить от этого дела шарлатанов, которые использовали их инициативу, например, тех,
кто с легкостью написал толстые учебники, полные ошибок, и получил преимущественное право на их переиздание, то есть воспроизведение ошибок.
Сами учителя были подготовлены интенсивной пропагандой... Методисты боятся потерять авторитет, если исправят допущенные ошибки.

Я прочел двум, сменившим один другого, министрам национального образования Франции основное содержание министерских инструкций, 
имеющих целью ошеломить наших детей научными определениями прямой... Они признали, что не понимают сами того,
что предлагают в качестве обязательных инструкций, однако инструкций не отменили.}

\sm 

{\bf\punct Бурбаки и мы.}
Реформаторские голоса из-за границы доносились и в Советский
Союз. Маркушевич, как высокий советский педагогический работник,
много ездил за границу и участвовал в разных конгрессах, 
посвященным проблемам образования. В своих русских статьях
он подробно отчитывался об услышанных им идеях (порой это сопровождалось вполне здравой критикой), см., например, [Мар1957]
и программную статью [Мар1964]. Ниже в списке литературы
эти статьи почти не отражены, но библиография Маркушевича легко доступна [АЛМШ1978], [ГЛШ1958].
Ограничимся одной цитатой. Из рассказа Маркушевича о Международной конференции по народному просвещению
в Женеве 1955 года  [Мар1957]:

\sm 

{\it 
 В связи
с математическим образованием в средней школе проф. Пьяже поставил вопрос о том, как некоторые общие идеи современной математики
(он сослался здесь на идеи Бурбаки) должны сказаться на построении
курса математики в средней школе. Он сообщил, что Международная
комиссия по математическому образованию, представляющая Международную математическую ассоциацию,
и Международная комиссия по изучению и усовершенствованию математического образования пришли
к выводу, что соответствующая реформа содержания математического
образования не только возможна, но и может облегчить обучение ма­тематике.}

\sm

Переводы бурбакистских проповедей разной степени радикализма издавались в СССР примерно с 1960 г, очевидно 
при чьем-то благожелательном к ним отношении.
В 1960 году в СССР была издана книга  {\it О преподавании математики} Ж. Пиаже, Э. Бет, Ж. Дьедонне,
А. Лихнерович, Г. Шоке, К. Гаттеньо (французский оригинал 1955, перевод А.~И.~Фетисова). 

Что касается официальной точки зрения советской педагогики
на колмогоровскую  реформу в момент ее начала, то она излагалась в статье [БМ1975]. Упоминались Международные конференции
по математическому образованию, конгрессы математиков в Амстердаме 1954 и Стокгольме 1962,

\sm 

Про Амстердам:

{\it 
Международная комиссия по математическому образованию представила доклад, в котором содержалось предложение радикально перестроить школьный курс
математики, {\bf положив в основу его понятия множества, преобразования и структуры}}. [слово <<структуры>> - визитная карточка Бурбаки]

\sm

Про Стокгольм:

{\it ... большинство стран предлагает введение в школьный курс элементарной теории множеств, элементов математической логики, понятий современной алгебры 
(группы, кольца, поля, векторы)%
\footnote{То есть Колмогоров был не первый, кто предложил ввести в школьную программу группы, кольца и поля, он просто следовал
лучшим зарубежным проектам.}, начальных сведений по теории вероятностей и математической статистике.}

Можно вспомнить яростную ругань В.~И.~Арнольда на Бурбаков (например, [Арн2000], [Арн2002]),
и вообще, и в связи с колмогоровской реформой, и в связи с образованием во Франции.
Цитируем Арнольда:

\sm 

{\it 
Мой учитель, Андрей Николаевич Колмогоров, очень меня убеждал, когда он начинал свою реформу, принять участие в этой реформе и переписывать все учебники,
делать их по-новому и излагать, как он хотел, бурбакизировать школьную математику и так далее. Я категорически отказался, прямо чуть не поссорился с ним, 
потому что, когда он мне стал рассказывать свою идею,
это был такой вздор, про который мне было совершенно очевидно, что пропускать его к школьникам нельзя. }

\sm 

Вот, что писал про Бурбаков яростный критик колмогоровской реформы Ю.~М.~Колягин [Коля2001],
его комментарии не содержат прямых ссылок, но как будто он выступает как очевидец.

\sm


{\it 
В 1966 г. очередное заседание Международного математического конгресса проходило в нашей стране.
Одна из секций конгресса была посвящена математическому образованию. В его работе официально участвовали и Н.~Бурбаки 
(пустое кресло с табличкой в зале). 
Вместе с профессором И.~К.~Андроновым я принимал участие в работе секции по математическому образованию. На секции речь шла о путях
и средствах коренной реформы школьного математического образования.

Выступавшие, в основном сторонники реформы, говорили о ней как о деле уже решенном в принципе, важном и нужном. Те трудности, которые
уже обнаружились на практике, объяснялись главным образом новизной подхода и неподготовленностью учителей. Следует заметить, 
что высшая школа оказалась в смысле реформы более консервативной и осторожной, чем средняя.

Подавляющее большинство отечественных математиков-педагогов и методистов (в том числе и автор данной книги) заразились этим новым «поветрием» с Запада. 
Никто тогда и не думал о том, какой урон нашей, отечественной средней школе нанесет эта реформа, как долго придется устранять ее последствия.}

\sm 

{\bf\punct Мы без Бурбаков.}
Но вопрос о западных источниках Реформы, быть может, вторичен, наши революционеры и реформаторы
всегда находят  на Западе то что им самим созвучно.

 Стремление изменить курсы математики в 1959 году было объективным и рациональным,
вопрос был в чувстве меры и способности просчитывать последствия. Насмотревшись за много лет на наш ученый и математический мир, автор думает,
что наличие мощной экстремистско-реформаторской
группы в 60е годы было неизбежным.
Проповеди Бурбаков и их  единомышленников  придавали нашим революционерам уверенность в своих силах, способствовали углублению их
радикализма, а также давали аргументы в споре с  оппонентами и в пробивании реформы в высших инстанциях. В этом смысле бурбакисты сыграли свою роль.

Автор думает, что  попытка реформы в духе проповедей обобщенных Бурбаков быстро остановилась бы при столкновении с действительностью, если бы не Маркушевич и Колмогоров.

\section{ Процесс пошел}

\COUNTERS

{\bf\punct Проект АПН СССР.}
В последнем номере <<Математики в школе>> за {\bf 1964 год} появляется программная статья Маркушевича
[Мар1964]. Алексей Иванович много говорит о взглядах Бурбаков и Пиаже, сохраняя при этом осторожный скептицизм
и не скрывая  возможных опасностей (понятно, что он был человеком умным и высококвалифицированным).
В той же статье сообщается, что {\bf АПН разработала проект реформы}
школьной математики. Описание не очень подробное (наверно, где-то в архивах еще лежат эти разработки),
но в целом это  близко к тому, что через 2-3 года явится на свет в качестве Колмогоровского проекта.

А дальше (см. [АГКЛШ80]):

\sm

{\it С  1965  по  1970  г.  он} [Маркушевич] {\it возглавлял  Комиссию  по  содержанию  образования  в  школах  СССР, 
  созданную  при  Президиумах  двух  академий  (АН  СССР  и  АПН), 
  в  составе  до    500  ученых  для  подготовки   проекта   
реформы  школьного   образования.}   

\sm

Надо сказать, что ученые тогда были в большой чести.  АН СССР едва ли могла вполне понимать
проблемы школьного образования. А вот АПН СССР это было положено по штату...

\sm

{\bf \punct Колмогоров вступает в игру.}
В какой-то момент (не позднее 1963 года) он становится председателем Комиссии по математическому образованию при
АН СССР (отметим, что этот предмет  шире, чем школьная математика, 
и Андрей Николаевич в этой должности выглядит естественно).

\sm

Первое документально известное выступление Колмогорова на тему школьных учебников состоялось 25 июня 1964 года
на совещании
в Минпросе РСФСР. Характерно, что основные докладчики - министр Е.И.Афанасенко и А.Н.Колмогоров.
Судя по реферату в  [Сов1964], выступление было реформистским, там есть ряд очень нехороших фраз, но  размеры  грядущих
потрясений из опубликованного текста все же не видны.

\sm

Следующее появление Колмогорова  [БМ1975]:
В декабре 1964 года была организована уже упомянутая Комиссия по  определению содержания среднего образования АН СССР и АПН СССР.
Математическую секцию комиссии возглавил А.Н.Колмогоров. Комиссия разработала новый учебный план,
по которому начальная школа
ограничивалась тремя годами (вместо четырех). Дальше предполагалось предметное преподавание.

\sm

Это было глобальное решение по всей школе,
в котором математики лишь принимали участие (впрочем, всю комиссию возглавлял Маркушевич). 
Интересно, что тогда было декларировано 
благое намерение о развитии факультативов в старших
классах,  не знаю, до какой степени его удалось реализовать. Но вообще факультативы в школах тогда в самом деле были.
Стоит иметь
в виду, что усложнение программ отбивает охоту к факультативам и у учеников, и у учителей.
Интересно, что когда факультативы 
явно отмирали, Колмогоров призывал к 140-часовому факультативному курсу математики в 9-10 классах 
(декабрь 1978 года,  [КС1978])%
\footnote{Буквально:
\it Факультативные занятия 
одно время получили довольно большое развитие, но потом интерес к ним 
стал угасать. Между тем последовательное проведение концепции 
общеобразовательного характера основного курса математики, избегающего не 
имеющих общеобразовательного значения технических осложнений, 
требует надлежащей организации занятий с теми учащимися, которые готовятся к работе в областях, тесно связанных с математическим аппаратом. В 
том числе следует самым широким образом удовлетворить и интересы 
учащихся, готовящихся к поступлению в высшее учебное заведение. 140-часовой дополнительный курс математики в IX-X классах мог бы содержать две большие дополнительные темы (комплексные числа с их применениями и комбинаторикой с началами теории вероятностей) примерно 20 
часов и углубленное изучение всего курса с решением более трудных задач (примерно на 100 часов). 
}.

\sm

Уже во втором номере (март-апрель) <<Математики в школе>> за 1965 год появляется материал для обсуждения
{\it Объем знаний по математике для восьмилетней школы} [Объ1965],
разработанный <<группой членов Комиссии
по математическому образованию математического отделения АН СССР
(И. М. Гельфанд, А. Н. Колмогоров,
А. И. Маркушевич, А. Д. Мышкис, Д. К. Фаддеев,
И. М. Яглом)>>. 
Предлагается уплотнить программу восьмилетки (с переносом туда
из старших классов показательной функции, логарифмов, 
плоскостей и прямых в пространстве и т.д)
и ввести некоторые элементы для школьной программы новые (векторы,...). 
По-моему, это бумага в поддержку
упомянутой чуть выше АПНовской программы, и это не удивительно, учитывая наличие среди авторов Маркушевича и Яглома.
Характерно, что про старшую школу в материале ничего не говорится,
но все это имело смысл прежде всего в отношении видов на старшую школу.
И главное - мы видим,
что Андрей Николаевич (который и был председателем Комиссии) в начале 1965 года 
уже присоединился к революционному проекту
(три другие фамилии  в дальнейшей истории реформы, кажется, не мелькали).
Издалека можно еще заметить, что эти лица, будучи знатоками профессионально-математического
образования, к массовой школе не имели отношения.

\sm

В том же номере журнала Колмогоров публикует оптимистичную статью [Колм1965]. 
Она выглядит чрезвычайно революционно (опора стереометрии на геометрические преобразования и векторы),
но, уже редактируя этот текст,
я с удивлением обнаружил,
что и ее, и статью [КЯ1965] во многом предвосхитил экспериментальное пособие Фетисова [Фет1963].
С этой статьи [Колм1965]
начинается поток публикаций Андрея Николаевича в журнале <<Математика в школе>>, см.
список в [Колм2003], больше 60 работ за годы реформы...

\sm

Колмогоров был человеком великих достоинств и великих заслуг, а тогда, в 1960е обладал огромным авторитетом.
Его имя стало знаменем реформаторов,
а он взял на себя ответственность за проведение реформы. Что касается деталей его участия, то они спорны.
Не вызывает сомнений, что Колмогоров  занимался учебником планиметрии. А
программа Реформы, видимо,  была в основном выработана до вступления
А.Н. в большую педагогическую  игру. Впрочем, Колмогоров  лично  поспособствовал радикализации программы. 

\sm

{\bf\punct Программа Колмогорова.}
Согласно [БМ1975], в 1965 году под руководством Колмогорова была разработана принципиально новая программа по математике.

\sm

Утверждается, что в 1966 она была опубликована (в качестве отдельной книги 
издательства <<Просвещение>>, я не нашел ни публикации, ни точной ссылки),
в 1967 появилась исправленная версия [Прог1967]. 

\sm

Последовал  комментарий
от общего собрания математиков АН СССР [Прог1967-м]. В ответе приветствовалось введение элементов мат.анализа в школе
(на самом деле в 1967 году эти <<элементы>> в школе были), отмечалась перегруженность программы
(указывались группы, кольца и поля), 
а также отмечалось, что для введения программы необходимо организовать переподготовку учителей.
Дабы улестить Колмогорова, одобрялось введение теории вероятностей в школу
(кроме трудностей с введением теории вероятностей, ей должна была предшествовать комбинаторика,
непростая для преподавания; комбинаторику  можно было бы и ввести, но тогда надо было многим жертвовать).
Отделение экономики АН СССР (Л.~В.~Канторович?) высказалось о программе положительно [Прог1967-э], но тоже отметило в программе странный
элемент в виде колец и полей. Редакция <<Математики в школе>> отметила, что замечания обоих отделений 
относились к предыдущей версии программы (1966), и они уже исправлены в 1967 году.
Итак, весьма сдержанные возражения были сняты, не успев быть опубликованными.

\sm

Так или иначе, программа  поддерживалась и АПН СССР и АН СССР. А.~Д.~Александров в 1978 году [Стен1978],
утверждал, что в 1966 году выступал против Колмогорова. Похоже (см. [Стен1978]),
что Общее собрание Отделения математики
АН СССР в 1966 году поддержало Колмогорова, не сознавая всего значения этого шага.

Были ли структуры, 
пытавшиеся оказать сопротивление, история на данный момент умалчивает.

\sm 

В трех номерах Математики в школе [О-прог1967-3], [О-прог1967-4],
[О-прог1968-1] были опубликованы подборки коротких статей учителей, методистов, работников педвузов, общий объем этих трех подборок
около 30 страниц (для обсуждения столь глобальной и судьбоносной  реформы не густо). Редакция журнала была к тому времени реформистской
и естественно предполагать, что она соответствующим образом фильтровала поток корреспонденции (в дальнейшем, во время кризиса реформы и начала контрреформации наличие решительной <<цензуры>> очевидно). Тем не менее несколько крайне
скептических статей было опубликовано. Несколько статей, начинавшихся с заздравиц,  имели весьма тревожное продолжение (и если бы у революционеров было бы желание подумать, даже эта подборка давала к этому достаточно поводов).

\sm 

В 1968г был  опубликован окончательный утвержденный Минпросом СССР вариант программы [Прог1968] 
(с комментариями Колмогорова).
Элементы теории вероятностей были изъяты 
 (<<с сожалением>> и <<ввиду неподготовленности нашей школы к их введению>>), комбинаторика оставлена.
Кажется, другие изменения были несущественны.

Составители \href{http://mat.univie.ac.at/~neretin/misc/reform/msh1967-1.pdf}
{{\it  программы}},
как явствует из публикации 1967 года  - это {\bf В.~Г.~Болтянский, А.~Н.~Колмогоров, Ю.~Н.~Макарычев,
	 А.~И.~Маркушевич,
Г.~Г.~Маслова, 
 К.~И.~Нешков, А.~Д.~Семушин, А.~И.~Фетисов, А.~А.~Шершевский, И.~М.~Яглом}
 (сведения об этих людях в последнем разделе статьи).
 
 \sm
 
 Окончательная редакция пояснительной записки - Колмогоров, Маркушевич (арифметика, алгебра и начала анализа), Яглом 
 (геометрия). 
 
 \sm
 
 Все идеи радужны и прекрасны, но при взгляде на  объем программы должны были встать древние вопросы, можно ли
 внедрить невнедряемое и впихнуть невпихуемое?
 И другой вопрос, который должны были задать уже
 вузовские преподаватели: а почему то, что должны будут понять все школьники, с такими трудами 
 и с такими затратами времени идет  в хороших 
 вузах со старательно отобранными студентами%
 \footnote{Я не видел таких голосов от математиков в печати. 
 	Учителя же недоумевали, почему та или иная тема, которая так тяжело идет в $n$-ом классе
 должна с легкостью преподаться в классе номер $k<n$?
Причем, если $n\ge 9$, а $k\le 8$, то это еще подразумевало разные множества школьников.}?...
 
 \sm
 
 Ничего хорошего этот проект не предвещал. В действительности вышло еще хуже.
  
  \sm
  
{\bf\punct Список учебников.}  Были пущены в дело и утверждены следующие учебники (списки их авторов варьировались,
в зависимости от издания,
содержание тоже):

\sm
 
 <<{\it Математика}>>, 4-5 класс, Н.~Я.~Виленкин, К.~И.~Нешков, ~С.~И.~Шварцбурд, А.~С.~Чесноков,
  А.~Д.~Семушин, Т.~Ф.~Нечаева
 — Под ред. А.И. Маркушевича, см. [ВНШСЧ]. 

 \sm
 
 <<{\it Геометрия}>>, 6-8 класс, 
 Колмогоров~А.~Н., 
 Р.~С.~Черкасов, А.~Ф.~Семен\'ович, Ф.~Ф.~Нагибин, [ А.~В.~Гусев], см. [КСНЧ]
 
 \sm
 
<<{\it Геометрия.}>> 9 - 10 классы
 В.~М.~Клопский,  З.~А.~Скопец,  М.~И.~Ягодовский (под редакцией З.~А.~Скопеца), см. [КСЯ1977]

\sm

 <<{\it Алгебра.}>>  Ю.~Н.~Макарычев,  Н.~Г.~Миндюк,  К.~С.~Муравин, [С.~Б.~Суворова] (под редакцией Маркушевича),
 см. [МММ]
 
 \sm
 
 <<{\it Алгебра и начала анализа}>>.
Был  учебник Б.~Е.~Вейц,  И.~Т.~Демидов
 9 класс (под редакцией Колмогорова)
Потом список авторов этого учебника уже за 9ый и 10 классы
многократно менялся, добавились А.~Н.~Колмогоров,  О.~С.~Ивашёв-Мусатов,  Б.~М.~Ивлев,  С.~И.~Шварцбурд,
уже на излете реформы в 1980г. добавился А.~М.~Абрамов, а Демидов выпал, см. [КВДШ].

\sm

Не берусь обсуждать учебник за 4-5 классы. Из остальных наиболее удачным (или наименее неудачным)
был учебник <<{\it Алгебры}>> под редакцией Маркушевича.
Наиболее неудачным был учебник стереометрии под редакцией З.~А.~Скопеца, об этом чуть ниже.

\sm

По всей видимости, коллективы авторов итогового комплекта учебников составлялись Колмогоровым и Маркушевичем лично,
и сделано это было до начала экспериментов, см [Вер2012].

\sm

{\bf\punct Программа-1968. Содержательное обсуждение.}

A){\sc Объем программы.}
Понятно, что программа была до невозможности перегружена.
 Форсированная гонка вперед
 вызвала не только перегрузку, она
 повлекла игнорирование учебных действий, которые были не нужны с логической точки зрения, но необходимы
 по психологическим причинам. 
 Мне не  хочется обсуждать это подробнее.

Думаю, что составители программы в действительности предполагали опустить
часть материала, если учебники <<не пойдут>>, см. [Колм1965-1]. О том, почему это не удалось, пойдет речь ниже. 

Позволю лишь привести цитату [Колм1967], которая показывает степень начального
энтузиазма (не знаю, была ли попытка реализовать именно
эту идею).

\sm

{\it 
Для IX-X классов при изложении стереометрии наиболее правильной кажется последовательно-векторная точка зрения.
В VII-VIII классах на уроках математики и физики учащиеся привыкнут к обращению с векторами (на практике, по преимуществу
лежащими в одной плоскости). Это позволяет в начале курса IX класса явно сформулировать аксиомы векторного пространства,
пригодные в любом числе измерений, обратить внимание на одномерный случай скалярных величин и сформулировать аксиомы, выделяющие двумерный
и трехмерный случаи.}

\sm

B) {\sc Теоретико-множественная идеология в учебниках.}
Цитирую Понтрягина [Понт1980]

\sm

{\it
Нет ничего предосудительного в том, чтобы в средней школе употреблялось <<множество>> как слово русского языка.
Так, определение окружности можно дать в двух вариантах.

Первый: <<Окружность состоит из всех точек плоскости, отстоящих от заданной точки на одном и том же расстоянии>>.

Второй: <<Окружность есть множество всех точек, находящихся на заданном расстоянии от заданной точки>>.

Второй вариант определения окружности ничем не хуже и не лучше первого%
\footnote{В старых учебниках писали <<геометрическое место точек>>.}.
И слово <<множество>> совершенно безвредно, а, в общем, бесполезно. Но в модернизированных учебниках и программах оно 
возведено в ранг научного термина, и это повлекло за собой уже серьезные последствия.}

\sm

Да, эта теоретико-множественная идеология была вещью более серьезной, чем просто словоупотребление
в стиле <<два притопа, три прихлопа>>. 
В \S1 уже приводились определения прямой и вектора.

 А вот я открываю учебник
<<Алгебры>> [МММ] за 6 класс, 1974. В его начале авторы вводят  соответствия
(бинарные отношения) между множествами. Они объясняют это на конечных множествах, а потом через
это определяют функцию. Должен ли 6-классник это понимать? Должен ли 6-классник воспринимать
кружочки со стрелочками как объекты однотипные с функциями вещественного аргумента?%
\footnote{А, скажем, использование в [МММ] скобок наизнанку, то есть $]1,2[$ для интервала от 1 и до 2,
из теоретико-множественной идеологии никак не вытекает.}

Другой пример из учебника под редакцией З.~А.~Скопеца [КСЯ1977], 6-ое издание, 1977г.

\sm

{\it
Преобразование $f$} [пространства] {\it можно рассматривать как множество всех
упорядоченных пар точек $(M,M_1)$, где $M_1= f (M)$. Множество
всех упорядоченных пар точек $(M,M_1)$ определяет другое 
преобразование, отображающее $M_1$ на $M$. Такое преобразование 
называется обратным преобразованием к $f$ и обозначается $f^{-1}$.
Если $M_1=f(M)$ то $f^{-1}(M)=M_1$.}

\sm

Да уж, понятней не скажешь... 

\sm

C) {\sc Перестройка курса геометрии.}
<<Теоретико-множественная идеология>> вызывала наибольшее раздражение в 70е годы,
но основную роковую роль сыграла не она.

Геометрия традиционно излагалась в духе Евклида. Реформаторы предполагали
опереть курс геометрии на геометрические преобразования и на векторы,
это была идея-фикс, существовавшая задолго до вступления Андрея Николаевича
в игру. То, что это  легко реализуется
как абстрактная логическая конструкция вполне очевидно, и каждый достаточно грамотный математик
может проделать такое упражнение. Но из этого никак не вытекает, что школьники поймут подобный курс.

Приведем длинную цитату из статьи [Колм1967] <<{\it Новые программы и некоторые основные вопросы усовершенствования
курса математики в средней школе}>>:

\sm

{\it 
Если я назвал действующие сейчас программы архаичными, то особенно это относится к геометрии. Составители проекта программ находились здесь
в затруднительном положении, так как по существу следовало бы начинать не с составления новых программ, а с работы по
выяснению желательной логической структуры курса геометрии. Ориентироваться здесь приходится по преимуществу на иностранный опыт,
а на отечественную экспериментальную работы в школах -- лишь в применении к младшим классам.

Основные тенденции перестройки школьного курса геометрии, уже нашедшие самое широкое признание, можно сформулировать
в виде трех положений.

1) Формирование начальных геометрических представлений происходит в младших классах.

2) Логическая структура систематического курса геометрии в средних классах заметно упрощается 
по сравнению с евклидовой традицией. Развитие привычки к строгим логическим доказательствам
на этом этапе соединяется с открытым признанием права принимать без доказательства избыточную
систему допущений.

3) Курс геометрии в старших классах строится на основе векторных представлений. При этом естественно
и обращение к координатному методу (однако в качестве вспомогательного, так как изложение не делается от этого
обращения менее <<геометрическим>>).

...................................

Проект программ предлагает начинать систематический курс геометрии в VI классе. Этот курс геометрии
задуман не в духе евклидовой традиции, а в соответствии с вторым из сформулированных выше положений.
По-моему мнению, хорошим образцом осуществления этого замысла остается французский учебник Эмиля Бореля,
написанный еще в 1905
году.}

\sm

Проблемны здесь п.2 и п.3.
В качестве точки опоры указывается учебник Бореля, но в воинственной книге Дьедонне 1964г
[Дье1975] евклидов подход к геометрии обличается с тех же позиций, что и у Колмогорова,
и с теми же призывами, что у Колмогорова.
Из чего вытекает, что французская школа по  учебнику Бореля не училась. Вначале говорится про иностранный
опыт, но он вроде бы относился к <<выяснению желательной структуры>>. Реализуемость колмогоровской декларации
на абстрактно-логическом уровне была несомненна, из чего никак не вытекало, что  она
пройдет  в реальной жизни.
Что касается, <<самого широкого признания>>, то учебники еще не были написаны. Нельзя даже сказать, что
широко признан был <<кот в мешке>>,  мешка с котом до 1970г. не было.

\sm

 Думаю, что именно это оказалось основной причиной неудачи учебника Колмогорова, Семеновича, Нагибина и Черкасова --
 краеугольного камня 
 всей реформационной программы. Несмотря на все таланты авторов и несмотря на все их усилия,
 школьники этот учебник не понимали.
 Он не поимел  упрощенных наследников, хотя 
 спрос на учебники геометрии  разного уровня  был, да и колмогоровское <<лобби>> оставалось мощным.
 Доселе используются
 не менее древний <<Погорелов>> [Пог1969], [Пог1982], и чуть более юные <<Атанасян>> [Ата1981]  и <<А.~Д.~Александров>> [АВЛ].
 Кстати, упрощенных наследников имели три учебника Колмогоровского проекта -- <<Математика>> 4-5,
 <<Алгебра>> 6-9, <<Алгебра и начала анализа>> 9-10 -- их удалось откатить на приемлемый уровень сложности.
 Учебник Геометрии был неоткатываемым.
 
 Хотя с <<Планиметрией>> все вышло совсем плохо, со <<Стереометрий>> получилось еще хуже.
 Во-первых, основная идея перестройки для Стереометрии была еще более порочной, во-вторых,
 авторы (вообще-то люди известные в педагогических кругах)
 по уровню своих талантов далеко не дотягивали до уровня авторов планиметрического учебника.
 De facto получился учебник, который Колмогоров
 в 1978 назовет <<определенно неудачным>>.
 
 \sm
 
 D) {\sc О математике, как об отчасти гуманитарном предмете.} В  советской школе
 геометрия традиционно была основным инструментом развития логического мышления у учеников.
 Стереометрия   была также важным  инструментом развития пространственного воображения%
 \footnote{Это рассматривалось  как одна  из целей предмета геометрии, см. [Чет1950].}.
 
 Процитируем записку Колмогорова к Программе-1968 [Прог1968]:
 
 \sm
 
 {\it
 Мы считаем, что практика преподавания геометрии в VI-VIII классах в настоящее время слишком
 часто направлена лишь на создание иллюзорной <<строгости>>. Преподаватели математики
 из курса педагогических институтов  хорошо знают, что все научные системы изложения геометрии на основе
 аксиом сложны, а список употребляемых при этом аксиом длинен. Но в школьной практике
 укоренился обычай указывать лишь <<примеры аксиом>>. Самый список этих примеров обычно до смешного короток.
 По-видимому, учащимся ни разу не предлагают провести анализ какого-либо доказательства 
 с выявлением всех лежащих в его основе аксиом. Между тем такое упражнение следовало бы настойчиво 
 рекомендовать...}
 
 \sm
 
  Здесь, как это ни печально, Колмогоров демонстрирует фундаментальное непонимание.
 Дореформенные учебники геометрии, особенно в начальной их части, не стремились к логической строгости:
 параллельно с обсуждением начальных утверждений они на наглядном и не допускающих
 различных интерпретаций материале учили человека логическому мышлению и показывали полезность
 логики... Ученики, способные понимать строгие логические рассуждения не сваливались с неба, их исподволь готовили
 на рассуждениях имитированной строгости. Согласно предложению Колмогорова, из обучения математике предполагалось 
 вынуть фундаментальную составляющую.
 
 \sm
 
 По поводу того, что получилось, процитируем А.Л.Вернера [Вер2012], соавтора А.Д.Александрова  по учебнику геометрии (я сам не считал...)
 
 \sm
 
{\it В Колмогоровской «Геометрии, 6–8»  в шестом классе \dots: из выделенных там 38 утверждений 
 почти две трети оставлены без доказательства, да
ещё среди недоказанных утверждений много и не выделено. Учебник первого
года систематического курса геометрии, который всегда двигался равномерно
и последовательно на уровне строгости, доступной ученикам этого возраста,
стал прерывистым по своему содержанию: то что-то докажем, то примем без
доказательства, затем снова что-то докажем, а затем снова примем без доказательства и т.д. 
\dots Того,
чего добивался А.Н. Колмогоров от своего курса геометрии — повышения его
уровня строгости (в том смысле, который вкладывал в слово строгость А.Н.
Колмогоров) и одновременного упрощения курса геометрии — у него не получилось.}

\sm
 
Стоит иметь в виду, что вынутая из курса <<Планиметрии>> составляющая
была не только необходима для обучения математике, прежде всего  на ней основывалось
гуманитарное значение математики. Процитируем слова  А.Д.Александрова
по  поводу непонятной школьной математики:
 
 \sm
 
 {\it 
 Вряд ли есть что-либо более вредное для духовного — умственного и морального — развития, 
 чем приучать человека произносить слова, смысл которых он толком не понимает
 и при необходимости руководствуется другими понятиями.}
 
 \sm
 
{\bf \punct График внедрения.}
Как я замечал выше, безопасное введение революционной школьной программы 
могло быть возможным не раньше, чем через 8-9 лет после начала широкомасштабных
экспериментов.
 В 1968 году большая часть революционных  учебников за 4-10 классы еще не была написаны,
 в 1968г. начались эксперименты по <<Математике>> 4 класса,
 а  всеобщий переход на новую программу начался в 
уже 1970 году (см. [К-на-1970]).
Осенью 1972 года до массовой школы дошли новые <<Геометрия>> и <<Алгебра>>, а 1975 году
новые <<Стереометрия>> и <<Алгебра и начала анализа>>.
Первый массовый  выпуск по новой программе тем  самым состоялся в 1977 году.

\section{ Успешный эксперимент}

\COUNTERS

{\bf\punct О времени эксперимента и экспериментировавшихся учебниках.%
\label{ss:experiment}}
В декабре 1978 учебники Колмогоровского проекта стали предметом обсуждения общего собрания Отделения
математики АН СССР,
 [Стен1978],
 
 \sm

КОРОТКОВ (зам министра просвещения СССР). {\it Достаточно сказать, что мы работали над учебниками 4-10 классов. Каждый учебник
рассматривался в трех-четырех вариантах 
со стороны разных групп авторского коллектива; затем в течение трех лет
экспериментировался  
 причем довольно широко, так как 7,5 тысяч школьников каждого класса учились по этим учебникам прежде, 
чем через три года он вводился как учебное пособие.}

\sm

Сказанное выглядит несколько загадочным. Дело в том, что ни подробнейшие списки  работ Колмогорова [Колм2003],
[Абр2016],
ни библиотечные каталоги не видят <<Геометрии>> под редакцией Колмогорова ранее 1970 года, а в массовую школу он пошел в 1972 году.
 До 1970г. не видно и <<Алгебры>> под редакцией Маркушевича, а до 1968 не видно учебника <<Математика-IV>>
Виленкина и др. В публикации [К-на-1970], посвященной началу введения новой программы с сентября 1970, речь шла об экспериментах  в 4 и 5 классах, но не в 6-ом.

\sm

В общем не было 3 лет экспериментов. Их было лишь два.

\sm

 Следующий вопрос, о каких 3-4 вариантах учебников идет речь?
 
 \sm

{\footnotesize
1) {\sf Учебник математики для IV классов.}
Здесь такие варианты находятся,
они упоминаются в [В-мин-1969], [К-но-1970]:

\sm 

И.~Г.~Баранова, З.~Г.~Борчугова, {\it Математика: 4 класс} : - Москва : Просвещение, 1968. - 398 с.

Н.~Я.Виленкин, А.~И.~Нешков, С.~И.~Шварцбурд,
{\it Математика: 4-й класс}, под редакцией А.И.Маркушевича,
 Москва : Просвещение, 1968. - 304 с.

 С.~А.~Пономарев, П.~В.~Стратилатов, Н.~И.~Сырнев,
{\it  Математика: 4-й класс}, под редакцией В.~И.~Левина,
  Москва : Просвещение, 1968. - 272 с.
  
   Н.~А.~Принцев,  Л.~Н.~Принцева,   М.~И.~Ягодовский {\it 
  	Математика 4 класс. Пробный учебник,} Просвещение,  1968  
  
  \sm 
  
 У всех этих%
 \footnote{Существовал также экспериментальный 
 учебник    И.~К.~Андронов,  Ю.~М.~Колягин,  Е.~Л.~Мокрушин,  Е.~С.~Беляева {\it Математика (множества, числа, фигуры, операции)}. 
 Новое экспериментальное учебное пособие для
 4 класса общеобразовательной школы, Просвещение,  1969, не знаю, ставились ли эксперименты по нему. Андронов руководил разработкой реформы младших классов.}
  пробных учебников были изданы продолжения за 5-ый класс,
 однако свой окончательный вердикт Коллегия Минпроса 
 СССР вынесла 25.07.1969 по итогам первого года экспериментирования.
 Победителем стал учебник под редакцией Маркушевича.
 
 \sm 
 
2. {\sf Систематические курсы планиметрии
 	и алгебры за 6-8 классы.}
  Было известно, что планиметрией занимается лично Колмогоров,
 а редактором учебника алгебры был лично Маркушевич. Писание учебника
 -- сложная и объемная работа. Если поставить себя на место тогдашних
 возможных авторов, то желающих
  писать учебник с гарантией того, что он не пойдет 
 в жизнь, должно было быть не много. Быть может, я искал сведения
 о конкурентах недостаточно настойчиво.
 Но я  смотрел в самые разные источники, 
 где они должны были бы упоминаться,
  и упоминаний о таковых 
 за 1968-1971гг. не нашел. В 1972г. (то есть, когда решение уже было принято и учебник
 Колмогорова пошел в массовую школу) появился учебник [БВС1972]:
 
 \sm 
 
 В. Г. Болтянский, М. Б. Волович, А. Д. Семушин. 
 {\it Геометрия, Эксперим. учеб. пособие для VI кл}.  Москва : Педагогика, 1972. - 110 с.  
 
 \sm
 
 Вряд ли он более удачен, чем учебник Колмогорова.  Но будь он и таковым,
 это уже не имело бы значения.
 Подведу итог этой части обсуждения: {\bf автору 
 не удалось найти пробных учебников, которые конкурировали бы с <<Геометрией>> Колмогорова и <<Алгеброй>> под редакцией Маркушевича.}
 
 \sm 
 
 {\sf 3.  <<Алгебра и начала анализа>> и <<Стереометрия>> 
 за 9-10 классы.} 
 По <<Алгебре и началам анализа>> Е.~С.~Кочетков и Е.~С.~Кочеткова
 (авторы действовавшего тогда учебника) выпустили пробный учебник, 1969, 1971.
 Однако на учебнике Вейца и Демидова [ВД1969] сразу появились слова  <<под редакцией
 Колмогорова>>, что уже означало предпочтение.
 
 Приведем цитату из академической стенограммы 1978 [Стен] (судя по контексту, описываемый эпизод относятся к 1974-75г.):
 
 \sm
 
  ШАБУНИН.
 {\it Несколько лет тому назад, когда
 	впервые готовился учебник для IX класса, мне, как частному лицу,
 	был послан на рецензию учебник для IX класса по Алгебре Вейца и Демидова.
 	Было заседание в Академии педагогических наук, где присутствовали Андрей Николаевич
 	и Алексей Иванович Маркушевич. Уже все было подготовлено к тому,
 	чтобы учебник рекомендовать, потому что был конкурирующий учебник
 	(тот учебник хотели, видимо, забраковать%
 	\footnote{Cf. [Колм1967-1].}). По сравнению
 	с ним теперешний учебник - сказка. Я не понимаю, почему такой
 	учебник, прошедший несколько рецензий, рекомендовали!? Грубейших
 	ошибок я нашел не меньше двух десятков и свою рецензию передал Маркушевичу.
 	Он схватился за голову, когда прочел. Было принято решение
 	взять этот учебник за основу, добавить в авторский коллектив 2-3 человек.
 	В течение месяца его переработали, грубые ошибки убрали, и теперь oн
 	появился как учебник IX классов.}
 
 \sm 
 
 Учебник Вейца--Демидова должен был пойти в эксперимент в 1973г.
 На мой взгляд, положительный исход 
 этого эксперимента невероятен (автору не приходилось работать в регулярной школе, но в хорошем техническом вузе я преподавал 20 лет,
 на глазок примерить его на студентов  могу).
 Авторами учебника, который спустя 2 года
 пошел в дело, были уже А.~Н.~Колмогоров, Б.~Е.~Вейц, И.~Т.~Демидов, О.~С.~Ивашёв-Мусатов, С.~И.~Шварцбурд. То есть это не учебник,
 который успешно прошел эксперимент, а другой учебник с другими авторами.

По <<Стереометрии>> существовал пробный учебник
К.~С.~Барыбина.

По-моему, все эти четыре учебника, конкурировавшие за два предмета, были неудачны, что во многом было предопределено
изначально правилами игры%
\footnote{Например, учебник Погорелова  не мог даже рассматривался как конкурент.}.
 Решающую роль
в выборе, надо думать, сыграл не эксперимент, а точка зрения Колмогорова.
Группу из З.~А.~Скопеца, В.~М.~Клопского и М.~И.~Ягодовского составил (согласно 
рассказу А.~Л.~Вернера) лично Колмогоров, а их учебник 
был жестко прицеплен к учебнику планиметрии,
что давало авторам дополнительное преимущество. 

} 

Вот такие данные об экспериментировании с альтернативными учебниками.

\sm 

{\bf \punct Эксперимент в московской экспериментальной школе.}
Я учился по этим, тогда еще экспериментальным,   учебникам с 4ого класса с 1968 года
(то есть был в первом потоке обучаемых).
Ничего  не помню про 4-5 классы, а  <<Геометрию>> Колмогорова помню.

 Это была
школа  710 г.Москвы рядом с метро Студенческая. Она была <<экспериментальной школой>> АПН СССР 
(это были ее лучшие времена,  вскоре со сменой директора многое изменилось, вышло, что я вовремя ушел в 91 школу),
состав учителей, конечно, 
 был  лучше, чем в среднестатистической школе.
Наш класс  был частью из района (район был, как сейчас принято говорить, 
хороший), частью был набран (не знаю, каким способом). Отношения между двумя половинами класса были 
нормальными, внутри класса никому не приходило в голову различать тех и других
(те, кто ездил издалека, в среднем учились лучше), 
а оппозиция учителям и обучению, разумеется, была, но носила здоровый 
характер.

В общем, была школа с хорошим набором учеником и улучшенным набором учителей. Но это не вся правда.

У нас был факультатив по математике, который вела выдающийся педагог Михайл\'овская Ариадна Юрьевна. 
Мне не удалось найти о ней
каких-либо данных,
кажется (но  могу ошибаться) она имела титул <<заслуженного>>, 
когда-то была учительницей в нашей школе, потом ушла куда-то вверх
в структуры АПН. Факультатив она вела очень ярко и понятно (причем понятным это оставалось и потом), сильно  выходя 
за рамки программы. Факультатив пользовался популярностью, как по причине своих достоинств, так и потому, что тогдашняя
школа 710  вовсе не старалась угрузить учеников до предела обязательными предметами - учебная нагрузка 
была весьма сдержанной. Кажется, Михайловская провела у нас пару 
уроков математики, и люди на нее повалили. Фактически мы имели учебник плюс Михайловская, и этот набор был и более разумен
и более понятен, чем просто учебник (в частности, она нормально рассказала про множества). На факультатив, 
конечно, ходили не все, но если предмет воспринимает часть учеников, то
и остальные подтягиваются тусовочным образом.

Так или иначе, юридическая школьная программа, по существу,
была встроена для нас в более широкие рамки. По мне это было очень хорошо, но  эксперимент
(у нас ведь проводился эксперимент) переставал быть точным экспериментом, и на такие вещи экспериментаторы обязаны были бы сделать поправку.

Не знаю, чему  в итоге научились мои
одноклассники. Мне учебник Колмогорова был  интересным и полезным.
Но здесь есть и другие <<но>>. Я читал много другого, начиная с конца 4 класса,
когда мне попала в руки книга Сергея Боброва <<Волшебный двурог>>.
Я читал какие-то пожелтевшие учебники геометрии, среди них был Киселев,
но было и что-то более подробное и продвинутое, по видимому, это были какие-то пособия 
для учителей, которые я сейчас не могу опознать, там там было много дополнительного материала мелким шрифтом. Читал справочник Выгодского,
читал Погорелова. Могу  засвидетельствовать,  что научно-популярная и педагогическая литература по математике
тогда лежала в Москве в обычных книжных магазинах, и что качество массово издававшейся тогда научно-популярной литературы было высоким.
Программный  материал по евклидовой геометрии%
\footnote{Помню, что при прохождении гомотетии меня удивило название, а не понятие,
про саму ее я где-то читал. Задачи на гомотетию были интересными.
Понятие сжатия для меня было новым.}
я обычно понимал до его прохождения, и учебник Колмогорова для меня был, отчасти, дополнительным чтением.
Вполне можно было подивиться
на упомянутое выше определение прямой (но определение вектора все равно было зубодробительно). Повторюсь, учебник Колмогорова был
для меня полезным и интересным, но он  не играл той роли, которую
должен играть учебник, он был ПОСЛЕ восприятия начал евклидовой геометрии,
а не ДО. А в смысле дополнительного образования эта книга,
разумеется, хорошая и качественная. 

\sm 

Школа гордилась тем, что вводит новую программу и от души старалась сделать это как можно успешнее. Скорее всего, 
в Министерство и АПН был отправлен победный рапорт. Возможно, что победный рапорт был отправлен и в отношении
учебника <<Стереометрии>>, хотя спасти этот курс обобщенная Михайл\'овская могла бы лишь прочтя вперед
курс нормальный.

\sm

Первый услышанный мной скептический отзыв (1972) об этой программе исходил от 
учителя 91 школы г.Москвы Владимира Мироновича Сапожникова...

\sm 

Из академической стенограммы [Стен1978]:

\sm 

М.И.ШАБУНИН (представитель физтеха) 
{\it  Представители министерства и Академии педагогических наук мне могут возразить,
что эксперименты проводились.
Но эти эксперименты
носили локальный характер, сплошь и рядом они проводились людьми, которые состояли в штате или полуштате Академии 
педагогических наук, и едва ли только на их опыте можно было основывать учебники...}

\sm 

{\bf \punct От Тосно до Белоярского.%
\label{ss:tosno}}
Что все-таки известно из литературных
источников об  эксперименте 1968-1969 и 1969-1971гг. по 4-5 классам?
Согласно информации [В-мин1969] о коллегии Министерства просвещения СССР от 25.07.1969,

\sm 

{\it Проверка проводилась во всех школах трех районов Российской Федерации
	(Суздальском районе Владимирской области, Тосненском Ленинградской области, Белоярском
	Свердловской области), во всех четвертых классах Севастополя, в некоторых школах Москвы,
	Новосибирска, Куйбышева и др., в небольшом числе классов почти всех союзных республик...}

\sm 

То есть относительно широкий эксперимент вроде бы ставился (хотя я встречал высказывания,
что список выше не вполне соответствует действительности). Понятно, что рапорты в <<Математике
в школе>> об учебнике Виленкина и др. были оптимистичны. Например, [К-на-1970]

\sm

{\it  Как показывает первый опыт работы в IV-V классах, новые учебники оказывают положительное влияние
	на развитие учащихся...}

\sm

Однако там же:

\sm 

{\it
	Так, например, в результате переработки пробного учебника IV класса
	оказалось возможным несколько упростить программу по математике для этого класса,
	найти более экономные методические приемы изложения материала, отработать 
	более эффективную систему упражнений. В результате переработки первого варианта
	пробного учебника по математике для IV класса он был сокращен на 15-20 процентов.}

\sm  

Свидетельствует ли это о больших успехах или это эвфемизм, необходимый в передовице
<<Математики в школе>>?
Снова ссылаемся на [В-мин1969].
При переходе от 1968-69 к 1969-70 учебного году

\sm 

{\it
	Из пособия исключены такие вопросы <<Как раскрывать скобки?>>,
	<<Сокращение частного>>,
	<<Распределительный закон деления>>, <<Геометрия вокруг нас>>,
	<<Измерения величин>>, <<Правила и формулы>>. Сокращено количество упражнений
	с 1613 до 1368. После доработки объем учебного пособия сократился на 1.5 листа.}

\sm

В том же номере Математики в школе, где объявлялось о начале обучения в сентябре 1970г. по новой программе, на следующей странице сообщалось [Вним1970]

\sm

{\bf Вниманию учителей математики  четвертых классов и руководителей методических объединений}

\sm

О ВРЕМЕННОМ СОКРАЩЕНИИ УЧЕБНОГО МАТЕРИАЛА ДЛЯ IV КЛАССА

\sm 

{\it  Главное управление школ Министерства просвещения СССР сообщило, что впредь до
	перехода в четвертые классы учащихся, прошедшим обучения по новым программам начальной школы, в изучение учебного материала по математике вносятся следующие изменения
	
	1) Исключаются разделы Окружность и круг, ... Измерения на местности}.

\sm 

[Далее сообщается об исключении задач с указанными номерами, всего 230]

\sm 

{\it   Время, высвободившееся за счет названных сокращений, должно быть использовано
	для   более глубокого изучения программного материала по математике.}

\sm 

Итак после двух прогонок учебника Математика-IV ряд  тем передвинулся в
следующий класс, а число задач сократилось на 30 процентов. Разумеется, задачи сокращались и по темам, которые д\'олжно было углубленно изучать.

Очевидно, что  первые два года эксперимента прошли, по меньшей мере, не вполне благополучно%
\footnote{Надо заметить, что систематические курсы геометрии и алгебры предполагалось опереть на эту самую пропедевтику 4-5 классов. Фактическое сокращение программы по 4-5 классам таким образом подмывало землю под систематическими курсами.}. 

Выше я не просто так обсуждал, было ли два года эксперимента или три.  Разница была принципиальной.
Через два года после начала экспериментирования по пропедевтическому курсу
IV-V классов  одновременно 1 сентября 1970 года началось экспериментирование
по систематическим курсам геометрии и алгебры в VI классе, а также началось всеобщее внедрение новой программы в IV.

 Мосты были сожжены.

\section{Столкновение с действительностью}

\COUNTERS

{\bf \punct Сентябрь, 1972.}
Вслед за началом  успешного эксперимента в 1968 году (но вовсе не после эксперимента),
с 1970 года новые учебники начали с 4 класса вводиться  по всем школам
[К-на1970]. Курс <<Математика>>
за 4-5 класс был пропедевтическим, и, скорее всего,  неплохим (ранних версий этого учебника у меня нет, а мои воспоминания крайне смутны).
Но учебник этот не оправдал возлагавшихся на него чрезмерных пропедевтических надежд.

\sm

Когда волна докатилась до 6 класса (до <<Геометрии>>),  стало ясно, что школьники
новых учебников не понимают. В сентябре-октябре 1972 года школа испытала первый шок.  Кстати, в тот год по свидетельству Колягина [Кол2001],
была отменена годовая оценка по геометрии, я натыкался на свидетельства, что кое-где отменялись
и оценки по геометрии  за первую четверть.

\sm

{\bf \punct Битва в пути.}
Был выдвинут лозунг: <<{\it учителя не готовы к работе по новой программе}>>.
Это соответствовало истине, так же, как истине соответствовали слова, которые вслух произносились не всегда,
но всегда подразумевались
<<{\it Учителя (часто) и сами этих учебников не понимают}>>. Но  это скрывало другую, более печальную, истину, что новый комплект 
учебников никуда не годился ни при каких учителях и на реальных детей рассчитан не был%
\footnote{И еще одну истину, не столь значимую, но забавную: даже авторы учебников не всегда
понимали, что они написали.}.

В физмат-школах учебники работали, но фактически физмат-школы учились по расширенной программе. 
Тамошние учителя вставляли содержание учебника Колмогорова (как до того и после того содержание иных учебников)
в более широкие рамки (тот же эффект, что и в описанном выше эпизоде с А.~Ю.~Михайловской, но в еще более рафинированной форме).
За исключением небольшого круга школ (или даже классов), этот подход  уже не мог работать (или требовал фактического раздувания
учебной нагрузки). В каких-то случаях расширение рамок могли провести родители.

\sm 

Предпринимались отчаянные усилия по усовершенствованию учителей, их надо было обучить новой программе, а потом
методике преподавания по оной. Для этого надо было усовершенствовать институты усовершенствования учителей.

Учебники исправлялись, программа сокращалась...

Ничего не помогало. С каждым годом в зону реформы втягивалось все больше
школьников,  срок пребывания ранее втянутых школьников в этой странной образовательной зоне увеличивался. 
Авторы учебников имели в виду развитие творческих способностей школьников. В действительности,
получилось наоборот, школьная математика стала превращаться в шаманский набор заклинаний.

А тут подоспел и  \href{http://mat.univie.ac.at/~neretin/misc/reform/Skopets.djvu}
 {\it учебник 
геометрии для 9-10 классов} под редакцией З.~А.~Скопеца. В итоге в 1977 году был отменен выпускной школьный экзамен по геометрии.
Отменен навсегда. В декабре 1978 года Колмогоров назвал этот учебник <<определенно неудачным>> [КС2012],
но ведь он был <<определенно неудачным>>
 с самого начала... (а первая версия вышла еще в 1969 году).

Учебники  алгебры были неудачны, но далеко не до такой степени, основной проблемой была всё же геометрия.

Мне попалась забавная книга Болтянского и Левитаса <<Математика атакует родителей>>, 1973, построенная в виде беседы
 родителей. Ее герои  обсуждают новую программу, сами ее постигают,
и в итоге проникаются. Переиздана в 1976 году. Была и такая попытка воздействия на общее мнение.
Книга эта была, наверно, в определенном смысле и неплохой. Более печальный памятник тогдашним настроениям -- сборник [ММЧ1978], 
вышедший в 1978 году.

По-видимому, часть учителей динамила программу и вела по старым учебникам (у нас не тюрьма народов!),
иные уходили из школ по выслуге лет...

Годы шли, а надежды на усовершенствование учителей и на постепенную отладку учебников не оправдывались...
Идея <<трудностей роста>> выглядела все более сомнительной.

\sm

{\bf\punct Тучи сгущаются.}
А профессиональные математики молчали (хотя не понимать происходящего было невозможно)...

Что касается учителей, то коль скоро они были обвинены в нехватке квалификации, то их никто не слушал (это кстати, известная стилистика
при проведении наших реформ, когда много позже вузовская профессура стала возражать против ЕГЭ, ее заклеймили в поголовной коррумпированности,
и больше на голоса из вузов можно было не обращать внимания).

Но другую сторону - миллионы возмущенных родителей,  которые, видимо, не достаточно внимательно изучили книгу Болтянского и Левитаса [БЛ1973],
 было невозможно не услышать.
 
 Летом 1977г. обучавшиеся по новой программе школьники
 начали поступать в вузы. Против реформы открылся новый фронт,
 теперь уже из вузовских преподавателей.

22 декабря 1977 г. ЦК КПСС и правительство выпустили постановление [ЦК1977]
<<О дальнейшем совершенствовании обучения, воспитания учащихся общеобразовательных школ и подготовки их к труду>>, где, в частности,
говорилось

\sm

{\it 
Школьные программы и учебники в ряде случаев перегружены излишней информацией и второстепенными материалами, что мешает выработке у учащихся навыков самостоятельной творческой работы. }

\sm 

 К каким именно предметам это относилось, в постановлении не уточнялось, но до сведения кого надо ЦК  довело опущенные в тексте
 постановления детали (вообще информацию, отсутствовавшую в газетах,
 было принято доводить до сведения на закрытых собраниях). Минпрос СССР задумался о смягчении программ. Минпросу РСФСР вся ситуация сильно не нравилась.
 Внешне все было тихо, новая программа стояла 
 непоколебимо.
 
 \sm
 
 {\bf \punct Реформаторы в ловушке.}
 Итак,  власти предложили математикам сократить программы. Вообще-то выкидывание
 тем в конце курса -- вещь очень неприятная для читающих, но с технической точки зрения
 она не представляет проблем (если только это не базовый курс, от которого
 многое потом зависит). У реформаторов, разумеется, не было желания отбрасывать такие темы. 
  Но действительным камнем преткновения была <<Геометрия>> Колмогорова 6-8 классов,
 начиная уже с 6-ого класса. Для упрощения учебника необходимо было
 отказаться от его идеологии и написать новый учебник. Но это означало бы
 открытую
 контрреволюцию.

\section{ Оверкиль}

\COUNTERS

{\bf \punct Выступление академиков.}
Весной 1978 года при Минпросе {\bf РСФСР} была создана контрреволюционная комиссия 
во главе с академиком А.~Н.~Тихоновым и методистом Ю.~М.~Колягиным. Колягин, как он и сам позднее признавал, долгое время был сторонником реформы, и приложил много сил к проведению реформы в жизнь (в частности он был автором многих
пособий по новой программе).

10 мая 1978 г Бюро Отделения математики
АН СССР издало постановление, где, в частности, говорилось

\sm 

{\it 
1. Признать существующее положение со школьными программами и учебниками по математике неудовлетворительным
как вследствие неприемлемости принципов,
заложенных в основу программ, так и в силу недоброкачественности школьных учебников.

2. Считать необходимым принять срочные меры к исправлению создавшегося положения,
широко привлекая, в случае необходимости, ученых-математиков, сотрудников АН СССР, 
к разработке новых программ, созданию и рецензированию новых учебников.

3. Ввиду создавшегося критического положения в качестве временной меры рекомендовать
рассмотреть возможность использования некоторых старых учебников.
}

\sm 

В декабре 1978 г. состоялось   общее собрание Отделения математики АН СССР [Стен1978].
На собрании выступал зам.министра просвещения 
СССР В.~М.~Коротков, который  предлагал несколько уменьшить школьную программу [КС2012] и продолжить реформу.
Выступали представители Минпроса РСФСР 
Г.~П.~Веселов и Ю.~М.~Колягин, которые были противниками Реформы. Кроме того, против реформы содержательно выступал 
представитель Физтеха М.~И.~Шабунин.

Со стороны академиков
главными атакующими были А.~Н.~Тихонов, В.~С.~Владимиров и Л.~С.~Понтрягин.
Жестче всех, разумеется, был Лев Семёнович (он был известен <<суровостью характера>>).
Но основное нападение, видимо, вел Тихонов, который рассчитывал перехватить в свои руки контроль над изданием учебников.
Критика
была  аргументирована, а общее мнение, по-видимому, выразил А.~Д.~Александров: {\it и так всем ясно, что в общем положение
неудовлетворительно}. Говорилось о необходимости срочно создавать новые учебники, объявлять конкурс, 
и о том, что нужны промежуточные учебники для тех,
кто уже учится по действующей программе. Академики не захотели возвращаться к старым учебникам
(хотя Владимиров на собрании о возможности их временного использования  говорил).

\sm 

Обсуждение достаточно интересно, но  я его не реферирую. На стороне Колмогорова
были Л.~В.~Канторович и С.~Л.~Соболев, которые говорили
много. Приведу лишь пару характерных фраз:

\sm

КАНТОРОВИЧ: {\it Так мне представляется, что большая работа по созданию новых учебников, которая
была проделана, это просто гражданский подвиг Колмогорова.}

\sm 

СОБОЛЕВ: {\it  Моя точка зрения в том, что новые изменения, это есть крупное достижение... Это все огромная
работа, проделанная Министерством, а также группой математиков во главе с А.~Н.~Колмогоровым.
Безусловно, очень много осталось погрешностей... Главным образом они объясняются,
видимо, слабой подготовкой учителей. }

\sm 

Еще один интересный момент

\sm 

ТИХОНОВ. {\it  Я по этому поводу вчера говорил с Министром высшего и специального среднего образования В.~П.~Елютиным.

Он говорил, что мы <<не предполагаем знаний математического анализа у поступающих в вузы, и вопросы
по разделу <<Математический анализ>>  включаются в билеты для поступающих в вуз только по моему указанию>>.
С его точки зрения, знания по этому разделу важны для тех учащихся, которые пойдут в вуз%
\footnote{Смысл этой фразы может показаться непонятным. Среди тогдашних прореформаторских речей
мелькала сентенция, что нельзя лишать не поступающих в вузы знания анализа.}.

В отношении всего курса математики он сказал так <<Лучше меньше, да лучше>>.}.

\sm 

 От себя скажу, что математики в вузах энтузиазма по отношению к матану в школе, как будто, никогда не испытывали
 (хотя на знакомство студентов с производной фактически, но не юридически, опирались, и знакомство это было полезным).

\sm 

Собрание проголосовало.
Предыдущая резолюция была существенно смягчена:

\sm 

{\it 1. Признать существующее положение со школьными программами и учебниками по математике неудовлетворительным.}

\sm 

Кроме того, было сказано

\sm 

{\it 2. Считать вновь представленную Министерством просвещения СССР программу по математике для средней школы неудовлетворительной.}
\sm

Вопрос о возможности использования старых учебников больше не поднимался (цели академиков были уже иные).

\sm 

За резолюцию было 26 человек, воздержались то ли трое, то ли двое. Сколько было всего голов
 было тогда в Отделении, я не считал.

\sm 

{\bf\punct Публичная полемика}. Цитирую Понтрягина

\sm 

{\it
В связи с развернувшейся на страницах упомянутого журнала} [<<Математике в школе>>]
{\it  дискуссией академик-секретарь Отделения математики АН СССР Н.~Н.~Боголюбов 
попросил журнал опубликовать полный текст решения общего собрания Отделения по этому вопросу (копия письма была послана министру просвещения СССР). 
Главный редактор журнала Р.~С.~Черкасов счел целесообразным ответить отказом...}

\sm

 Колягин утверждает, что резолюция АН СССР не была опубликована  из-за позиции
Минпроса СССР. В третьем номере журнала все же была опубликована  статья академиков
Владимирова, Понтрягина и Тихонова [ВПТ1979], в четвертом номере был ответ академиков Канторовича и Соболева [КС1979].

\sm 

В июне 1979 года умер Алексей Иванович Маркушевич (1908-1979).
Его последняя посмертная статья в защиту реформы была опубликована в <<Математике в школе>> в том же 1979 году [Мар1979].

\sm 

Математическая общественность не прекращала песни об отсталых учителях.
Из телевизора и из уличных матюкальников гремел \href{https://www.youtube.com/watch?v=Ti_2vCkMRl8}{\it (cссылка)}{ встречный голос}  
{Аллы Пугачевой} (чьими только устами не может
глаголать истина).

\sm 

{\bf \punct Позиционная война.}
События происходили и в непубличном пространстве. 

Тихоновцы начали писать новый комплект учебников.

Министр просвещения СССР М.~А.~Прокофьев предложил А.~Д.~Александрову исправить учебник стереометрии, см. [Вер2012].
Сделать это было невозможно, и Александр Данилович
 засел в 1979 году писать учебник стереометрии.

Относиться к этому можно двояко: с одной стороны в этом можно видеть борьбу академиков 
за хлеб насущный, с другой стороны ведь необходимо было что-то делать
и делать быстро. 

В Харьковской области начались опыты
с учебником А.~В.~Погорелова (он сам был из Харькова). Тихоновцы написали учебники и приступили к их  экспериментальной обкатке. А.~Д.~Александров издал свой
учебник как серию препринтов, и он начал обкатываться в нескольких
ленинградских школах.

Так или иначе, никаких крупных решений не было,
власти, надо думать, надеялись, что учебники упростят,
 реформаторы технически не могли этого сделать без отказа
 от <<Геометрии>> Колмогорова.
 Наверно позиционная война могла бы тянуться еще годами.

Возможностей для публичного обсуждения в профессиональных изданиях
не было. Главред (1958-1991)   журнала <<Математика в школе>> Черкасов входил
в число авторов комплекта учебников,  лично
Колмогоров был членом редакции, да и б\'oльшая часть редакции была реформистской. Триплет академиков смог пробиться в журнал, но
вообще это было <<не место для дискуссий>>. <<Успехи математических наук>> не стали бы публиковать статей против Колмогорова
(да и главредом  был ближайший друг Колмогорова - П.~С.~Александров).

\vspace{22pt}

Понтрягин  каким-то образом (см. [Пон2008]) нашел ход в главный журнал ЦК КПСС - <<Коммунист>>.
Смертельная статья 
<<\href{http://http://mat.univie.ac.at/~neretin/misc/reform/pontryagin-communist.html}{\it  О математике и качестве её преподавания}>>
вышла в сентябре 1980 года.

\sm 

Реформе настал конец.

\section{ Головокружение от успехов (вместо заключения)}

\COUNTERS

{\bf\punct Реформация.}
Посмотрим на всю эту историю <<с птичьего полета>>. Откуда же  проистекла эта роковая стратегическая ошибка?

\sm

В 50-60е годы были временем расцвета дополнительного школьного математического образования,
кружки, олимпиады, факультативы, лекции для школьников, массовый выпуск разнообразной популярной литературы.
Тогда же в Москве появились матшколы.

Участники этого движения были разными. Самый яркий и известный
 деятель  -- Николай Николаевич Константинов --  был весьма умеренным, не рвался кому-либо силой навязывать математику,
к раздуванию программ не стремился, а мат.школьный <<мат.анализ>> организовал по принципу <<бега на месте>>. Это нетривиальное
педагогическое решение
оказалось более удачным, чем казавшиеся более естественными движение вверх или развитие элементарной математики 
на ее собственной основе%
\footnote{Забавно, что такая точка зрения в чем-то была близка к старым традициям 
русской и советской школы: курсы математики были ориентированы
на развитие разума и до начала 60х не стремились раздуванию формального объема
информации.}. Но это был именно методика <<Константиновской системы>> (в конце 70х в нее входили
179, 91, 57 школы Москвы, ранее также 7 и 444).

Судя по всему (см. [Мар1964]), основной замысел реформы сложился в конце 50х-начале 60х годов внутри определенной  тусовки, которая частью была
 связана с дополнительным школьным математическим образованием, а частью входила в советскую педагогическую элиту
 (причем эти два подмножества сильно пересекались).
 Если посмотреть на мелькающие в связи с реформой фамилии,
 то они, в основном,  были известными и заслуженными  людьми педагогического мира,
 авторами учебников (в том числе массовых), популярных книг, методических пособий,
 наименее известные из деятелей реформации однако обнаруживаются среди
 победителей конкурса учебников 1962-1964 (см. следующий раздел статьи).

 Колмогоров не имел опыта работа ни в массовой
 школе, ни в высшей школе на нематематических специальностях.
Он мог  что-то недопонимать.
 Но почему эту разницу не понимали участники тусовки?
  Они в самом деле общались со школьниками, но обычно 
 не со случайными школьниками,
 и обычно в обстановке, отличной от регулярного урока. Идея, что в школах неправильно учат, была для этой тусовки вполне
 естественна. Идея была отчасти справедливой (в тусовке были и методисты, имевшие возможность наблюдать это воочию).
 Но конструктивный ответ <<Как правильно?, если речь идет не  о 0.2 процентах
 школьников, а и обо всех прочих тоже,
 дать  было совсем
 не  просто.

 В принципе, школьникам можно рассказывать много элементарных сюжетов из неэлементарной математики, если те слушают добровольно и ничем
 никому не обязаны
 (или если они тщательно отобраны). При переходе сюжета в обязательный он, в среднем, теряет привлекательность, потому
 что у ученика не остается выбора,
 понимать его до конца или пропустить. То, что было бы двигателем познания, превращается в препятствие.
 Точно так же в случае добровольности рассказ не окажется ловушкой для рассказывающего - не поняли, так и не поняли.
А.~Ю.~Михайл\'овская (о которой  говорилось выше), среди прочего, на факультативе рассказала  семиклассникам
 (неспециализированного класса...) про производную (лучше, чем это обычно делают в хороших вузах).
 Это в принципе можно, если школьники не случайны, если их не обязывают
 это в дальнейшем знать, и если ты Михайловская. 
 Если, если, если...

 Элементы экстремизма в Программе-1959 показывают влияние этой тусовки. Всё могло бы кончиться благополучно,
 остановившись
 на уровне учебников Кочетковых и Погорелова,  нереальные мечты были бы оставлены при столкновении  с действительностью...
 Но  даже авторы Программы-1959 рассматривали ее как недостаточно решительную и мечтали о революции...
  
 Маркушевич в 1964 году по каким-то  причинам оставляет пост заместителя министра просвещения
 и снова становится вице-президентом АПН. Настроения в АПН, судя по разработанной к 1964 году,  были  революционными,
  в 1965-1968 году в математическую группу вводят
  новых членов, и она вполне становится революционным центром (см. ниже список членов).

 Колмогоров уверовал в уже существовавшую программу элитной педагогической тусовки, подкрепленной
  аргументами от группы Бурбаки
 (а там  люди были грамотнее некуда, как впрочем и Маршалл Стоун в США).
 И он имел все основания считать, 
 что опирается на лучших представителей советского педагогического мира (оно так и было, но это была
 неслучайная выборка из подобных представителей), 
 что многое уже просчитано и продумано. 
 
Так или иначе, без участия Андрея Николаевича Колмогорова Революция не имела бы шансов на воплощение в жизнь, и, как это ни грустно дополнительно
признать,
он внес в революционную программу и свой собственный вклад.

 \sm

{\bf\punct Контреформация.}  
 Почему математики молчали? Потому что человек, открыто выступивший против, стал бы мишенью для научных кланов, поддерживавших Колмогорова,
  (даже без команды <<Ату!>>, а просто по общепринятым правилам игры). 
 Параллельно, выступивший вышел бы  из под зонтика своего клана, потому как на объявление клановой войны  его никто бы не уполномочил...
  
 Впрочем, выступления, наверно, были. Известно [КС2012] <<Мнение совета отделение математики механико-математического
 факультета МГУ>> о том, что учебники нуждаются в серьезной доработке (к сожалению, документ в публикации 
 не датирован,  по содержанию это похоже уже на осень 1978 года, то есть Бюро отделения математики АН СССР уже высказалось).
 И при нем особое мнение С.Б.Стечкина: {\it Вопрос о том, что...
 черного кобеля не отмоешь добела.} Но это уже осень 1978 года, когда встречный ветер начал подниматься.
 Похоже на то, что А.~Г.~Курош (1908-1971) относился к школьно-педагогической деятельности Колмогорова
 отрицательно (в [Шир2003] приводится обширное письмо Колмогорова Курошу от 5.1.1964, но из письма не понятно, на что в точности А.~Н. 
 возражает).

 <<Бунт>> академиков случился уже после постановления ЦК, при поддержке Минпроса РСФСР,  при поддержке вузов.
 Вот тогда-то академики и не испугались.
 
 Выступление академиков было общим,  понятно, что почти всем реформа сильно не нравилась, но, конечно, 
 дискуссия декабря 1978 года имела  клановую структуру.... Все равно 
 надо хоть и за это академикам сказать <<спасибо>>.
 
 Положение усугублялось тем, что в ходе реформы было произведено <<сжигание мостов>>. Если в отношении <<Алгебры>>
 и <<Элементов анализа>> еще был возможен отход на позиции, близкие к предреформенным, путем упрощения реформистских 
 учебников (и в итоге это было проделано), то упоминавшаяся выше <<коренная перестройка сложившегося курса геометрии>> влекла необходимость
 обратной <<коренной перестройки>>.
 
 \sm 
 
 Начавшийся после выступления Л.~С.~Понтрягина выход из ловушки оказался мучительным.
 Встал вопрос, как быть с учебниками.
 Казалось бы решение напрашивалось: 
 вернуть <<Кочеткова>> (который  еще за четыре года до того использовался в школе) и запустить <<Погорелова>>, 
 который с 60х годов лежал без дела.
 Оба учебника были умеренно реформистские, а с Реформой Колмогорова были  напрямую не связанны. 
 По алгебре средних классов можно было бы оставить реформистский учебник под редакцией Маркушевича, основательно  его почистив.
 Некоторое время зализывать раны и параллельно объявлять конкурс... 

  \sm 
  
 Но ни крупных деятелей, ни общественных структур, ищущих разумного выхода, не оказалось. АПН, в чьем прямом ведомстве были учебники,
 погрязла по уши в реформе.
 Академия наук СССР, которая 7 лет доброжелательно наблюдала за начинавшейся реформой, а потом
 пять лет спокойно взирала на безумие в школе, такой структурой не оказалась.
 Началась схватка академических кланов на ниве просвещения. Тихонов стал
 захватывать математику в школе в свою сферу влияния. Колмогоров не хотел уходить. Великий геометр А.~Д.~Александров
 написал свой учебник по Стереометрии. В этой битве бизонов (при всех
 личных достоинствах ее участников)
 не могло быть честного судейства, и никто, кроме оных,
 не имел шансов вмешаться. Уже позже, в 1986-88 гг проводился новый конкурс учебников. 
 Он, как и конкурс 1962-1964 годов, был официально
 конкурсом тайным, но теперь на нем под  девизами <<скрывались>>
  уже опубликованные учебники от академиков%
  \footnote{На этом конкурсе, завершившимся в 1988г.,
 	учебник Колмогорова занял пятое место, уступив достойным противникам
 	-- <<Атанасяну>>, Погорелову, <<Александрову>> и Болтянскому.}..

 М.~А.~Прокофьев, который в течение всей этой истории  был министром просвещения СССР (1966-1984),
 напоследок сделал разумный шаг - в 1982 году пустил <<Погорелова>> в дело
 (написанный в конце 60х  учебник Погорелова совсем не идеален, см. комментарий при ссылке
 [Пог1982]%
 \footnote{Этот учебник, будучи реформистким в смысле Программы-1959, бескровно вводил векторы и геометрические преобразования. Однако он унаследовал и установку на чрезмерную
 строгость в начале курса, выдвинутую в 60х годах Колмогоровым.
В итоге начало курса сложновато, а школьники имеют право не понимать, почему одни очевидные высказывания  надо доказывать,
а другие не надо. Здесь учитель должен выруливать в непростой ситуации. По-видимому, нужна была серьезная переработка начала учебника с участием хороших педагогов, что сделано не было,
быть может из-за недостаточной гибкости самого Погорелова.}, но в целом он был  добротен, уже был <<под рукой>> и был известен). Решение, по-видимому, было продвинуто
 Виноградовым и Понтрягиным.
   
   В остальном нива просвещения была поделена между тихоновцами и колмогоровцами (правда учебники имени колмогоровского проекта
    сменили авторов  и содержание), А.~Д.~Александров остался, но оказался оттесненным в уголок.
 
 Процитируем В.~И.~Арнольда
 
 \sm 
 
 {\it 
 Уважаемые мною люди, А.~Д.~Александров, Погорелов, Тихонов, Понтрягин} [Понтрягин -- оговорка] {\it — все приняли участие и все написали плохо.
 Я могу точно сказать, что плохо написал Колмогоров, скажем, ну и про других тоже знаю; учебники,
 которые они предложили, могу критиковать, но не могу предложить своего учебника... 

}

 \sm

 Катастрофа была остановлена, но отыграть  математику в школе на дореформенный  уровень не удалось. Учебники математики стали сложнее, чем были до реформы. 
 Семь лет вгона школьной математики в ступор не могли пройти даром. Кроме вещей понятных и без слов,
 были трудно оцениваемые количественно эффекты типа <<деградации инфраструктуры>>...
 До реформы математика была благополучным предметом, она развивала мышление школьников
 и подтягивала школу вверх. На время реформы ее роль  в массовой школе  стала превращаться 
 в противоположную... На исполнение ошибочного  стратегического плана
 были истрачены силы многих достойных людей,  силы эти заслуживали лучшего применения,
 и при ином сценарии, вероятно, были бы лучшим способом применены...
 Вряд ли эта история могла бы обойтись без последствий для пед.вузов...
 Отношение общества к математике ухудшилось, не без помощи  Реформы,
 хотя и по разным иным причинам тоже... Не пошла Реформа и на пользу 
 науки вообще
 (кто может просчитать, какой вклад Реформа внесла в общественную реакцию против естественных
 наук, грянувшую в Перестройку).
  
  \sm 
  
 Выпускной школьный экзамен по геометрии с 1977 года канул в Лету%
 \footnote{Если мне не изменяет память, в середине 80х усилиями академика А.~П.~Ершов в школу был введен новый предмет <<Информатика>>, что было действием вполне разумным.
 Это мероприятие сопровождалось изыманием <<часов>> у школьной математики, что на пользу школьной математике не шло, но тоже было действием разумным. С одной стороны, о программировании до этого говорилось в курсе <<Алгебра и начала анализа>>, с другой - математики не сумели разумно распорядиться выданными им <<часами>>. Я не пытался выяснять, как менялось число <<часов>> на математику в 1959-1990гг.}

 \sm 
 
 Навалилось много иных бед, напрямую не связанных со школьными  учебниками. В вузах  шла реакция против математиков,
 которые достали всех в предыдущее десятилетие.
 За разными внутренними баталиями профессиональное сообщество просмотрело, что математика вступительных экзаменов превратилась 
 в отдельную науку и стала давить на школу. Свои проблемы вносило  всеобщее среднее образование, введение которого по времени совпало с Реформой.
 Уже позже это серьезно усугубилось в связи с развалом профессионально-технического образования, который привел
 к объединению всех молодых людей в единой старшей школе 
 (не надо экстраполировать эту последнюю проблему на первую половину 80х).
 
 \sm
 
 Вскоре после остановки Реформы настала эпоха социальных потрясений. На этом фоне скорее стоит удивляться,
 что школьная математика оказалась на удивление живучей,
 и что до Кузьмин\'овских реформ, начавшихся в 2001 году, она  сохраняла относительно приличный уровень.

\sm

{\bf \punct Двадцать лет спустя.}
 Приведем слова И.~Ф.~Шарыгина, сказанные в 2001 году при виде начавшихся тогда новых реформ:
 
 \sm 
 
{\it  В течение тридцати с лишним лет в Советской России и Советском Союзе, медленно, но не мучительно, 
формировалась система математического образования,
 которую потом назвали Советской. Пожалуй, лишь к началу пятидесятых годов эта система сформировалась полностью. 
 Следующие два десятилетия Советское математическое образование развивалось и совершенствовалось.
 Вероятно, главным итогом этого развития явились немногочисленные пока еще специализированные математические школы и классы.
 В начале эти классы были явлением безусловно прогрессивным. Но одновременно с их появлением начался раскол некогда единой системы
 школьного математического образования. Начавшийся на верхних этажах школьного здания этот раскол пошел вниз
 и сегодня почти достиг начальной школы.

В начале семидесятых годов по инициативе выдающегося математика А.~Н.~Колмогорова в Советском Союзе началась реформа математического образования - 
первая из до сих пор непрекращающейся вереницы реформ. На мой взгляд, эта реформа была недостаточно обоснованной,
плохо продуманной и совсем скверно реализованной. По мнению других, большею частью близких к Колмогорову реформаторов%
\footnote{Cм., например,  [Абр1988], [Абр2003], [Aбр2010],  [Гнед1993], [Чер1988], [Чер1993], [ГЧ1993], см. также [Тих2009].} 
реформа была необходимой и хорошо проведенной. Не буду спорить. Но если мы хотим указать точку отсчета,
с которой началась, вначале очень медленная, деградация системы математического образования Советского Союза и России, 
то она приходится примерно на середину семидесятых годов. Забавно также, что период реформирования в системе образования начался 
с реформирования самого благополучного предмета - математики, 
и инициировали это сами математики. (Не ведаем, что творим?).}

\sm 

{\bf\punct Сорок лет спустя.} Здесь можно сказать еще меньше утешительного. Последние 10 лет школьная математика ползет вниз на глазах, а геометрия вообще уходит в прошлое.

Математическое сообщество 20х - начала 70х годов XX века в целом успешно вело конструктивную образовательную политику.
К сожалению, как мы видели, в 1959-1968гг. часть математической тусовки   сформулировала программу с положительными целями, хоть и авантюристичную, а в 1968-1980 приложила недюжинные усилия для проведения ее в жизнь. В 1980-1985
другая часть сообщества сумела оздоровить ситуацию. 

В последние 30 лет российская профессионально-математическая тусовка в отношении школьного и
вузовского образования (за исключением узкопрофессионального сектора) не была способна ни к конструктивным действиям наступательного характера, ни к обороне, ни к постановкам проблем, ни даже к анализу постепенно  ухудшавшегося  положения.

 \section{Дополнение.\\ Соратники  Колмогорова и Маркушевича}
 
 \COUNTERS

 \begin{center}
 \bf\large
A. Данные о педагогическом статусе до 1970г. авторов программы-1966-1968 года и авторов комплекта реформистских учебников
\end{center}
 
  Ниже список, включающий в себя авторов программы 1965-1968гг и авторов ударного комплекта учебников.
 Меня интересует их предшествующая биография {\bf как педагогов} (данные о научных заслугах Виленкина и Болтянского
 никак не отражаются). Приводимые ниже данные очень обрывочны
 и в основном содержат списки их научно-популярных и педагогических {\bf книг} до начала реформы, т.е. до 1969-1970гг.
 Я использовал [Кат-РНБ], [Кат-РГБ] и [Кат-МГУ]. Учитывались издания тиражом выше 1 000экз. Стоит иметь в виду, что 
 в те времена книги издавались немалыми тиражами, но наименований издаваемых книг было немного. Сам факт публикации
 книги свидетельствовал об определенном общественном положении автора.

 Педагогические и научно-популярные публикации легко различаются по их названиям.

{\bf Болтянский  Владимир Григорьевич} (род. 1925),
член-корреспондент АПН РСФСР (1965), член-корреспондент АПН СССР (1968),
автор учебника

{\footnotesize Болтянский В.Г. , Яглом И.М.
{\it Геометрия. Учебное пособие для 9 класса средней школы.} - М.: Учпедгиз. 1963 (тираж  2 300 000), 1964 (тираж 2 000 000), 
издана также в Киеве, 1963, тираж 230 000.

}

Также автор книг

{\footnotesize
Яглом И. М., Болтянский В. Г. {\it  Выпуклые фигуры.} ГИТТЛ, 1951 (тир. 25 000)

Болтянский В. Г. {\it Равновеликие и равносоставленные фигуры,} Гостехиздат, 1956 (тир. 40 000)

Болтянский В. Г. {\it Что такое дифференцирование,} Гостехиздат, 1955 (тир. 50 000), Физматгиз 1960 (тир 35 000)

Болтянский В.Г {\it Огибающая,} Физматгиз, 1961 (32 000)

В.Г. Болтянский, И.М. Яглом. {\it Векторы в курсе геометрии средней школы,} Учпедгиз, 1962 (тир. 25 000)

В.Г. Болтянский, И.М. Яглом. {\it Преобразования. Векторы.} (для учителей) Просвещение 1964 (тир. 56 000)

Гохберг И. Ц., Болтянский В. Г.  {\it Теоремы и задачи комбинаторной геометрии.} — М.: Наука, 1965; (тир. 23 000)

В.Г. Болтянский и др. {\it Сборник задач московских математических олимпиад,} Просвещение, 1965 (тир. 122 000)

Болтянский В. Г., Виленкин Н. Я. {\it Симметрия в алгебре,} Наука, 1967 (тир. 50 000)

Болтянский В.Г. (редактор) {\it Комплексы учебного оборудования по математике,} Педагогика, 1971 (тир. 25 000)

Гохберг И. Ц., Болтянский В. Г. {\it Разбиение фигур на меньшие части,} Наука, 1971

 В. Г. Болтянский, Ю. В. Сидоров, М. И. Шабунин.
{\it Лекции и задачи по элементарной математике} - Москва : Наука, 1971 (тир. 100 000)

}

{\bf Вейц Борис Ефимович} (1921-2007),  к.ф.м.н., 1965, Мурманский пед.

{\bf  Виленкин Наум Яковлевич}  (1920-1991)

{\footnotesize
Виленкин Н.Я. {\it Метод последовательных приближений}, Физматгиз, 1961 (тираж 30 000), 1968 (тираж 100 000)

Виленкин Н.Я. {\it Рассказы о множествах,} 1965 (тираж 50 000), 1969 (тираж 100 000)

Болтянский В. Г., Виленкин Н. Я. {\it Симметрия в алгебре,} Наука, 1967 (тир. 50 000)

Н. Я. Виленкин, Р. С. Гутер, С. И. Шварцбурд и др.
{\it Алгебра: Учеб. пособие для IX-X классов сред. школ с матем. специализацией}  - Москва : Просвещение, 1968
(тир. 40 000)

Виленкин Н.Я. {\it Комбинаторика,} М.Наука, 1969 (тир. 100 000)

Виленкин Н.Я.,  Шварцбурд С.И. {\it Математический анализ. Учебное пособие для школ с математической специализацией.} Просвещение, 1969

Н. Я. Виленкин, А. А. Кочева, И. В. Стеллецкий 
{\it Задачник-практикум по элементарной алгебре: Для студентов заочников физ.-мат. фак. пед. ин-тов}  Просвещение, 1969 (тир. 30 000)

Виленкин Н.Я., Литвиненко В.Н., Мордкович А.Г. {\it Элементарная математика. Учебное пособие для студентов-заочников пед-институтов}, Просвещение 1970(тир. 35 000)

Виленкин Н.Я., Михайловская А.Ю. {\it Элементы теории множеств. Факультативный курс математики для VII класса.}
Методические указания. М., 1970 (тир. 10 000)

}

{\bf  Гусев Валерий Александрович}, (1942-), работал в колмогоровском интернате, кандидатская диссертация,  1971.
Появлялся в качестве соавтора колмогоровского учебника геометрии 8 класса, 1972-1980

{\bf 
Демидов Иван Тимофеевич} (1909-1975) , Мурманский пед.

{\footnotesize
Демидов И.Т. {\it Основания арифметики,} Учебное пособие для пед-вузов, Учпедгиз,  1963

}

{\bf Ивашёв-Мусатов Олег Сергеевич}  (род 1927).
Пасынок Колмогорова,
доцент кафедры математического анализа мехмата МГУ.
В 1958—1969 годах преподавал также на химическом и геологическом факультетах.
Появляется в числе авторов колмогоровского учебника <<Алгебра и начала анализа>> с 1975г.

{\bf Ивлев Борис Михайлович} (1946-1990). 
Выпускник Колмогоровского интерната. В 1964-1972 студент, потом аспирант МГУ.
Появляется в числе авторов колмогоровского учебника <<Алгебра и начала анализа>> с 1976г.

{\bf Клопский Владимир Михайлович} (1926-1982),
защитил кандидатскую в 1972 году в Ярославле, работал в Курском пед.институте

{\footnotesize
Клопский В.М., Ягодовский М.И  {\it Геометрия, 9-10 классы}, Просвещение, 1966 (тираж 15 000 экземпляров)

}

{\bf Макарычев Юрий Николаевич}  (1921-2007)   Кандидат с 1964 года.

{\footnotesize 
Макарычев Ю.Н. {\it Система изучения элементарных функций в старших классах средней школы.}
Учебно-методическое пособие для учителей, Просвещение, 1964 (219с, тираж 25 000)

  Ю. Н. Макарычев, К. И. Нешков,
 {\it Алгебра: Учеб. материалы для VI класса}, АПН РСФСР.   Москва, 1966 (тир 800), 1967 (тир. 1 500). 

Макарычев Ю.Н., Нешков К.И. {\it Математика в начальных классах} (под редакцией А.И.Маркушевича), Педагогика, 1970 (260 000 экз)

}

{\bf Миндюк Нора Григорьевна,}  (1933-2016),
с 1964г. работала в секторе математики Института методов обучения АПН СССР,  кандидат с 1966 года

{\footnotesize
Миндюк, Н. Г. {\it 
Математика: Учеб. материалы для V класса} / АПН РСФСР.  - Москва : Просвещение, 1965 (тир. 2 000)

}

{\bf   Маслова, Галина Герасимовна.} (1920-) Кандидат, 1954.
 В качестве места работы упоминается	Научно-исследовательский институт содержания и методов обучения АПН СССР

 {\footnotesize 
Н. Н. Никитин, Г. Г. Маслова.   {\it Сборник задач по геометрии для 6-8 классов восьмилетней школы,}
(основной школьный задачник 1957- 1971), в том числе
Планиметрия 6-7 классы,
1957 (тир. 2 000 000),
1958 (тир. 2 000 000),
1959 (тир.  600 000),
Геометрия,
1961 (тир. 1 000 000),
1962 (тир.  1 925 000),
1963 (тир. 1 500 000),
1964 (тир. 1 000 000),
1965 (тир.  1 100 000),
1966 (тир. 1 000 000),
1967 (тир. 1 200 000),
1967 (тир. 100 000),
1968 (тир. 1 200 000),
1969 (тир. 1 100 000),
1970 (тир. 1 200 000)
1971 (тир. 1 000 000)

 Маслова   Г. Г.  {\it Методика обучения решению задач на построение : в восьмилетней школе.} АПН РСФСР, 1961. (тир. 41 300)
 
 Маслова  Г. Г.    {\it О программированном обучении математике,} Просвещение 1964 (тир. 37 000)
   
    Маслова  Г. Г. (редактор) {\it Повышение эффективности обучения математике,} 1971

}

  {\bf  Муравин Константин Соломонович} (1920—1993),  кандидат 1967 года, несколько мало-тиражных изданий, а также
  
  {\footnotesize 
   Муравин К.С.,
  {\it Математика: Учеб. материалы для V класса},
   Акад. пед. наук. Ин-т общего и политехн. образования. - Москва : [б. и.], Ч.2, 1966 (тир. 800), Ч.2, 1965, (тир. 1 500), Ч.3, 1965 (тир. 2000)
    
 Муравин К.С.,  {\it Самостоятельные и контрольные работы по алгебре для 8-летней школы. Пособие для учителя,}
 Просвещение, 1965(тираж 175 000), 1971 (тираж 300 000)

 Муравин К.С., Фрейдлин Е.Г. {\it Сборник задач по алгебре для 6-8 классов,} Просвещение, 1964 (тираж 77 000), 1968 (тираж 250 000)

}

  {\bf  Нешков Константин Иванович} (1923-)  кандидат с 1956 года
   
   {\footnotesize
  Нешков К.И.    {\it Система изложения курса арифметики в V классе}, Изд. АПН РСФСР, 1963 (293с., тираж 55000)
  
  Ю. Н. Макарычев, К. И. Нешков ;
  {\it Алгебра: Учеб. материалы для VI класса},  Акад. пед. наук РСФСР. Науч.-исслед. ин-т общего и политехн. образования. - Москва, 1966 (тир 800), 1967 (тир. 1 500). 
    
  Нешков, К. И.,  К. А. Краснянская 
{\it Математика  : Учеб. материалы для IV класса} Акад. пед. наук СССР. Ин-т общего и политехн. образования.  - Москва : Просвещение, 1966, Ч.1.(тир. 3 000),  Ч.2.(тир. 1 300).
  
     Нешков К.И. Пышкало, А.М., {\it Самостоятельные работы в курсе арифметики V класса. Дидактический материал.}
     Просвещение, 1964 (330с, тираж 75 000)
     2 изд, 1967 (328с, тираж 75 000)
    
     Нешков К.И. {\it 
     Математика: Учеб. материалы для IV класса} / Акад. пед. наук СССР. Ин-т общего и политехн. образования. - Москва : [б. и.], 1963, Ч.1, Ч.2 (тир. 3 000), Ч.3, 1964 (тир. 2 000)  1967
     
    Нешков К.И. Пышкало, А.М.,  {\it Математика в начальных классах} (под редакцией А.И.Маркушевича), Просвещение, 1968  (тир. 200 000)

    Макарычев Ю.Н., Нешков К.И. {\it Математика в начальных классах, ч.1} (под редакцией А.И.Маркушевича), Педагогика, 1970 (260 000 экз), часть 2 (1970) (тир. 260 000)
    
}

  {\bf    Нагибин Федор Федорович}  (1909-1976),
  Работал в Вятском пединституте, видимо, с 1939 года был деканом (в Войну короткое время - ректором%
  \footnote{Формально должность называлась директор.})
    
    {\footnotesize 
 Е. С. Березанская, Ф. Ф. Нагибин, 
 {\it Сборник вопросов и упражнений по алгебре и тригонометрии для 8–10 классов,} Учпедгиз, 1951 (тир. 50 000), 1955 (тир. 60 000)
 (переведена на китайский язык)
 
 Е. С. Березанская, Ф. Ф. Нагибин, {\it Упражнения для устных занятий по алгебре,} Учпедгиз, 1949 (тир. 25 000)

    Нагибин Ф.Ф. {\bf Математическая шкатулка.} Москва: Учпедгиз, 1958 (тираж 100 000), 1961 (тираж 150 000), 1964 (тираж 245 000), 
    Знаменитая книга, многократно переиздавалась,
    переведена на японский и китайский языки

Е. С. Березанская, Н. А. Колмогоров и Ф. Ф. Нагибин, {\it Сборник задач и вопросов по геометрии,} Учпедгиз, 1962 г. [Николай Андреевич Колмогоров, не путать с Андреем Николаевичем]

Нагибин Ф.Ф. {\it Экстремумы.  Пособие для учащихся старших классов,} Просвещение, 1966 (тир. 110 000)

Ф.Ф. Нагибин, А.Ф. Семенович и Р.С. Черкасов, {\it Геометрия, для 6-8 классов,} Просвещение 1967, (тираж 11 000)

Н. А. Колмогоров, Ф. Ф. Нагибин, В. В. Чудиновских 
    {\it Сборник задач для подготовки учащихся средних школ к математическим олимпиадам};
    Волго-Вятское кн. изд-во, 1968. (5000 экземпляров) [Николай Андреевич Колмогоров]
    
}

 {\bf Семен\'ович Александр Федорович}   (1920-?)
 
 {\footnotesize
  Семенович А.Ф., Воробьев Г.В.  {\it Первые уроки геометрии (из опыта работы учителя)}, Учпедгиз, 1958 (тир. 15 000)
  
  Семенович А.Ф. {\it Задачи на доказательство по готовым чертежам,} Свердловск, 1960 (29с.,тираж 8500)
   
 Семенович А.Ф.  {Учебное пособие по проективной геометрии (для студентов-заочников пединститутов)}, Учпедгиз, 1961
 (тир. 25 000)
  
 Семенович, А.Ф. {\it 
Геометрия  : Пробный учебник для седьмого класса} / Ульян. гос. пед. ин-т им. И. Н. Ульянова. - Ульяновск : Кн. изд-во, 1962. 
 
Семенович, А.Ф. {\it 
Геометрия  : Пробный учебник для восьмого класса : (Пособие для учителей)}.
 Ульян. гос. пед. ин-т им. И. Н. Ульянова. - Ульяновск : Кн. изд-во, 1963 [вып. дан. 1964] (тир 1 500)
  
 Семенович А.Ф.  {\it Геометрия Пробный учебник для восьмого  класса (пособие для учителей),} Ульяновск, 1963 (тираж 1500)

  А.Ф. Семенович, Ф.Ф. Нагибин, и Р.С. Черкасов, {\it Геометрия, для 6-8 классов}, Просвещение 1967, (тираж 11 000)
  
}

{\bf Семушин Алексей Дмитриевич} (1915—?), канд. пед наук, 1955
 
{\footnotesize 
 Семушин А.Д. (редактор) {\it Изготовление наглядных пособий по геометрии и их применение на уроках.} Сборник статей, Изд АПН РСФСР, 1953 (тир. 15 000)
 
 Семушин А.Д. (редактор) {\it Вопросы методики математики в средней школе.} Сборник статей, Изд АПН РСФСР, 1954
 (тир. 20 000)
 
 Семушин А.Д. (редактор) {\it Вопросы повышения качества знаний учащихся по математике,}  Сборник статей, Изд АПН РСФСР, 1954  (тир. 20 000)

 Семушин А.Д. (редактор) {\it Политехническое обучение в преподавании математике}. Из опыта работы в V-X классах.  Сборник статей, Изд АПН РСФСР, 1956 (тир. 20 000)

 Семушин А.Д. {\it Методика обучения решению задач по стереометрии,}  Изд АПН РСФСР, 1959 (тир. 28 000)

 Семушин А.Д. {\it О преподавании математике в школе в 1959/1960 году} (учебное пособие),  Изд АПН РСФСР, 1961
 (тир. 42 300)

 Семушин А.Д. (редактор)  {\it О преподавании математики в восьмилетней школе,} Изд АПН РСФСР, 1961
  (тир. 52 300)
 
 Гибш И. А., Семушин А. Д., Фетисов А. И. {\it Развитие логического мышления учащихся в процессе преподавания
 математики в средней школе: пособие для учителей.}
 1958 (тир. 30 000), второе издание
 
}
 
{\bf  Скопец  Залман Алтерович}  (1917-1984).    Заведующий кафедрой в Ярославском пединституте,
геометрии с 1962 (до этого - зав.кафедрой элементарной математики), дфмн, 1961

   {\footnotesize 
  Майоров В.М.,  Скопец З.А. {\it Задачник-практикум по векторной алгебре, Для студентов заочников физ-мат фак. пед институтов,} М. : Учпедгиз, 1961. (тир. 20 000)
 
 Жаров В.А.,  З. А. Скопец.    {\it Задачи и теоремы по геометрии. Планиметрия: пособие для пед. ин-тов}. - М. : Учпедгиз, 1962 (тир. 38 000). 
 (первый вариант, Ярославль, 1958, тир. 2 000)
  
{\it Вопросы совершенствования преподавания в средней школе} [Сборник статей] / Под ред. З. А. Скопеца и А. И. Голубева. - Ярославль : [б. и.], 1963.  (Доклады на научных конференциях/ Яросл. гос. пед. ин-т им. К. Д. Ушинского. Педагогика, методика; Т. 2. Вып. 1. Ч. 1) (тир. 1000)
 
 Жаров В.А., Марголите П.С., Скопец З.А. {\it Вопросы и задачи по геометрии. Пособие для учителей}. Просвещение, 1964
 (тир.106 000)
  
 Майоров, В. М.,   Скопец З.А.
{\it Векторное решение геометрических задач  (Задачник-практикум по спецсеминару)}:
Для студентов-заочников физ.-мат. фак. пед. ин-тов / В. М. Майоров; Москва : Просвещение, 1968. (тир. 20 000)

}

 {\bf     Суворова Светлана Борисовна}.
  Появляется в числе авторов <<Алгебры>> с 1976 года. Кпн с 1982
 
{\bf  Фетисов, Антонин Иванович} (1891-1979), окончил сельхоз-училище 1919, 
сдал экстерном экзамены за курс обучения в МГУ, 1928, кандидат 1946
  
  {\footnotesize
 Делоне Б.Н., Житомирский О.К., Фетисов А. И. {\it Сборник геометрических задач,}
 Пособие для учителей, Учпедгиз, 1941, Учпедгиз, 1951 (тир. 30 000)
 
Фетисов А. И. {\it Опыт преподавания геометрии в средней школе}. – 1946, 
 
Гибш И.А., Фетисов А. И., {\it Исследование решений задач с параметрическими данными,}  Изд.АПН РСФСР, 1952 (тир.10 000)
  
Фетисов А. И., {\it О доказательстве в геометрии.} Гостехиздат,  1954 (тир. 50 000),
переведена на польский, болгарский, английский, немецкий и французский языки.
 
Фетисов А. И., И.Н.Шевченко, В.Л. Гончаров,  Гибш И. А.  {\it Преподавание математики в школе в свете задач политехнического обучения.
Материалы в помощь учителю.} АПН РСФСР, 1953 (50 000),  Алма-Ата, Казучпедгиз, 1954 (тираж 3 000), АПН РСФСР, 1954 (тираж 30 000)
  
 Никитин Н.Н., Фетисов А. И. {\it Геометрия. Учебник для семилетней и средней школы, Ч.1,} Учпедгиз, 1956 (тираж 2 000 000),  1957 (тираж 900 000) 
  
 Фетисов А. И. Геометрия. {\it Учебник для 8-9 классов средней школы, 2 изд.,} 1957 (тираж 10 000)
  
Фетисов А. И. Геометрия. {\it Пробный учебник для средней школы, Ч.2, Стереометрия,} Учпедгиз, 1957 (тираж 10 000)

   Фетисов А. И.(ред.) {\it Преподавание математики.} Сборник статей, Изд.АПН РСФСР, 1957 (26 000)
 
 Гибш И. А., Семушин А. Д., Фетисов А. И.
 {\it Развитие логического мышления учащихся в процессе преподавания математики в средней школе: пособие для учителей.}
 Учпедгиз, 1958 (тир. 30 000)

 Фетисов А. И. Геометрия. {\it Учебное пособие по программе старших классов.} Изд.АПН РСФСР, 1963 (тир. 34 800); 
  
Фетисов А. И. {\it Очерки по евклидовой и неевклидовой геометрии.} Просвещение,  1965 (тир. 17 000);
  
Фетисов А. И. {\it Учебные материалы по геометрии для V класса. Ч.1., Ч.2,} Просвещение, 1965 (тир. 1 500)
  
Фетисов А. И. {\it Учебные материалы по геометрии для VI класса. Ч.1., Ч.2,} Просвещение, 1965 (тир. 2 000)
 
  Фетисов А. И.(ред.)
{\it Методика преподавания геометрии в старших классах средней школы,} Просвещение, 1967 (тираж 125 000)
  
Фетисов А. И. Геометрия. {\it Учебное пособие по программе старших классов.} 1963; 

 Перевел программную книгу [ПБДЛШГ1960]

 }

 {\bf    Черкасов Ростислав Семенович}  (1912-2002).
 После окончания аспирантуры МГПИ в 1948 направлен в аппарат Мин-ва просвещения РСФСР; работал инспектором вузов, 
 нач. отд. инспекции вузов, нач. отд. общенауч. дисциплин (1948–57). 
   Заведующий кафедрой методики преподавания математики Моск. гор. пед. ин-та им. В.П.Потемкина.
 После слияния ин-та с МГПИ оставался зав. каф. до 1985.
 С сер. 1960-х до начала 1970-х – декан мат. ф-та МГПИ,  проф.
 МГПИ с 1958 года по 1983 год.
 Главный редактор журнала «Математика в школе» (1958-1991).
 
 {\footnotesize
  Черкасов Р.С. {\it Сборник задач по стереометрии: Пособие для учителей сред. школы,} Учпедгиз, 1952 (тир 25 000), 1956. (тир. 35 000)
  
  А.Ф. Семенович, Ф.Ф. Нагибин, и Р.С. Черкасов, {\it Геометрия, для 6-8 классов,} Просвещение 1967, (тираж 11 000)
 
 А. И. Маркушевич, К. П. Сикорский, Р. С. Черкасов ;
 {\it Алгебра и элементарные функции:} Учебное пособие по математике /  Под ред. А. И. Маркушевича. - Москва : Просвещение, 1968 (тир. 100 000)
 
}

{\bf Чесноков Александр Семенович} (1926- ),
Кандидат пед. наук, 1979
  
{\bf   Шварцбурд, Семён Исаакович} (1918-1996), канд. пед. наук, 1961, докт. пед. наук 1972.
   Заслуженный учитель школы РСФСР (1962),  чл-корр. АПН СССР с 1968.
 В 1959 году организовал в 444 школе первый
 математический класс.

 {\footnotesize
 Шварцбурд С.И. {\it Системы уравнений., Методическая разработка курса алгебры VIII класса,} Изд. АПН РСФСР,  1955 (тир. 40 000)
  
 Шварцбурд Б.И., Шварцбурд С.И. {\it Задачи по математике для школ с машиностроительной специализацией,} Пособие для учителей IX-X классов.
 М.Учпедгиз, 1962 (тир 54 000)
  
 Шварцбурд С.И. {\it Математическая специализация учащихся средней школы. Из опыта работы школы 444 г. Москвы.} Изд. АПН РСФСР, 1963 (тираж 11300) 
 
 Шварцбурд С.И., Монахов В.М., Ашкинузе В.Г. (составители) {\it Обучение в математических школах} (сборник статей). Просвещение, 1965 (тир. 10 000)
  
 Шварцбурд С.И. (составитель) {\it Математический анализ и алгебра} (сборник статей), Просвещение, 1967
 (тир. 30 000)
  
 Н. Я. Виленкин, Р. С. Гутер, С. И. Шварцбурд и др.
{\it  Алгебра: Учеб. пособие для IX-X классов сред. школ с матем. специализацией} - Москва : Просвещение, 1968
(тир. 40 000)
  
 Шварцбурд С.И. {\it Математика и естествознание. Проблемы математической школы.} (сборник статей), Просвещение, 1969 (?)
   
 Шварцбурд С.И. {\it Математика и естествознание. Проблемы математической школы.} (сборник статей), Просвещение, 1970 (?)
  
 Виленкин Н.Я.,  Шварцбурд С.И. {\it Математический анализ. Учебное пособие для школ с математической специализацией.} Просвещение, 1969 (тир 30 000), 1973 (тир 100 000)

 }

{\bf  Шершевский Александр Абрамович (-1973).}
 Из учительской газеты 30.11.1973.
  В 1937г. А.~А.~Шершевский поступил на физико-математический факультет МГУ... 
  В 1964 году Александр Абрамович становится научным сотрудником лаборатории обучения математике АПН СССР. 
  Одновременно с этим он преподает математику в физико-ма\-те\-ма\-ти\-че\-ской школе при МГУ. Являясь членом программной комиссии АН СССР и АПН СССР,
  А.~А.~Шершевский принимает активное участие в создании новой программы по математике и совершенствовании новых учебников...
А.~А.~Шершевский участвовал в составлении сборников серии <<Математическая школа>>, им непосредственно создан факультатив <<Множества и операции над ними>>, 
получивший широкое распространение в нашей стране.

  {\bf  Яглом  Исаак Моисеевич}  (1921—1988), д.ф.м.н, 1965.
Московский педагогический институт, доцент 1956, в 1965-1968 профессор

 {\footnotesize
 Болтянский В.Г. , Яглом И.М.
{\it Геометрия. Учебное пособие для 9 класса средней школы.} - М.: Учпедгиз. 1963 
М.: Учпедгиз. 1963 (тираж 2 300 000), 1964 (тираж 2 000 000), издана также в Киеве, 1963, тираж 230 000.
 [см. комментарии в разделе 3]
 
Также автор книг
 
 Болтянский В. Г., Яглом И. М. {\it Выпуклые фигуры}. ГИТТЛ, 1951 (тир. 25 000) 
 
Яглом И. М., Яглом А. М. {\it Неэлементарные задачи в элементарном изложении.} Гостехиздат, 1954.
(тир. 35 000)
 
Яглом И. М. {\it Геометрические преобразования. Том 1-2.} Гостехиздат 1955—1956 (тир. 25 000+15 000)
 
Головина Л. И., Яглом И. М. {\it Индукция в геометрии.}  Гостехиздат 1956(тир. 35 000), Физматгиз,  1961 (тир. 35 000). 
 
Яглом И. М. {\it  Теория информации.} Знание, 1961 (тир. 33 000)
 
В.Г. Болтянский, И.М. Яглом. {\it Векторы в курсе геометрии средней школы,}  Учпедгиз, 1962 (тир. 25 000) 
 
Яглом И. М, Ашкинузе В.Г. {\it Идеи и методы аффинной и проективной геометрии.
Учебное пособие для педагогических институтов в трех частях,} Учпедгиз, 1962 (тир. 15 000)
 
В.Г. Болтянский, И.М. Яглом. {\it Преобразования. Векторы.} (для учителей) Просвещение 1964 (тир. 56 000) 
 
Яглом И. М. {\it Комплексные числа и их применение в геометрии.} Физматгиз, 1963 (тир. 43 000)
 
Соминский В.Г.,  Л. И. Головина,  Яглом И. М. {\it О математической индукции,} М. Наука, 1967 (тир. 75 000)
 
Яглом И. М. {\it Как разрезать квадрат?} Наука, 1968. (тир. 125 000)
 
Яглом И. М. {\it Герман Вейль.} М.: Знание, 1967. (тир.  42 100)
 
Яглом И. М. {\it Необыкновенная алгебра.} М. Наука, 1968. (тир. 240 000)
 
Яглом И. М. {\it Геометрия точек и геометрия прямых.} М., Знание, 1968. (тир 31 100)

Яглом И. М. {\it Принцип относительности Галилея и неевклидова геометрия.} Наука, 1969. (тир. 50 000)
 
Яглом И. М. {\it О комбинаторной геометрии,} Знание, 1971 (тир. 46 510)
 
 Яглом И. М. 
{\it  Элементарная геометрия прежде и теперь}  Москва : Знание, 1972.

 }
 
А также соавтор известного зубодробительного задачника  для матшкол 
 
 {\footnotesize 
Д.О. Шклярский, Г.М. Адельсон-Вельский, Н.Н. Ченцов, А.М. Яглом, И.М. Яглом. 
{\it Избранные задачи и теоремы элементарной математики. Часть 1. Арифметика и алгебра.} --- М.-Л.: ГТТИ, 1950. --- 296 с. (тир 25 000)
 
 Шклярский Д. О., Ченцов Н. Н.,  Яглом И. М. {\it Избранные задачи и теоремы элементарной математики: Арифметика и алгебра,}  М.: ГТТИ, 1954 (тир. 25 000), М.: Физматгиз, 1959  (тир. 30 000),  М.: Наука. Физматлит, 1965.
 (тир. 50 000),  М.: Наука. Физматлит, 1976(тир 100 000)
 
 Шклярский Д. О., Ченцов Н. Н., Яглом И. М. {\it  Избранные задачи и теоремы элементарной математики: Геометрия (планиметрия)}.  М.: ГТТИ, 1952 (50 000 экз.), М.: Наука, 1967. (тир. 25 000)
 
 Шклярский Д. О., Ченцов Н. Н.,  Яглом И. М. {\it  Избранные задачи и теоремы элементарной математики: Геометрия (стереометрия)}. М.: ГТТИ, 1954
(тир. 50000 экз.)

 Шклярский Д. О.,Ченцов Н. Н., Яглом И. М. {\it Геометрические неравенства и задачи на максимум и минимум.} 1970. 
 М., Наука, 1970 (тир. 75000)

}

 {\bf 
 Ягодовский Михаил Ильич}  (1919?-?),  защитил кандидатскую диссертацию в Курске в 1968, доцент Курского пединститута
  
  {\footnotesize
  Принцев Н.А.,   Ягодовский М.И., Зотов
  {\it 
 Повышение эффективности преподавания математики в общеобразовательной школе: Метод. рекомендации} 
 // Курский ин-т усовершенствования учителей. Курское обл. отд-ние Пед. о-ва РСФСР. - Курск, 1962 (тир. 2 000)

 Клопский В.М., Ягодовский М.И  {\it Геометрия, 9-10 классы,} Просвещение, 1966 (тираж 15 000)
     
  Принцев Н.А.,   Ягодовский М.И 
  Арифметика. {\it Учебник для 5-6 х классов сред. школы},  Просвещение,  1966. (тир. 20 000)
   
     
}
    
    \medskip
    
    Наконец, среди соратников Маркушевича и Колмогорова были и профессиональные математики,
    (см., хотя бы [ММЧ1978]), напрямую со школьной педагогикой не связанные. 
    
   {\bf  Ляпунов Алексей Андреевич}
    
    {\bf Гнеденко Борис Владимирович}
    
    {\footnotesize
    Гнеденко Б.В. {\it Краткие беседы о зарождении и развитии математики,} Изд АПН РСФСР, 1946 (тир. 25 000)
     
    Гнеденко Б.В. {\it Как математика изучает случайные явления.} Львов, Изд АН УСССР, 1947 (тир. 8 000)
     
    Гнеденко Б.В. {\it  Выдающийся русский ученый М.В. Остроградский.} М.Знание, 1952 (тир. 69 000)
   
    Гнеденко Б.В.  
{\it Языком математики,} Знание, 1962 (тир.40 000)
     
    Гнеденко Б.В. {\it Беседы о математической статистике.} М. Знание, 1968 (тир. 60 000).

}

   {\bf Соболев Сергей Львович}
   
   \smallskip

    Вот так. Тоже люди знаменитые и заслуженные. Кто б спорил...О том и статья.

   \bigskip

   \begin{center}
   \bf \large
 B. Победители конкурса учебников 1962-1964гг 
\end{center}
 
 См. [ГП1965]. Жирным выделены авторы будущих  реформистских учебников.
 
  \medskip
  
 АРИФМЕТИКА
  
 Вторая премия: Н.~А.~Принцев, {\bf М.~И.~Ягодовский}
  
 Поощрительная премия: С.~А.~Пономарев, П.~В.~Стратилатова, Н.~И.~Сырнева
  
  Поощрительная премия: С.~Ф.~Моисеев (учитель)
   
   \medskip
   
  АЛГЕБРА ДЛЯ 8-ЛЕТНЕЙ ШКОЛЫ
   
   Поощрительная премия: Н.~А.~Принцев, П.~А.~Ларичев
    
    Поощрительная премия: М.~Ф.~Клюквин (учитель)
   
   \medskip 
   
  АЛГЕБРА И НАЧАЛА АНАЛИЗА ДЛЯ СРЕДНЕЙ ШКОЛЫ
   
 Первая премия: Е.~С.~Кочеткова, Е.~С.~Кочетков
  
 Вторая премия: Н.~И.~Худобина, А.~И.~Худобина, М.~Ф.~Шуршалов (учителя)
  
  Поощрительная премия: {\bf А.~И.~Маркушевич, Р.~С.~Черкасов,} К.~П.~Сикорский
   
   Поощрительная премия:   В.~Е.~Андреев, {\bf  Б.~Е.~Вейц, И.~Т.~Демидов}
    
    \medskip 
    
   ГЕОМЕТРИЯ ДЛЯ ВОСЬМИЛЕТНЕЙ ШКОЛЫ
    
 Вторая премия: {\bf  А.~Ф.~Семенович, Ф.~Ф.~Нагибин, Р.~С.~Черкасов}
  
  Поощрительная премия: К.~С.Барыбин
   
   Поощрительная премия: П.~Я.~Великина
    
    \medskip 
   
   ГЕОМЕТРИЯ ДЛЯ СРЕДНЕЙ ШКОЛЫ
    
    Поощрительная премия: К.~С.~Барыбин
     
     Поощрительная премия: {\bf  В.М.Клопский,  М.И.Ягодовский}
      
     \medskip 
     
     Черкасов и Ягодовский присутствуют в списке дважды.
      
     Председатель комиссии - Б.~В.~Гнеденко. Председатели предметных комиссий: арифметика - В.~И.~Левин, алгебра - А.~Г.~Курош, 
     геометрия - Н.~Ф.~Четверухин.
    
     Из этих учебников в середине 60х пошел в дело учебник Кочетковых как массовый.
     Были изданы как  пробные учебники Семенович-Нагибин-Черкасов, Клопский-Ягодовский
     (позже в реформистских учебниках [КСНЧ],  [КСЯ1977] мы видим тех же авторов, но это были совсем другие учебники),
     Маркушевич-Черкасов-Сикорский,
      Принцев-Ягодовский, Принцев-Ларичев, Клюквин, Барыбин. Книга Худобин-Худобина-Шуршалов была издана как сборник задач.
     
      \medskip 
     
     Академиков АН СССР среди авторов учебников пока нет...
     
     Учебник Кочетковых потом вышел под редакцией О.~Н.~Головина, работавшего на кафедре у Куроша.

 Напомню, что учебник Болтянского и Яглома пошел в дело в 1963-64гг без конкурса.
 
\begin{center}\bf\large
 C. Математики - действительные члены и члены-корреспонденты АПН РСФСР и потом АПН  СССР в 1960-1980 гг 
 \end{center}
 
См. [Коля2001].
  
 Александров Павел Сергеевич (1896-1982), акад. 1945 года, друг Колмогорова.
    
 Колмогоров Андрей Николаевич (1903-1987), акад. с 1965.  
  
 Маркушевич Алексей Иванович (1908-1979), чл-корр. с 1945, акад. с 1950,
  вице-президент 1950-1958, 1964-1975.
  
 Четверухин Николай Федорович (1891-1974), чл.-корр с 1945, акад. с 1955.
  
 Андронов Иван Козьмич (1894-1975), чл.-корр. с 1957. Руководил разработкой реформы образования в 1-3 
 классах [БМ1975]
 (это была автономная часть реформистского  проекта, не вызвавшая таких нареканий, как его основная часть).
 
 Болтянский Владимир Григорьевич (1925), чл.-корр. с 1965.
  
 Брадис Владимир Модестович (1890-1975), чл.-корр с 1955.
 
 Бровиков Иван Семенович (1916-1971),  чл.-корр с 1965, <<Сторонник изучения в школе элементов теории вероятностей и математической статистики>>
 [Кол2001]
 
  Верченко Иван Яковлевич (1907-1996), чл.-корр с 1968 (ученик Колмогорова).
  
 Ларичев Павел Афанасьевич (1892-1963), чл.-корр с 1950.
  
 Шварцбурд Семен Ицкович (1918 - 1996), чл.-корр с 1968. 
 
 \sm
  
 Всего 11 человек. Четверухин, Брадис и Ларичев 1890-1892гг. рождения, к появлению Колмогорова
 Ларичев умер, двум другим было 75 лет (об их отношении к реформе мне ничего не известно). Все остальные, так или иначе,
 связаны с Реформой или Колмогоровым.
 
 \sm
 
 Стоит упомянуть членов АПН РСФСР, умерших  до 1960г., то есть до начала большого реформирования:
 
 Арнольд Игорь Владимирович (1900-1948), чл.кор. с 1947;
 
 Хинчин Александр Яковлевич (1894-1959), акад. с 1944;
 
 Гончаров Василий Леонидович%
 \footnote{Кстати, крупный математик, почему-то выпавший из современных святцев.}(1896-1955),
  чл.кор с 1944;
 
 Перепелкин Дмитрий Иванович (1900-1954), чл.кор с 1950.
 
\label{rus-end}

\section*{Литература (Bibliography%
	\footnote{according Russian alphabet order.})} 

\label{bib-begin}

По ссылкам ниже открывается много учебников, 1888-1988. Почти все из них (за парой исключений) являются памятниками педагогической мысли, 
хотя во многих случаях и заблудшей. Кое-где я после ссылки привожу характерные и интересные цитаты из статей, а также свои комментарии, там, 
где это необходимо. {\bf В  почти всех случаях соответствующие тексты снабжены работающими link'ами.} {\bf Usually references are equipped with working links.}
 
 \bigskip

\hangindent=0.3cm \noindent
[Абр1988] 	А. М. Абрамов,
\href{http://www.mathnet.ru/php/archive.phtml?wshow=paper&jrnid=rm&paperid=2043&option_lang=rus}
{\it О педагогическом наследии А.Н. Колмогорова}, УМН, 43:6(264) (1988), 39-74.
English translation:
A. M. Abramov, {\it A.N. Kolmogorov's pedagogic legacy}, Russian Math. Surveys, 43:6 (1988), 45-88.

\hangindent=0.3cm
{\footnotesize [Одна из самых последовательных апологий Реформы. Цитата:
	\newline
{\it Переходу школы на новые пособия математики предшествовало экспериментальное обучение. В конце 70-х годов иногда высказывались мнения
об отсутствии практической проверки введенных в ходе реформы учебников. На самом деле эксперимент проводился ряд лет;
непосредственная проверка осуществлялась во всех школах четырех экспериментальных районов
(Суздальский район Владимирской области, Тосненский — Ленинградской,
г. Севастополь, Белоярский район Свердловской области). При этом проверке подлежали различные варианты учебников: 
пособия по алгебре и началам анализа Б. Е. Вейца, И. Т. Демидова и Кочетковых, пособия по геометрии для VI — VIII классов А. Ф. Семеновича,
Р. С. Черкасова, Ф. Ф. Нагибина и коллектива под руководством В. Г. Болтянского; по стереометрии
конкурировали учебники К. С. Барыбина и коллектива под руководством
3. А. Скопеца. Окончательное решение принималось после сопоставления
результатов эксперимента. Можно обсуждать, насколько удачной оказалась
схема проведения экспериментов, но не следует отвергать их наличия: это
противоречит фактам.}
\newline 
 Никто не обвинял реформаторов в том, что экспериментов вообще не было. 
Возражения состояли
в некорректности постановки экспериментов и анализа их результатов,  в частности, 
в малом промежутке времени между началом экспериментов и запуском
всеобщей реформы. Подробное обсуждение см. 
в п.\ref{ss:experiment}, \ref{ss:tosno}.].

}

 
 \hangindent=0.3cm \noindent
  [Абр2003]
 Абрамов А. М.
 \newline
  \href{http://www.mathedu.ru/lib/books/abramov_o_polozhenii_s_matematicheskim_obrazovaniem_2003/#1}
 {\it О положении с математическим образованием в средней школе (1978—2003).} Фазис, 2003
 
 \hangindent=0.3cm
 {\footnotesize{Другой вариант сентенции, цитированной в связи со ссылкой [Абр1988]:
 		\newline
 {\it Обсуждение различных промежуточных вариантов прогамм происходило чрезвычайно широко в период с 1965 по 1968 годы. В последующие годы
 учебники рецензировались и обсуждались во всех республиках, областях и краях. Эксперимент в 4 районах
 проходил много лет и позволил выявить победителей (по каждому предмету и на каждой ступени конкурировали два учебника).
 \newline 
 Собственно очевидно, что без оснований и одобрений серьезные решения о массовом переходе школы на новые программы и учебники
 в советское время не могли быть приняты. Проблема лишь в том, что в 1978 году резко изменилась точка зрения
 внутри отделения} [математики].
\newline
Разумеется, между 1972 и 1978г. не могла случиться ничего такого, что могло бы изменить 
мнение Отделения математики АН СССР.
 Еще цитата: 
 \newline 
 {\it Одно из объяснений интереса вузов к школе носит вполне прагматический характер. В связи с появлением выпускников, обучающихся по новым программам
 и учебникам, возникла острая проблема: Как экзаменовать абитуриентов? Строго говоря, прошедшая перестройка школьного курса
 требовала и перестройки системы вступительных экзаменов, перемен в системе экзаменационных заданий. Но существовал и иной вариант
 -- резкая критика нововведений и сохранение устоявшейся системы.}
\newline
 Не совсем только понятно, был ли выпускной школьный экзамен по 
 геометрии в тот год отменен  потому, что не придумали, как его проводить?]}

}
 
 \hangindent=0.3cm \noindent
 [Абр2010] A.M. Abramov,
 \newline 
  \href{http://mat.univie.ac.at/~neretin/misc/reform/Abramov-Karp.pdf}{\it Toward 
 a History of Mathematics Education Reform in Soviet Schools (1960s–1980s).}
In {\it Russian Mathematics Education.
History and World Significance}, Wold Scientific, 2010, pp.87-140.

\hangindent=0.3cm
{\footnotesize 
[Другая апология реформы, на этот раз  по-английски.]

}
 
 \hangindent=0.3cm \noindent
[Абр2016] Абрамов А.М
\newline 
 \href{https://math.ru/lib/files/pdf/Kolmogorov-biblio.pdf}{\it Андрей Николаевич Колмогоров
Полная библиография  его трудов}
{\it и список публикаций, ему посвящённых.} Москва, МЦНМО, 2016
 
 \hangindent=0.3cm \noindent
  [АА2002]
 Авдеев Ф.С., Авдеев Т.К. \href{http://mat.univie.ac.at/~neretin/misc/reform/Avdeev-Kiselev.djvu}
 {\it Андрей Петрович Киселев,}  Издательство Орловской государственной телевещательной компании, 2002.

  \hangindent=0.3cm \noindent
[Але1980] Александров А.Д. \href{http://mat.univie.ac.at/~neretin/misc/reform/Alexandrov1.pdf}{\it О геометрии.}
 Математика в школе. 1980. № 3. С. 56–62.
 
 \hangindent=0.3cm \noindent
 [Але1980-1] Александров А.Д. \href{http://mat.univie.ac.at/~neretin/misc/reform/Alexandrov2.pdf}
 {\it О состоянии школьной математики.}
 Доклад на заседании ученого совета 
 института математики СО АН СССР 25.12.1980
 
 \hangindent=0.3cm
{\footnotesize [Далее была принята резолюция СО АН СССР [РезСО1980].]}
   
   \hangindent=0.3cm \noindent
[Але1981] Александров А.Д.  \href{http://mat.univie.ac.at/~neretin/misc/reform/stih.html}{\it Лев на ниве просвещения.}

\hangindent=0.3cm
{\footnotesize
Эта басня, современная реформе,  всплыла много лет спустя, по-видимому, уже после смерти Александра
Даниловича...	
\newline
\it
Однако курс был слишком гадок
\newline
И никому невпроворот:
\newline
Кого от той науки рвeт, кого проносит, —
\newline
Одни шакалы Льва возносят.

}

\hangindent=0.3cm \noindent
 [АВЛ] Александров А.Д.,   Вернер А. Л.,  Рыжик В. И. {\it Геометрия}. 
 У учебника было много клонов, у меня в закромах есть ранние версии:
  \href{http://mat.univie.ac.at/~neretin/misc/reform/Aleksandrov-Verner6(1984).djvu}{\it 6 класс,  (1984)},  
 \href{http://mat.univie.ac.at/~neretin/misc/reform/Aleksandrov-Verner8(1986).djvu}  {\it 8 класс (1986)},  
 \href{http://mat.univie.ac.at/~neretin/misc/reform/Aleksandrov-Verner9(1981).djvu} {\it 9 класс (1981)}, 
 \href{http://mat.univie.ac.at/~neretin/misc/reform/Aleksandrov-Verner10(1982).djvu} {\it 10 класс (1982)}.

\hangindent=0.3cm \noindent
[АГКЛШ1980]
П.~С.~Александров, Б.~В.~Гнеденко, А.~Н.~Колмогоров, М.~А.~Лаврентьев, Б.~В.~Шабат,
    \href{http://www.mathnet.ru/php/archive.phtml?wshow=paper&jrnid=rm&paperid=3577&option_lang=rus}{\it Алексей Иванович Маркушевич (некролог)},
УМН, 35:4(214) (1980), 131–133. English translation: P.~S.~Aleksandrov, B.~V.~Gnedenko, A.~N.~Kolmogorov, M.~A.~Lavrent'ev, B.~V.~Shabat, 
\newline 
\href{http://mr.crossref.org/iPage?doi=10.1070%2FRM1980v035n04ABEH001863}{\it Aleksei Ivanovich Markushevich (obituary)}, 
Russian Math. Surveys, 35:4 (1980), 153-155.

\hangindent=0.3cm \noindent
[АЛМШ1978]
П. С. Александров, М. А. Лаврентьев, Д.~Е.~Меньшов, Б.~В.~Шабат, 
\newline 
{\it \href{http://www.mathnet.ru/php/archive.phtml?wshow=paper&jrnid=rm&paperid=3511&option_lang=rus}
{Алексей Иванович Маркушевич (к семидесятилетию со дня рождения)}}, УМН, 33:4(202) (1978), 229–235;

\hangindent=0.3cm \noindent
  [Анд1954] Андронов И. К.
{\it Арифметика натуральных чисел.} Экспериментальное пособие. Учпедгиз, 1954

\hangindent=0.3cm
{\footnotesize [пособие по арифметике (видимо, все же не для первого класса) начиналось с понятия множества и взаимно однозначного соответствия,
По-видимому, это был первый эксперимент в этом роде.
См. разгромную рецензию Виленкина и Яглома в   [ВЯ1955]. На основе этого пособия был издан экспериментальный учебник  [АБ1962],
из которого были, в частности, изъяты раздражающие элементы.]

}

\hangindent=0.3cm \noindent
[Анд1967]
Андронов И.К.
\newline
 \href{http://mat.univie.ac.at/~neretin/misc/reform/Andronov.djvu}{\it Полвека развития школьного математического образования в СССР.}
- М.: Просвещение, 1967. - 180 с.
 
 \hangindent=0.3cm \noindent
  [АБ1962]
  Андронов И.К, Брадис В.М. {\it Арифметика.  Пособие для средней школы}, Учпедгиз, 1957
   

  \hangindent=0.3cm \noindent
[Арн2000] Арнольд В.И.  \href{http://scepsis.net/library/id_649.html}{\it Нужна ли в школе математика?}
Доклад на Всероссийской конференции «Математика и общество. 
Математическое образование на рубеже веков» в Дубне 21 сентября 2000 года

\hangindent=0.3cm
{\footnotesize [сам доклад не про то, см. там ответ на вопрос М.А.Цфасмана]}
 
 \hangindent=0.3cm \noindent
[Арн2002]  Арнольд  В. И.
\href{http://vivovoco.astronet.ru/VV/PAPERS/NATURE/BURBAKI.HTM}
{\it Математическая дуэль вокруг Бурбаки}, Вестник АН СССР, том 72, № 3, с. 245-250 (2002)
 
  \hangindent=0.3cm \noindent
  [Ата1981]
Атанасян Л.С., Бутузов В.Ф., Кадамцев С.Б., Поздняк Э.Г.,
\newline 
\href{http://mat.univie.ac.at/~neretin/misc/reform/Atanasyan(1981).djvu}{\it Геометрия.
Пробный учебник для 6-8 классов.} Просвещение, 
1981
 
 \hangindent=0.3cm \noindent
[АЛС1959]
Ашкинузе B. Г.,  Левин В. И.,  Семушин А. Д.
\newline 
   {\it О    перестройке    
программ  по  математике  в  свете  новых  задач   средней   школы»},  «Математика  в  
школе»,   1959,  No   1,  стр.  40—51

\hangindent=0.3cm \noindent
[АЛС1960]	 Ашкинузе B. Г., Левин  В. И.,  Семушин А. Д., 
\newline
	\href{http://www.mathnet.ru/php/archive.phtml?wshow=paper&jrnid=mp&paperid=646&option_lang=rus}
	{\it Некоторые замечания к проекту программы по математике для средней школы}
	Математика, ее преподавание, приложения и история, 5,  ГИТТЛ, 1960, 127–132
 
 \hangindent=0.3cm \noindent
[Барс] Барсуков А.Н. \href{http://mat.univie.ac.at/~neretin/misc/reform/Barsukov.djvu} {\it Алгебра, 6-8 классы}, Просвещение, 1966
  
  \hangindent=0.3cm \noindent
  [БД1975]  Белый Б. Н., Дербенева К. Ф.
  \newline
\href{http://mat.univie.ac.at/~neretin/misc/reform/BeMa3}
{\it Учебники и учебные пособия для начальной и средней школы, 1917-1972.} cc. 319-332, в книге
Штокало И.З. (ред) История математического образования в СССР. - Киев: Наукова думка, 1975.

  \hangindent=0.3cm \noindent
[БМ1975] Белый Б. Н., Маслова Г. Г.,  Беспамятных Н. Д., {\it Развитие преподавания математики в общеобразовательной средней школе}, 15-92, в книге
Штокало И.З. (ред) История математического образования в СССР. - Киев: Наукова думка, 1975. - 383 с;
\newline 
\href{http://mat.univie.ac.at/~neretin/misc/reform/BeMa1}
{\it 4. Период после введения закона <<Об укреплении связи школы с жизнью} {\it и о дальнейшем развитии народного образования в СССР>> (1958-1965гг.)}, 
\newline 
 \href{http://mat.univie.ac.at/~neretin/misc/reform/BeMa2}
{\it 5. Период реформы школьного образования (1966-1972гг)}

\hangindent=0.3cm
{\footnotesize[К сожалению, в файле по ссылке отсутствуют некоторые страницы. Цитаты: 
	\newline
 {\it Анализ программы 1968г. показывает, что ее авторы при опеределении содержания школьного курса не пошли по пути резкой
модернизации математического образования. При составлении программы учитывался опыт, накопленный в массовых и экспериментальных школах, 
в сочетании с оригинальным решением ряда методических проблем...
\newline
Таким образом разработанная под руководством А.Н.Колмогорова программа отличается большой продуманностью.
Ее авторы бережно отнеслись к прогрессивному наследию прошлого и без лишнего увлечения внесли в школьный курс то новое,
что характерно для современных тенденций модернизации школьного математического образования.}]

}

 \hangindent=0.3cm \noindent
[БВЯ1959]  Болтянский В. Г.,  Виленкин Н. Я.,  Яглом И. М., 
\newline 
\href{http://www.mathnet.ru/php/archive.phtml?wshow=paper&jrnid=mp&paperid=529&option_lang=rus}
{\it О содержании курса математики в средней школе},
Математика, ее преподавание, приложения и история,
 4, 1959, 131–143
  
  \hangindent=0.3cm \noindent
  [БВС1972]
  Болтянский В.Г., Волович М.Б., Семушин А.Д.
{\it Геометрия 6-8.} Экспериментальное учебное пособие. - М.: Педагогика.  Эти экспериментальные учебники
издавались с 1972 года, 6, 7, 8 класс. Полный комплект был издан в 1979 году,
он \href{http://mat.univie.ac.at/~neretin/misc/reform/Bolt-Vol.djvu}{\it вот тут}.
    
    \hangindent=0.3cm \noindent
 [БЛ1973]
 Болтянский В.Г., Левитас Г.Г. 
 \newline 
  \href{http://mat.univie.ac.at/~neretin/misc/reform/Bolt1973.djvu}{\it Математика атакует родителей.}
 Педагогика, 1973; второе издание, 1976.

\hangindent=0.3cm \noindent
 [БЯ1963]
 Болтянский В.Г. , Яглом И.М.
 \newline 
\href{http://mat.univie.ac.at/~neretin/misc/reform/Bolt1963.djvu}{\it Геометрия. Учебное пособие для 9 класса средней школы.} - М.: Учпедгиз. 1963
  
  \hangindent=0.3cm \noindent
[Бор1914] Э. Борель,
\newline \href{http://www.mathnet.ru/php/archive.phtml?wshow=paper&jrnid=mp&paperid=493&option_lang=rus}
{\it Как согласовать преподавание в средней школе с прогрессом науки.} (Перевод с французского под ред. М.3.Кайнера; 
предисловие Я.С. Дубнова)”, 
Математика, ее преподавание, приложения и история,  3,  1958, 89–100. Translated from
 Е.  Borel,  {\it L'adaption   de   renseignemeint   
secondaire   aux  progres    de    la   Science},  L'Ens.    math.    16,   {\bf  1914},       198 -210.

\hangindent=0.3cm \noindent
 [БЛ1959]
 И. Н. Бронштейн, А. М. Лопшиц, 
 \newline 
 \href{http://www.mathnet.ru/php/archive.phtml?wshow=paper&jrnid=mp&paperid=531&option_lang=rus}
 {\it Реплики:
 Не изгонять из школы идей аксиоматического метода},
 Математика, ее преподавание, приложения и история,
 Матем. просв., сер. 2, 4, 1959, 151–152

\hangindent=0.3cm \noindent 
[ВД1969]
Вейц Б.Е., Демидов И.Т.  \href{http://mat.univie.ac.at/~neretin/misc/reform/Veits-Demidov.djvu}
{\it Алгебра и начала анализа, 9 класс,} 
Пробный учебник под редакцией Колмогорова
Просвещение, 1969, Был также учебник для 10ого класса, Просвещение, 1971.

\hangindent=0.3cm
{\footnotesize Это учебник, который предполагалось пустить в дело
	в 1975г. В последний час (см. п. \ref{ss:experiment}) он был преобразован в [КВДШ].
	\newline 
	{\it После этого предварительного рассмотрения нетрудно понять следующее
	определение предела...} 
\newline
Дальше идет (правильное) определение предела на языке $\epsilon$
и $\delta$. Я лично видел очень мало студентов за рамками очень узкого круга,
которым это определение было бы легко понять, и твердо могу сказать, что на основе
объяснений, предлагаемых авторами, его не поймёт НИКТО. Не ясно, понимали ли
его сами авторы, вот что они пишут про  максимумы:
\newline	
{\it Точка $x_0$ из области определения функции называется точкой максимума этой функции,
если найдется такая $\delta$-окрестность $(x_0-\delta, x_0+\delta)$ этой точки,
 что в интервале $(x_0-\delta,x_0)$ функция возрастает, а в интервале
 $(x_0,x_0-\delta)$ убывает...} 
\newline
Пример в другом роде. Задача. {\it Материальная точка совершает прямолинейное движение 
	по закону 
	\newline
	$s(t)=5t+2t^2-2/3t^3$
	\newline
	где $s(t)$- путь в метрах, а $t$ -- время в секундах. В какой момент времени $t$
	скорость движения будет наибольшей и какова величина этой скорости?}
\newline
Надо сказать, что текст такого уровня не улучшится, даже если исправить в нем все ошибки...
Однако на этой книге стоял гриф под <<редакцией Колмогорова>>, и эта книга, надо думать,
была многократно положительно отрецензирована (хотя любой грамотный человек
должен был бы <<схватиться за голову>>).

}

\hangindent=0.3cm \noindent 
[Вер2012] Вернер А.Л. \href{https://cyberleninka.ru/article/v/a-d-aleksandrov-i-shkolnyy-kurs-geometrii}{\it  А.Д. Александров и школьный курс геометрии}.
Математические
структуры и моделирование
2012, вып. 25, с. 18–38.

\hangindent=0.3cm
{\footnotesize{[{\it Я помню как в июне 1967 года, приехав в Петрозаводск на Всесоюзный симпозиум 
по геометрии <<в целом>>, 
А.В.Погорелов гордо сказал мне: <<Я написал курс элементраной геометрии. Я ввел в нем аксиомы расстояния.
Меня похвалил Колмогоров.>>
\newline
да это были классы 9–10, теперь это классы 10–11).
Подбирая авторские коллективы для различных учебников математики,
А.Н. Колмогоров ездил по педагогическим вузам страны и встречался с математиками. 
Приезжал он и в Герценовский институт}[Ленинградский пединститут], {\it  и я помню, как в кабинете
ректора мы встречались с А.Н. Колмогоровым, и речь шла о реформе школьного курса математики. Наверное, наши взгляды не подходили А.Н. Колмогорову:
никого из герценовцев он в свою команду не взял. Учебник по геометрии для
старших классов А.Н. Колмогоров поручил писать профессору Ярославского
пединститута З.А. Скопецу и доцентам Курского пединститута В.М. Клопскому
и М.И. Ягодовскому.}
\newline
Приводится письмо А.Д.Александрова от 10.05.1979:
{\it Мне прислали из Министерства рукопись нового издания (нового варианта)
пособия} [КСЯ1977]. {\it Министр написал мне предложение стать научным редактором.
Но по ознакомлении с сочинением, я пришёл к выводу, что редактировать
его — напрасный и невозможный труд; нужно — и это проще — переписать
сочинение заново. Вот я и хочу это сделать и притом совершенно срочно....
Революция в средней школе — злодейство. Одно уже было.
Второго допустить ни в коем случае нельзя. Виноградово-Тихоновская революция или
контрреволюция может быть ещё хуже Колмогоровской. Надо не дать им ходу.
А для этого надо захватить инициативу, т.е. надо взяться за улучшение дела
реально, без широковещательных деклараций, без лишней ругани и пр.
}]}

}

\hangindent=0.3cm \noindent
[ВВ2014] Вечтомов  Е. М., Варанкина В. И., 
\newline 
 \href{http://mat.univie.ac.at/~neretin/misc/reform/Nagibin.pdf}
{\it Профессор Федор Федорович Нагибин,}
Вестник Вятского государственного гуманитарного университета,  2014, 5, 170-176

\hangindent=0.3cm
{\footnotesize[Цитата:{\it 
На Всесоюзном конкурсе учебников геометрии для 8-летней школы в 1964 г. учебник
группы авторов Ф. Ф. Нагибина, А. Ф. Семеновича, Р. С. Черкасова <<Геометрия. Учебник для 6–8 классов>>
был признан лучшим, получил вторую премию (первая ни одному учебнику присуждена не была)
и издан в 1967 г. издательством «Просвещение» объемом 384 с. Академик
А. Н. Колмогоров, возглавлявший комиссию АН СССР по реформе математического образования,
предложил авторскому коллективу этого учебника создать под его руководством 
современный учебник геометрии для 8-летней школы. С этого момента началось сотрудничество
Федора Федоровича Нагибина с выдающимся математиком А. Н. Колмогоровым. Работа
по созданию учебника была выполнена.}]

}
 
\hangindent=0.3cm \noindent 
  [Вил1964] Виленкин Н.Я.
  \newline 
    \href{http://mat.univie.ac.at/~neretin/misc/reform/Vilenkin1965.ps}
  {\it О некоторых аспектах преподавания математики в младших классах,} Математика в школе, 1965, 1, 19-29
 
 \hangindent=0.3cm \noindent
 [ВНШСЧ]
  Виленкин Н.Я., Нешков К.И., Шварцбурд С.И., Чесноков А.С., Семушин А.Д.,
  Математика, 4-5 класс (под редакцией Маркушевича)
  Учебник, по-видимому, имел разные версии, вот \href{http://mat.univie.ac.at/~neretin/misc/reform/Vilenkin1971.pdf}{\it 5ый класс за 1971 год.}
  \href{http://mat.univie.ac.at/~neretin/misc/reform/[Vilenkin_N.YA.,_Neshkov_K.I._i_dr.]_Matematika._4(z-lib.org).djvu}{\it 4ый класс за 1977 год.}
 
 \hangindent=0.3cm \noindent
 [ВЯ1955]
 Н. Я. Виленкин, И. М. Яглом, 
 \newline 
 \href{http://www.mathnet.ru/php/archive.phtml?wshow=paper&jrnid=rm&paperid=7993&option_lang=rus}
{\it И. К. Андронов, “Арифметика натуральных чисел” (рецензия)}, УМН, 10:2(64) (1955), 225–228

\hangindent=0.3cm
 {\footnotesize [довольно любопытный отсвет, в тексте моей статьи на это ссылок нет]
 
}
 
 \hangindent=0.3cm \noindent
 [Вин2015]
 Э. Б. Винберг,
 \newline \href{http://www.mathnet.ru/php/archive.phtml?wshow=paper&jrnid=mp&paperid=848&option_lang=rus}
 {\it О концепции учебника геометрии А. В. Погорелова}, Матем. просв., сер. 3, 19, Изд-во МЦНМО, М., 2015, 199–205
 
 \hangindent=0.3cm
 {\footnotesize[Чрезвычайно резкая критика учебника Погорелова:
 \newline
 {\it Таким образом, концепция учебника Погорелова полностью несостоятельна. Она приводит к тому, что учебник не только не способен
пробудить интерес к геометрии, но может вызвать её неприятие, особенно
на решающем начальном этапе обучения. Конечно, он сообщает некоторые
полезные сведения (которые, впрочем, можно найти и в справочнике), но он
не решает задач интеллектуального и духовного воспитания учащихся.
Поразительно, что этот учебник в течение столь долгого времени поддерживался и продолжает поддерживаться Министерством просвещения,
а затем Министерством образования (и науки) РФ.}

}

\hangindent=0.3cm \noindent
[В-мин1969] \href{http://mat.univie.ac.at/~neretin/misc/reform/V-min.ps}
{\it В Министерстве просвещения СССР},
Математика в школе, 1969, 5, 18-21.

\hangindent=0.3cm \noindent
[Вним1970] \href{http://mat.univie.ac.at/~neretin/misc/reform/v-pomoshch.ps}
{\it Вниманию учителей математики  четвертых классов} {\it  и руководителей методических объединений},
Математика в школе, 1970, 4, с.5

\hangindent=0.3cm \noindent
[ВПТ]
В. С. Владимиров, Л.С.Понтрягин, А. Н. Тихонов
\newline 
\href{http://mat.univie.ac.at/~neretin/misc/reform/VPT.pdf}{\it О школьном математическом образовании.}
 Математика в школе. — 1979. N 3. — С. 12–14. 

\hangindent=0.3cm \noindent
[ВММО1937]
\href{http://www.mathnet.ru/php/archive.phtml?wshow=paper&jrnid=mp&paperid=809&option_lang=rus} 
{\it В Московском математическом обществе}, Сборник статей по элементарной и началам высшей математики, Матем. просв., сер. 1, 13, 1938, 69

\hangindent=0.3cm \noindent
[ГЛШ1958]	А. О. Гельфонд, А. Ф. Леонтьев, Б. В. Шабат,
\newline 
\href{http://www.mathnet.ru/php/archive.phtml?wshow=paper&jrnid=rm&paperid=7517&option_lang=rus}
{\it Алексей Иванович Маркушевич (к пятидесятилетию со дня рождения)}, УМН, 13:6(84) (1958), 213–220
  
  \hangindent=0.3cm \noindent
 [ГКР1973]
Галкина М. С., Колягин Ю. М., Ройтман П. Б.
{\it Уроки геометрии в VII классе }. Пособие для учителей. - М.: «Просвещение», 1973. Первое полугодие.
\href{http://mat.univie.ac.at/~neretin/misc/reform/Galkina-Kolyagin.djvu}{\it  Второе  полугодие.}

\hangindent=0.3cm
{\footnotesize[Из первой части:
	\newline
	 {\it Авторы хотели бы отметить замечания многих учителей (приславших отзывы на аналогичное пособие авторов по VI классу),
что материал, предлагаемый на тот или иной урок, часто слишком велик по объему и не может быть изучен за время, 
отводимое ему на данном уроке. Авторы признают эту критику правильной. Однако рамки программы и учебного плана 
исключают возможность дать иное планирование и менее насыщенно распределить учебный материал. 
Основной материал программы должен быть усвоен школьниками в намеченный срок. Поэтому авторы сразу оговариваются,
что данное пособие для VII класса имеет тот же недостаток — объем материала для отдельных уроков, 
может быть, несколько завышен.}
\newline
Вторая часть этой книжки была сдана в набор 27.7.1973 года]

}
  
  \hangindent=0.3cm \noindent
 [Гла] Глаголев Н.А. \href{http://mat.univie.ac.at/~neretin/misc/reform/Glagolev-plan.djvu}{\itПланиметрия,} Учпедгиз, 1954, 
 \href{http://mat.univie.ac.at/~neretin/misc/reform/Glagolev-ster.djvu}{\it Стереометрия,} Учпедгиз, 1948
 
  \hangindent=0.3cm \noindent
  [Глад2009]
  А. В. Гладкий,   \href{http://www.mathnet.ru/php/archive.phtml?wshow=paper&jrnid=mo&paperid=115&option_lang=rus}
  {\it О преподавании алгебры и начал анализа в школе}, Матем. обр., 2009, № 3(51), 7–16 
  
  \hangindent=0.3cm
  {\footnotesize[
  Статья эта интересна и содержательна, но, как положено, во всех провалах повинны российские власти 
  (и, наверно, в параллельной французской реформе они же):
  \newline
{\it   Неудача была обусловлена прежде всего системой управления
образованием и системой подготовки учителей. Важным элементом советской системы упра­вления образованием, 
к 70-м годам окончательно закосневшей, были «стабильные учебники»,
по которым обязаны были преподавать учителя во всей огромной стране. Переходу на новый
стабильный учебник предшествовало его «экспериментальное опробование», но оно было чисто
формальным: в нескольких городах и районах всем учителям данного предмета предписывалось
работать по пробному учебнику и писать отчеты, причем отрицательные отзывы во внима­ние не принимались... 
Геометрия и при традиционном изложении по Киселеву была трудна для школьников,
а новый учебник, возникший в результате сотрудничества А. Н. Колмогорова с талантливыми
педагогами-математиками А. Ф. Семеновичем и Р. С. Черкасовым, был намного труднее. Но это
книга добротная, тщательно продуманная, богатая новыми методическими идеями. В нормаль­ных
условиях этот учебник был бы сначала взят на вооружение небольшим числом учителей, а
со временем на его основе теми же или другими авторами были бы созданы учебные пособия, 
ко­торые получили бы более широкое распространение.} 
\newline 
Ю.Н. В стабильных учебниках, безусловно, есть минусы (а возврат к ним невозможени и нежелателен),
но  1938-1976годы стабильные учебники математики хорошо отработали (и были, кстати, дважды за это время  прореформированы, 1956, 1966).
Экспериментальных учебников после войны издавалось  довольно много разных. С запуском реформы безусловно произошла управленческая ошибка,
решение о запуске асфальтового катка в 1970 году было принято в 1968, когда {\it комплекта} 
новых учебников еще не было.
В том же 1970г., когда когда упомянутый  учебник [КСНЧ]
еще не начал проходить экспериментальную проверку, его запуск в 1972 был уже почти неизбежен.
Однако в 1968 по вопросу реформы было единодушие АН СССР, АПН СССР, разных комиссий, лучших экспертов-педагогов,
знаменитых математиков... То, что все они ошибаются, и что вместо кота из мешка выпрыгнет неукротимая тигра, управленцы все же 
имели право и не предвидеть. С другой стороны, уже были неудачи с геометрическими учебниками Фетисова-Никитина (1956)
и Болтянского-Яглома (1963),  и, казалось бы, именно с математико-педагогической стороны 
следовало бы проявить хотя бы минимальную осторожность.]

}
 
\hangindent=0.3cm \noindent 
  [ГЧ1993]
 Глейзер Г.Д., Черкасов Р.С. 
 \newline 
 \href{http://mat.univie.ac.at/~neretin/misc/reform/rao.doc}
{\it Центр творческих усилий педагогов} {\it [к 50-летию Российской академии образования]}  Математика в школе. 1993. № 5. С. 2–8; 1993. № 6. С. 2–6.

  \hangindent=0.3cm \noindent
 [Гнед1965]
 Б. В. Гнеденко, \href{http://mat.univie.ac.at/~neretin/misc/reform/Gnedenko1965.pdf}{\it О перспективах математического образования},
  Математика в школе, 1965, 6, c.2-11

   \hangindent=0.3cm \noindent
 [Гнед1993]
 Б. В. Гнеденко, \href{http://mat.univie.ac.at/~neretin/misc/reform/gnedenko.pdf}{\it  Учитель и друг}, Колмогоров в воспоминаниях
 учеников, МЦНМО, 2006
 
 \hangindent=0.3cm
  {\footnotesize[Цитаты: {\it Конечно, учебники, написанные коллективом под руководством Колмогорова,
требовали серьезной доработки. ... Учебник мало написать, его необходимо 
 выстрадать и многократно к нему возвращаться. Такой возможности Колмогорову
дано не было. На него свалилась резкая и далеко не всегда справедливая критика.
\newline 
Я считаю критику Л. С. Понтрягина, опубликованную в журнале «Коммунист», 
необъективной и не содержащей никаких положительных предложений
и решений. Поскольку в ту пору на эту статью нельзя было возразить, то иные
точки зрения и оценки не увидели света. В результате наша средняя школа была
дезориентирована и пошла по плохо продуманному пути.
\newline
Вообще, моя точка зрения на школьные реформы состоит в том, что их 
предварительно следует осмыслить всесторонне, проверить экспериментально 
и только затем вводить в широкую практику. Каждая ошибка в такого рода делах
тиражируется в десятках миллионах душ и умов и сказывается, по меньшей мере,
в течение жизни целого поколения.}]

}
    
      \hangindent=0.3cm \noindent 
 [ГП1965]
 Б. В. Гнеденко, И. С. Петраков
 \newline
  \href{http://mat.univie.ac.at/~neretin/misc/reform/Konkurs1965.pdf}{\it Итоги открытого
конкурса на учебники по математике.} Математика в школе, 1965, 2, 4-9

   \hangindent=0.3cm \noindent
 [Гон1955]
 Гончаров В.Л. {\it Начальная алгебра.}  М.: Издательство АПН РСФСР, 1955,
 \href{https://mat.univie.ac.at/~neretin/misc/reform/goncharov_v_l_nachal_naya_algebra.pdf}{2-е издание}.  Москва: Издательство АПН РСФСР, 1960. — 452 с. 
 
 \hangindent=0.3cm
 {\footnotesize [Учебник (по виду качественный) был написан, чтобы несколько реформировать тогдашний курс алгебры. Был издан под грифом
 <<Пособие для учителей>>. См. рецензию на него 
 \href{http://www.mathnet.ru/php/archive.phtml?wshow=paper&jrnid=mp&paperid=448&option_lang=rus}{Г.Б.Гуревич,
 {\it Учебник алгебры В.Л.Гончарова},  Матем. просв., сер. 2, 1, 1957, 243–250}.

}
 
    \hangindent=0.3cm \noindent
[Дуб1960]
Я. С. Дубнов,
\newline 
 \href{http://www.mathnet.ru/php/archive.phtml?wshow=paper&jrnid=mp&paperid=639&option_lang=rus}{\it Содержание
и методы преподавания элементов математического}  {\it  анализа и аналитической геометрии в средней школе}, 
Математика, ее преподавание, приложения и история,  5,  Гос. изд-во технико-теоретической литературы, 1960, 17–55

   \hangindent=0.3cm \noindent 
[Дье1975] Ж. Дьедонне \href{http://mat.univie.ac.at/~neretin/misc/reform/Dieudonne.djvu}{\it Линейная алгебра и элементарная геометрия}, 
М.: Наука, 1975 (перевод с французского Дорофеева под редакцией И.~М.~Яглома).
French original: Dieudonn\'e, Jean
{\it Alg\`ebre lin\'eaire et g\'eom\'etrie \'el\'ementaire.}
Enseignement des Sciences, VIII Hermann, Paris 1964 223 pp
editions: 1964, 1968, 1968, English translation: {\it Linear algebra and geometry.}
 Houghton Mifflin Co., Boston, Mass. 1969 

   \hangindent=0.3cm \noindent
[КС1979]
Канторович Л.В., Соболев С. Л.
\href{http://mat.univie.ac.at/~neretin/misc/reform/KS1979.pdf}
{\it Математика в современной школе}. Математика в школе
— 1979. — No 4. — С. 6–11.

   \hangindent=0.3cm \noindent
[КК1972] \href{http://vivovoco.astronet.ru/VV/PAPERS/BIO/KOLMOGOR/KOLM_KAP.HTM}
{\it Переписка А.Н.Колмогорова и П.Л.Капицы}. Вопросы философии, 1972,7, 16-24.

\hangindent=0.3cm
{\footnotesize [Капица:{\it Я себе представляю задачи специальной школы по сравнению с обычной аналогично тем,
которые преследует клиника по сравнению с больницами.
\newline 
Клиника изучает и отрабатывает новые методы диагностики и лечения и для этого имеет наиболее квалифицированный персонал,
и ее задача - внедрить передовые методы в жизнь и этим поднять уровень медицинского обслуживания больных 
в обычных больницах.
При этом, конечно, клиники должны быть специализированными по определенным видам заболевания. 
Полезность и необходимость такой организации в здравоохранении общепризнаны и не вызывают сомнений. 
То же должно иметь место и при развитии образования.
\newline 
Задача специальных школ - изучать и разрабатывать передовые методы обучения воспитания. 
Спецшколы должны иметь хорошо подобранные кадры преподавателей, образцовую организацию. 
Конечно, такие школы не могут охватывать обучение по всем областям знания
и должны быть специализированы по отдельным дисциплинам, как математика, физика, биология и проч.
\newline 
Тут мы, по-видимому, несколько расходимся с Вами во взглядах. В Вашем письме, характеризуя деятельность Ваших школ,
Вы определяете их значимость по научным успехам Ваших питомцев. Это, конечно, показывает,
что Ваши методы преподавания математики действительно являются более совершенными. Но Вы не говорите о том, 
что Вы предпринимаете, чтобы эти методы обучения распространились более широко, 
и как они влияют на качество преподавания математики в обычных школах. Я считаю,
что повышение уровня преподавания в стране в широких масштабах и должно быть основной задачей спецшкол.
Если это так, то из этого следует, что характер организации этих школ, отбор преподавателей и учеников
должны быть согласованы с этой задачей.
\newline 
Существуют еще специальные школы, в которых основной предмет обучения практически отсутствует в обычных школах. 
Например, это балетные школы, цирковые училища, музыкальные и художественные школы и т. п.
Поскольку такого рода специалисты требуют обучения смолоду и малочисленны, то существование подобных школ
вполне оправдано.}]

}

   \hangindent=0.3cm \noindent 
 [Кат-МГУ]
\href{http://nbmgu.ru/catalogs/alcats/books_rus/}{\it Каталог Библиотеки Московского государственного университета.}
 
    \hangindent=0.3cm \noindent
 [Кат-РГБ]
\href{http://search.rsl.ru/#ff=24.12.2016&s=fdatedesc}{\it Каталог  Российской государственной библиотеки.}
 
    \hangindent=0.3cm \noindent
 [Кат-РНБ]
 \href{http://www.nlr.ru/e-case3/sc2.php/web_gak}{\it  Каталог Российской национальной библиотеки.}

   \hangindent=0.3cm \noindent
 [Кис1892]
Киселёв А.П. \href{http://mat.univie.ac.at/~neretin/misc/reform/Kiselev1892.djvu}{\it Геометрия.}  Лашкевич, Знаменский и К, 1892

   \hangindent=0.3cm \noindent
 [Кис1906]
Киселёв А.П. \href{http://mat.univie.ac.at/~neretin/misc/reform/Kiselev-algebra-1906.djvu}{\it Элементарная алгебра,}
17ое издание,  1906,  Печатня Яковлева, 1906.

   \hangindent=0.3cm \noindent
 [Кис1912]
Киселёв А.П. \href{http://mat.univie.ac.at/~neretin/misc/reform/arifmetika1912-2002.djvu}{\it Систематический курс арифметики.}
 Репринтное издание. Издательство Орловского государственного университета, 2002 
 
    \hangindent=0.3cm \noindent
 [Кис1914]
Киселёв А.П. \href{http://mat.univie.ac.at/~neretin/misc/reform/kis-geom.djvu}{\it Элементарная геометрия,}
изданиe 23е, типография Рябушинского.
  
     \hangindent=0.3cm \noindent
[Кис1925]
Киселёв А.П. \href{http://mat.univie.ac.at/~neretin/misc/reform/Kiselev1925.pdf}{\it Элементы алгебры и анализа}, 1925

   \hangindent=0.3cm \noindent
 [Кис1938-ал]
Киселёв А.П. {\it Алгебра}, Физматлит, 2006, \href{http://mat.univie.ac.at/~neretin/misc/reform/algebra-kiselev1938-2006-1.djvu}{\it Часть 1},
\href{http://mat.univie.ac.at/~neretin/misc/reform/algebra-kiselev1938-2006-2.djvu}{\it Часть 2.} Перепечатка <<Алгебры>> Киселева
1938 года, переработанной А.~Н.~Барсуковым.

   \hangindent=0.3cm \noindent
 [Кис1938-ге]
Киселёв А.П. \href{http://mat.univie.ac.at/~neretin/misc/reform/geometr-plaste-2004.djvu}
{\it Геометрия. Планиметрия. Стереометрия.} Физматлит, 2004. Перепечатка издания 1938 года, переработанного Глаголевым.

   \hangindent=0.3cm \noindent
 [Кис1940]
Киселёв А.П. \href{http://mat.univie.ac.at/~neretin/misc/reform/arifmetika1940-2002.djvu}{\it Арифметика},
2002, перепечатка издания 1940 года, переработанного А.Я.Хинчиным.

   \hangindent=0.3cm \noindent
\hangindent=0.3cm \noindent
 [Кла1973]
Kline, Morris,  \href{http://www.rationalsys.com/mk_johnny.html}{\it Why Johnny Can't Add: The Failure of the New Math.}
New York: St. Martin's Press,
1973

\hangindent=0.3cm \noindent
[КСЯ1977]
 Клопский В. М., Скопец З. А., Ягодовский М. И.
 \newline 
\href{http://mat.univie.ac.at/~neretin/misc/reform/Skopets.djvu}{\it Геометрия, 9-10 классы}. Под редакцией З.А.Скопеца,
Просвещение, 1977

\hangindent=0.3cm
{\footnotesize[Школьный учебник стереометрии, 1976-1983. Определение
многогранника.
\newline
{\it 
\textbf{Простой многогранной поверхностью} называется объединение
конечного числа многоугольников, удовлетворяющее cледующим
условиям:\newline
1) для любых двух вершин этих многоугольников существует
ломаная, составленная из их сторон, для которой взятые вершины
служат концами;
\newline
2) произвольная точка объединения многоугольников либо 
является точкой только одного из данных многоугольников, либо 
принадлежит общей стороне двух и только двух многоугольников, либо
является вершиной только одного многогранного угла, плоскими 
углами которого служат углы данных многоугольников.
\newline
Если каждое ребро многогранной поверхности содержится
в двух ее гранях, то эту многогранную поверхность называют
\textbf{замкнутой}.
\newline
...
Замкнутая многогранная поверхность 
разбивает множество всех не принадлежащих
ей точек пространства на два подмножества.
Для одного из них существуют прямые, 
содержащиеся в этом подмножесте; для 
другого — таких прямых не существует. Первое
из указанных подмножеств называют внешней
областью замкнутой многогранной 
поверхности, а второе — ее внутренней областью.
\newline
\textbf{Определение}. Объединение замкнутой многогранной поверхности и ее внутренней
области называется многогранником.
\newline
Грани, ребра, вершины 
поверхности многогранника называют соответственно гранями, 
ребрами и вершинами многогранника.}
\newline
{\it Как и многоугольники, многогранники могут быть 
выпуклыми (рис. 19) и невыпуклыми (рис. 20). Мы будем изучать только
выпуклые многогранники.}
\newline
Нельзя сказать, что сказанное было исключительно понятным (я еще не все процитировал), и непонятно зачем все
это сказано, если будут рассматриваться только выпуклые многогранники. 
Кстати, поверхность октаэдра, из которой удалены две противоположные
грани, подпадает под определение простой многогранной поверхности.
А если удалить лишь одну грань, то не подпадает. Зато, если удалить 
одну грань из поверхности куба, то снова подпадет.
\newline 
Интересно это сопоставить с определением многогранного угла в том же учебнике:
\newline 
{\it Пусть даны многоугольник $\Phi = ABC\dots$ и точка $S$, не 
принадлежащая его плоскости (рис. 157, а). Объединение всех лучей,
имеющих общее начало $S$ и пересекающих данный многоугольник
$\Phi$ (рис. 157, б), называется многогранным углом.}
\newline
Авторы стараются ради определения невыпуклых многогранников, но в их определении все грани оказываются
выпуклыми (очевидное упражнение).
То есть, согласно З.~А.~Скопецу, пирамида с невыпуклым основанием многогранником не является...
Желающие могут поразмышлять над тем, что в итоге названо <<гранью многогранника>>...
Очевидно, что авторы учебника сами не понимали собственного мудрёного определения.
Фраза: <<многогранник - это тело, ограниченное конечным набором плоских многоугольников}%
\footnote{
У Киселева: {\it Многогранником называется тело, ограниченное со всех сторон плоскостями.}
В учебнике Киселева в редакции Глаголева:
{\it Многогранником называется тело,
ограниченное плоскими многоугольниками}. Ради научности
можно добавить слово <<конечным числом>>, полезность этого добавления неочевидна.}{\footnotesize
>>
и более понятна, и более правильна.
]

}

   \hangindent=0.3cm \noindent
[К-на-1970] \href{http://mat.univie.ac.at/~neretin/misc/reform/nov-progr.ps}{К началу работы по новой программе}. Математика в школе, 1970, 4, с.2-3.

   \hangindent=0.3cm \noindent
 [Колм1965]
А. Н. Колмогоров
\newline 
 \href{http://mat.univie.ac.at/~neretin/misc/reform/Kolmogorov1965.pdf} {\it Геометрические преобразования:
в школьном курсе геометрии}. Математика в школе, 1965, 2, 24-29 (перепечатано в [ММЧ1978])

   \hangindent=0.3cm \noindent
 [Колм1965-1] А. Н. Колмогоров
  \href{http://lyapunov.vixpo.nsu.ru/?int=VIEW&el=1003&templ=VIEW_TYPE}{\it Письмо А.~А.~Ляпунову}, 
видимо, сентябрь-ноябрь 1966. 

\hangindent=0.3cm
{\footnotesize[{\it программы эти составлялись в августе несколько наспех и за их окончательную редакцию ответственны лишь 
несколько сотрудников АПН (Семушин, Фетисов)... В этой обстановке нам пришлось довольно быстро работать
(кроме сотрудников АПН реально работают Болтянский, Виленкин, Яглом, я, Маркушевич)... Я стараюсь сделать формулировки достаточно широкими,
чтобы не слишком стеснять работу авторов учебников. Следующий этап уточнения программ разумно предпринимать уже после того,
как коллективы, работающие над учебниками, представят мотивированные свои пожелания...
Естественно, что программы, подготовляемые к 15 ноября, вновь будут опубликованы для «широкого обсуждения».
Но тем, кто не будет с ними радикально расходиться, важно поскорее переключиться на подготовку учебников, 
так как в программе легко пишется что угодно. Я пока с двумя учителями веду в Болшеве опыт преподавания
начал анализа в 9­ом классе вполне «среднего» состава.}
\newline
Ю.Н.: Интересно, что Виленкина нет среди авторов Программы-1967, [Прог1967]]

}

   \hangindent=0.3cm \noindent
[Колм1967] Колмогоров  А.  Н.
\newline 
\href{http://mat.univie.ac.at/~neretin/misc/reform/kolmogorov-1967.ps}
{Новые программы и некоторые основные вопросы усовершенствования} {\it курса математики в средней школе}.
Математика в школе, 1967, 2, 4-13.

\hangindent=0.3cm \noindent
[Колм1967-1] Колмогоров  А.  Н.
\newline 
 \href{http://mat.univie.ac.at/~neretin/misc/reform/kokoko1967.pdf} {\it Об учебниках на 1967-1968 учебный год. Алгебра и элементарные функции} {\it  Е.С.Кочеткова и Е.С.Кочетковой.}
Математика в школе, 1967, 1, 43-48.

   \hangindent=0.3cm \noindent
[Колм2003] \href{http://mat.univie.ac.at/~neretin/misc/reform/Kolmogorov-pedagogika.pdf}{\it  Список работ А.Н.Колмогорова по педагогике.}
 В книге
{\it Колмогоров. Юбилейное издание в трех книгах. Книга первая. Истина - благо.} Редактор-составитель А.Н.Ширяев.
Москва, Физматлит, 2003 (стр. 286-301)
 
    \hangindent=0.3cm \noindent
 [КМЯ1967] Колмогоров  А.  Н.,  Маркушевич  А.  И.,  Яглом  И.  М.,  
 \newline 
\href{http://mat.univie.ac.at/~neretin/misc/reform/msh1967-1.pdf}{\it Проект  программы  средней  школы  по  математике.}   1967,  No  1, 
стр.  4—23.

   \hangindent=0.3cm \noindent
 [КСНЧ]
 Колмогоров А.Н., Семен\'ович А.Ф., Нагибин Ф.Ф., Черкасов Р.С.
{\it Геометрия. 6-8 класс}. Под редакцией А.Н. Колмогорова.
Учебник имел много  различных версий, начиная с 1970 года (как пробный учебник и с 1972 как общий), 
и разный состав авторов 
(появлялся временно Гусев,
Нагибин в итоге выпал). У меня в закромах есть уже поздние издания  \href{http://mat.univie.ac.at/~neretin/misc/reform/Kolmogorov-Nagibin.djvu}
{\it 8ой класс за 1976 год, 4ое издание},   \href{http://mat.univie.ac.at/~neretin/misc/reform/Kolmogorov-geometry-7-1977.djvu} 
{\it 7ой класс за 1977 год, 6ое издание}, 
\newline
 \href{http://mat.univie.ac.at/~neretin/misc/reform/Kolmogorov-geometry-6-8-1979.djvu}
{\it 6-8 классs за 1979 год}. {\footnotesize [Последние издания , 1981 (тир. 3 500 000, на 1981-1982 учебный год он оставался основным учебником, 
это был десятый заход учебника в 6-ые классы общей школы) и 1982 (тир. 773 000)]}.

   \hangindent=0.3cm \noindent
[КВДШ]
Колмогоров А.Н., Вейц Б.Е., Демидов И.Т., Ивашёв-Мусатов О.С., Шварцбурд С.И.
\newline 
{\it Алгебра и начала анализа}. У этого учебника было много версий, начиная с 1975 года, предварительный вариант 1969 года [ВД1969],
у меня есть
Колмогоров А.Н., Вейц Б.Е., Демидов И.Т., Ивашёв-Мусатов О.С., Шварцбурд С.И.
\href{http://mat.univie.ac.at/~neretin/misc/reform/AA9-1976.djvu}{\it 9 класс}, Просвещение, 1976
 Колмогоров А.Н., Ивашёв-Мусатов О.С., Ивлев Б.М., Шварцбурд С.И.
\href{http://mat.univie.ac.at/~neretin/misc/reform/AA10-1976.djvu}{\it 10 класс}, Просвещение, 1976.

\hangindent=0.3cm
{\footnotesize [В дальнейшем учебник под тем же <<брендом>> продолжал существовать.
С 1980г. среди авторов появляется А.М.Абрамов.
У последнего
издания 1987 при жизни Колмогорова авторы: Колмогоров А.Н., Абрамов А.М., Вейц Б.Е., Демидов И.Т.,
Ивашев-Мусатов О.С., Шварцбурд С.И.,.
С 1990г. авторы: Колмогоров А.Н., Абрамов А.М., Дудницын Ю.П., Ивлев Б.М., Шварцбурд С.И. С этим же составом авторов
учебник издается поныне.]

}

   \hangindent=0.3cm \noindent
 [КЯ1965]  Колмогоров А. Н.,  Яглом  И.  М.
 \newline 
  \href{http://mat.univie.ac.at/~neretin/misc/reform/Kolmogorov-Yaglom.pdf}
{\it О содержании школьного курса математики,} Математика в школе, 1965, 4,  c.53-62.
 
    \hangindent=0.3cm \noindent
 [Колм-проект] \href{http://mat.univie.ac.at/~neretin/misc/reform/content.pdf}
 {\it Проект издания педагогических трудов Колмогорова},
 \href{http://www.math.ru/ank-ped/index.htm} 
 {  МЦНМО}
 
 \hangindent=0.3cm
 {\footnotesize [очень характерный список, в него не вошла основная масса публикаций А.Н., связанных с Реформой, ср. со списком
 педагогических публикаций А.Н. в
 [Колм2003]]

}
 
    \hangindent=0.3cm \noindent
[Коля001]
Колягин Ю.М.
\newline 
\href{http://mat.univie.ac.at/~neretin/misc/reform/Kolyagin.djvu}
{\it Русская школа и математическое образование. Наша гордость и наша боль,} Просвещение, 2001
 
    \hangindent=0.3cm \noindent
 [КС2012]
 Колягин Ю.М., Саввина О.А.
 \newline 
  \href{http://mat.univie.ac.at/~neretin/misc/bunt.pdf}
 {\it Бунт российского министерства и отделения математики АН СССР.} {\it 
 (Материалы по реформе школьного математического образования 1960-1970-х гг.)} - Елец: ЕГУ им. И.А. Бунина, 2012.
 
    \hangindent=0.3cm \noindent
 [Кост]
И.П.Костенко, 
\href{http://www.mathnet.ru/php/archive.phtml?wshow=paper&jrnid=mo&paperid=454&option_lang=rus}{\it Почему надо вернуться к Киселеву?}
 Матем. обр., 2006, № 3(38), 12–17;
 \newline 
\href{http://www.mathnet.ru/links/b49c67355ac8731c73a8c011a8537a74/mo104.pdf}
{\it Корни, ветви и <<ягодки>> реформы-1970.}{  Матем. обр., 2009, № 2(50), 14–23};
\newline
\href{http://www.mathnet.ru/php/archive.phtml?wshow=paper&jrnid=mo&paperid=164&option_lang=rus}
{\it Динамика качества математического образования. Причины деградации.}
Матем. обр., 2011, № 2(58), 2–13;  
\newline	
\href{http://www.mathnet.ru/php/archive.phtml?wshow=paper&jrnid=mo&paperid=220&option_lang=rus}{\it 1918 – 1930 гг.
 Первая коренная реформа русской школы.} Матем. обр., 2012, № 4(64), 2–10;
\newline
\href{http://www.mathnet.ru/php/archive.phtml?wshow=paper&jrnid=mo&paperid=515&option_lang=rus}{\it 1930–1956 гг. Возрождение и рост русской школы.}
 Матем. обр., 2013, № 1-2(65-66), 14–36;
\newline
 \href{http://www.mathnet.ru/php/archive.phtml?wshow=paper&jrnid=mo&paperid=28&option_lang=rus}{\it 
 1956–1965 гг. Подготовка второй <<коренной>> реформы советской школы: <<перестройка>>
  программ и <<научное>> обоснование ложных идей.}
  Матем. обр., 2014, № 2(70), 2–17;
 \newline \href{http://www.mathnet.ru/php/archive.phtml?wshow=paper&jrnid=mo&paperid=33&option_lang=rus}
 {\it 1965 – 1970 гг. Организационная подготовка реформы-70: МП, АПН, кадры, программы, учебники.} Матем. обр., 2014, № 3(71), 2–18
\newline
\href{http://www.mathnet.ru/php/archive.phtml?wshow=paper&jrnid=mo&paperid=8&option_lang=rus}{\it  1970–1986 гг. Реализация реформы-70,
удержание её результатов.}
Матем. обр., 2015, № 2(74), 2–17;
\newline
\href{http://www.mathnet.ru/php/archive.phtml?wshow=paper&jrnid=mo&paperid=63&option_lang=rus}{\it Уроки <<ВТУ-реформы>>.}
Матем. обр., 2015, № 4(76), 2–21

\hangindent=0.3cm 
{\footnotesize[Стоит быть очень осторожным в отношении приводимых в этих статьях фактам. Малая толика того, что можно сказать есть в [Шевк];
кое-что было сказано выше в разделе 2 настоящей статьи. Замечу еще, что <<объективная>> оценка состояния школы с помощью оценок вступительных экзаменов 
(применяемая Костенко, и не им первым) по многим
причинам некорректна: процент тех или иных оценок (2,3,4,5) выставленных на экзамене обуславливается потребностями 
экзаменационных комиссий, а не объективным уровнем поступающих, кроме того, сложность вступительных экзаменов в 60-70е годы возрастала,
 одновременно увеличивалось количество выпускников старшей школы, а также число мест в вузах.
Открывающиеся link'и настоящей библиографии в принципе дают возможности для проверки самых разнообразных заявлений.
Литературные <<зацепки>> в этих статьях интересны и  широко использовались в настоящей работе.
Есть также книжная версия этих статей: Костенко И.П. 
\newline
 \href{https://russianclassicalschool.ru/pdf/kostenko-mono.pdf}{\it
Проблема качества математического образования в свете исторической ретроспективы} 2-е изд., доп. — М.: РГУПС, 2013]

}
 
    \hangindent=0.3cm \noindent 
  [КК1968]
  Кочетков Е.С., Кочетков Е.С. 
  \newline
   \href{http://mat.univie.ac.at/~neretin/misc/reform/Kochetkov.djvu}
 {\it Алгебра и элементарные функции. 10 класс}, Просвещение, 1968. 
  \newline
    \href{http://mat.univie.ac.at/~neretin/misc/reform/Kochetkov-9.djvu}
 {\it Алгебра и элементарные функции. 9 класс}, Просвещение, 1969.
   
      \hangindent=0.3cm \noindent
  [Кур2007] 
  Курдюмова Н.А.
  \newline 
\href{https://refdb.ru/look/2089731.html}{\it  Былое: Воспоминания учительницы о Колмогоровской реформе.}
Архимед: Научно-методический сборник. Вып. 3. С. 20-44.

\hangindent=0.3cm
 {\footnotesize[Интересные, хотя и несколько эмоциональные картинки.
Цитаты: {\it ...в начале 60-х гг. школьное содержание математического образования 
не критиковали разве что первоклашки...}, 
\newline 
 {\it Всем хотелось революции, хотя бы в области математического образования}.
 \newline 
  Понятно, что 
в провале реформы виноваты все: учителя,
ЦК КПСС, общество, последовавшая через много лет после реформы  компьютеризация и даже <<особисты>>
...Еще цитата:
\newline  
{\it В середине 60-х г. никто и помыслить не мог, что с началом Реформы начнется холодная гражданская война, в которую будут вовлечены
и учителя математики, и ученые, и дети, и их родители. Война будет проходить под флагом усовершенствования школьной математики,
но сам курс математики станет только полем битвы, в которой чувства детей по большому счету никого не интересовать не будут.
И эта война истребит лучшие силы педагогов и математиков.}
\newline 
В отношении математиков это явное преувеличение (хотя страсти были изрядные),
а  педагоги в массовом предъявлении кузькиной матери друг другу преуспели
(но все же приводимые в статье факты о преждевременной смерти авторов учебников верны лишь частично).
Думаю (опять-таки, насмотревшись на образовательно-научный мир), 
что у части революционной тусовки желание показать оную мать было одним из воодушевляющих мотивов.]

}

     \hangindent=0.3cm \noindent 
     [Лар1958]
  Ларичев П.А. \href{http://mat.univie.ac.at/~neretin/misc/reform/Larichev6-7-1958.djvu}{\it  Сборник задач по алгебре для 6-7 классов.} 
  издание 10е,
  Издательство Министерства просвещения РСФСР, 1956
   
      \hangindent=0.3cm \noindent
   [Лев1959]
В. И. Левин, 
\newline 
\href{http://www.mathnet.ru/php/archive.phtml?wshow=paper&jrnid=mp&paperid=530&option_lang=rus}
{\it Некоторые вопросы преподавания математики в средней школе}, Математика, ее преподавание, приложения и история,
Матем. просв., сер. 2, 4, 1959, 145–150  

   \hangindent=0.3cm \noindent
[Луз1920] Н.Н.Лузин 
\href{https://www.mat.univie.ac.at/~neretin/misc/luzin/luzin-preface-graneville.pdf}
{\it Предисловие к русскому переводу книги Грэнвилль} {\it  <<Курс дифференциального и интегрального исчисления>>} М. 1926
 
    \hangindent=0.3cm \noindent
   [Ляп1959]
   А. А. Ляпунов,
   \newline  \href{http://www.mathnet.ru/php/archive.phtml?wshow=paper&jrnid=mp&paperid=532&option_lang=rus}
  {\it Реплики: О роли математики в среднем образовании}, Математика, ее преподавание, приложения и история,
   Матем. просв., сер. 2, 4, 1959, 152–154 
   
   \hangindent=0.3cm
{\footnotesize[Из статьи: {\it Если у кого-нибудь возникнет вопрос о загруженности школьников
и о недостатке учебного времени, то я предложу ему положить стоп­кой на стол школьные учебники и обязательные учебные пособия по
истории и литературе. Думаю, что разумное сжатие этой стопки поз­волит найти время для математики.}]

}

      \hangindent=0.3cm \noindent
   [МММ]
 Макарычев Ю.Н., Миндюк Н.Г., Муравин К.С. {\it Алгебра, 6-8 классы.} Под редакцией Маркушевича.
 Было много версий этой книги с разным составом авторов, начиная с 1970 года,
 у меня есть  \href{http://mat.univie.ac.at/~neretin/misc/reform/Algebra-6.djvu}{\it 6 класс за 1974 год}
 и
 \href{http://mat.univie.ac.at/~neretin/misc/reform/Makarychev-7.djvu}{\it  7 класс за 1976 год}.
 
 \hangindent=0.3cm 
 {\footnotesize[В дальнейшем учебник под тем же <<брендом>> с меняющимся составом авторов   продолжал существовать.
 После 1985 года слова <<под редакцией Маркушевича>>
 исчезли.]

}
  
     \hangindent=0.3cm \noindent
  [Мар1957]
 А. И. Маркушевич,
 \newline \href{http://www.mathnet.ru/php/archive.phtml?wshow=paper&jrnid=mp&paperid=375&option_lang=rus}
 {\it На XIX международной конференции по народному просвещению}, Математика, ее преподавание, приложения и история, Матем. просв., сер. 2, 1, 1957, 9–15 

      \hangindent=0.3cm \noindent
  [Мар1964]
 Маркушевич А.И.
 \newline
  \href{http://mat.univie.ac.at/~neretin/misc/reform/Markushevich1964.pdf}
{\it К вопросу о реформе школьного курса математики}.  Математика в школе, 1964, 6, c.4-8.
 
    \hangindent=0.3cm \noindent
 [Мар1979]
 Маркушевич А.И.
 \href{http://mat.univie.ac.at/~neretin/misc/reform/Mar1979.pdf}{\it О школьной математике.}
 Математика в школе, 1979, 4, 11-16.

   \hangindent=0.3cm \noindent
 [МатШ]  \href{http://mat.univie.ac.at/~neretin/misc/reform/V-min.ps}{В Министерстве просвещения РСФСР}

    \hangindent=0.3cm \noindent
 [ММЧ1978]
  А. И. Маркушевич, Г. Г. Маслова, Р. С. Черкасов (составители)
 \href{http://mat.univie.ac.at/~neretin/misc/reform/Puti.djvu}{\it На путях обновления школьного курса математики.}
 Сборник статей и материалов.
 Пособие для учителей, М., «Просвещение», 1978, 
  
     \hangindent=0.3cm \noindent
  [МСЧ1968]
 Маркушевич А. И., Сикорский К.П., Черкасов  Р. С.  {\it Алгебра и элементарные функции}, Просвещение, 1968 г.
 
   \hangindent=0.3cm \noindent  
 [Маш2006]
 Mashaal, M.
Bourbaki.
{\it A secret society of mathematicians.} Tran\-sla\-ted from the 2002 French original. American Mathematical Society, Providence, RI, 2006.
 \href{http://mat.univie.ac.at/~neretin/misc/reform/newmath.pdf}{\it  New Math in classroom, 134-145}.
 
    \hangindent=0.3cm \noindent
 [Нер2016]
Неретин Ю.А.
 \href{http://www.mat.univie.ac.at/~neretin/misc/Kiselev.html}
{\it Второе пришествие Киселёва (сага в резолюциях)}.
 
    \hangindent=0.3cm \noindent
 [Ник1971]
 Никитин Н.Н. \href{http://mat.univie.ac.at/~neretin/misc/reform/Nikitin1971.djvu}{\it  Геометрия, 6-8 класс.} Просвещение, 1971,
 
 \hangindent=0.3cm
{\footnotesize [В первой версии это был учебник Никитина-Фетисова  [НФ1956],
он подвергся жесткой критике (см. прицеп к ссылке [Обс1957]), и дальше это стал учебник планиметрии Никитина. 
Фетисов опубликовал   <<Стереометрию>>	[Фет1963]].

}

   \hangindent=0.3cm \noindent
 [НМ1971]
 Никитин Н.Н., Маслова Г.Г.
 \newline 
  \href{http://mat.univie.ac.at/~neretin/misc/reform/Nikitin-Maslova-1971.djvu} 
 {\it Сборник задач по геометрии для 6-8 классов,}  15 издание, Просвещение,  1971,
  
     \hangindent=0.3cm \noindent
 [НФ1956]  Никитин Н.Н., Фетисов А.И. Геометрия. {\it Учебник для семилетней и средней школы}, Учпедгиз, 1956 
      \hangindent=0.3cm \noindent
 
  \hangindent=0.3cm \noindent
 [Нов1950]
 Новосёлов С.И.
 \newline 
  \href{http://mat.univie.ac.at/~neretin/misc/reform/Novoselov1950.pdf}{\it К вопросу о введении элементов дифференциального и интегрального} {\it  исчислений в курс средней школы.}
 Математика в школе, 1950, 2, 35--39.
 
 \hangindent=0.3cm
 {\footnotesize Аргументированные возражения против введения элементов мат.анализа в школе.
 Многие аргументы устарели, как-никак прошло 70 лет, но текст этот интересен.
 
}
  
     \hangindent=0.3cm \noindent
  [Нов1962]
  Новосёлов С.И. \href{http://mat.univie.ac.at/~neretin/misc/reform/Novoselov.rar} {\it  Тригонометрия}, Учпедгиз, 1962
  
     \hangindent=0.3cm \noindent
  [Обс1957]
  \href{http://www.mathnet.ru/links/ed01564e66d102f29e54d96a70e4c7f1/mp440.pdf}{\it Обсуждение новых стабильных учебников по математике},
  Математика, ее преподавание, приложения и история, Матем. просв., сер. 2, 1, 1957, 195–209
  
  \hangindent=0.3cm  
 {\footnotesize [Заседание Московского математического общества, Секция средней
 	школы, обсуждение учебника Никитина и Фетисова, 13.11.56:
 	\newline
 	{\it 
 	А. И. Ф е т и с о в излагает соображения, положенные авторами в основу
 	при составлении учебника: близость школьного курса к жизни и приложениям; знакомство с чертежными и измерительными инструментами (политехнизация); идея геометрического преобразования, в частности гомотетия как
 	исходный пункт при изложении подобия фигур. До и после напечатания учебник обсуждался в ряде педагогических коллективов, на заседаниях кафедр педвузов (некоторые названы) и наряду с критикой встречал положительные
 	оценки.
 	\newline
 	Я. С. Д у б н о в считает, что авторами руководили некоторые добрые намерения: дальнейший отход от евклидовских традиций в описательной части... Однако свежие педагогические идеи
 	ослабляются, хуже того — компрометируются крайне несовершенным их
 	воплощением...}[Дальше много примеров] {\it Предложение: на
 	два-три года примириться с преподаванием по любым учебникам, пока не будет создана полноценная книга.
 		\newline
 	И. М. Я г л о м присоединяется к Я. С. Дубнову в положительной оценке
 	общих идей новой книги; поэтому он считает, что учебник нельзя ставить на
 	одну доску с отвергнутой в свое время книгой Гангнуса и Гурвица. Однако
 	правильные идеи не реализованы как следует из-за спешки в работе авторов;
 	поэтому в книге оказалось большое число серьезных дефектов, включая прямые ошибки (приводит примеры)... 
 	\newline
 	Однако авторам учебника следует предоставить возможность работать
 	над его улучшением.
 	\newline
 	В. А. У с п е н с к и й, возражая И. Я. Танатару, заявляет, что учебник
 	ему понравился тем, что содержит явно выраженные идеи и в этом выгодно
 	отличается от книги Киселева. Соглашается с тем, что ошибок много; в частности, изложение вопроса об измерении совершенно неудовлетворительно.
 	Считает, что в качестве стабильного этот учебник не годится, но он может
 	быть временно допущен наряду с учебником Киселева. Общество должно
 	добиваться того, чтобы были наказаны конкретные виновники выпуска миллионным тиражом учебника с такими серьезными недостатками.
 	\newline
 	Н. М. Б е с к и н, не касаясь вопроса о качестве учебника, заявил, что
 	Министерство просвещения и Учпедгиз неправильно интерпретируют постановления ЦК КПСС и Правительства о стабильных учебниках, рассматривая
 	их как директиву к изданию единственных учебников по каждому предмету. 
 	\newline
 	В. А. Е ф р е м о в и ч резко выступает против учебника и тех, кто пытается взять его хотя бы частично под свою защиту. Выпуск Учпедгизом этой
 	кииги считает тяжким нарушением государственных интересов. Необходимо
 	издание курса геометрии специально для учителей средней школы; вместе
 	с тем нужен новый учебник для учащихся. (Резкий тон выступления
 	В. А. Ефремовича вызвал сдерживающее замечание председателя.)
 \newline
 	И. С. Г р а д ш т е й н приводит пример (стр. 193) громоздкости в формулировках теорем. Нужно издать несколько учебников. 
 	\newline
 	Е. Н. Е в з е р и х и н а (учительница школы № 58) считает, что с методической стороны учебник совершенно неудовлетворителен. Доказательства ряда
 	теорем не обладают достаточной полнотой и вызывают затруднения не только
 	у учащихся, но и у преподавателей. Как курьез, можно отметить, что в
 	Институте усовершенствования учителей существует специальный семинар,
 	имеющий целью <<расшифровку>> этого учебника. Не ожидая санкции свыше,
 	многие учителя вернулись к преподаванию по учебнику Киселева. Под аплодисменты собравшихся оратор шутливо предлагает наказать лиц, ответственных за появление учебника, потребовав от них в течение одного года преподавать по этой книге.
 	\newline
 	Н. Я. В и л е н к и н приходит к выводу, что дело подготовки учебников
 	нельзя доверять Министерству Просвещения и Учпедгизу, допустившим непростительную халатность. Перед Министерством следует поставить вопрос
 	о персональной ответственности за выпуск этого учебника и об усилении руководства редакцией математики Учпедгиза. При сложившемся положении должно
 	быть позволено пользоваться старыми учебниками.} 
 \newline
 Я привел эту подборку по разным причинам. Во-первых, меня забавляет оценка идеи низвергнуть
 Киселева как доброго намерения. Безусловно добрым намерением было бы написать
 учебник лучше киселёвсого, но ведь всем было очевидно, что это не  достигнуто. 
 С другой стороны забавляет степень живости, с которой в 1955-56гг. происходило 
 обсуждения новых учебников (упомянутый мельком учебник Гангнуса и Гурвица
 после нескольких лет писания писем в инстанции стал предметом еще более живого обсуждения
 в 1936-37гг., см. [Нер2016]). После обсуждения учебник Никитина-Фетисова был существенно изменен,
 а Фетисов исчез из числа авторов.
 Интересно, а как на самом деле происходило обсуждение в 1966-68гг.??  Понятно, что молчали не все
 (отзвуки  есть в [О-про1967-4]).
 \newline
 \qquad
 На другую тему.
 	Из резолюции объединенного семинара  кафедр высшей математики
московских втузов, 14.XI.1956:
\newline
{\it 
  Особое место в ряду указанных учреждений принадлежит Академии
педагогических наук РСФСР (а также специальным институтам методов
«обучения в союзных республиках), смысл существования которых только
и состоит в своевременной разработке разумных методов преподавания
и в создании серьезной учебной литературы. За много лет функционирования, в «состоянии постоянного возбужденного бездействия», эти
научно-педагогические институты, как оказалось на поверку, ничего не
смогли дать ни для математической педагогики, ни для практики преподавания, за исключением разве лишь нескольких 
прожектерских планов коренной ломки всей схемы обучения, вызвавших, естественно, резкое
противодействие со стороны учительства и породивших недоверие учительства к участию научных работников в делах средней школы.}]

} 
  
   \hangindent=0.3cm \noindent  
 [Объ1965]
 \href{http://mat.univie.ac.at/~neretin/misc/reform/volume1965.pdf}{\it Объем знаний по математике для восьмилетней
школы,} «Математика в школе», 1965, N2, c.21-24
 
   \hangindent=0.3cm \noindent 
 [О-про1967-3] \href{http://mat.univie.ac.at/~neretin/misc/reform/OP1967-03.pdf}{\it О проекте программы средней школы}.
  Математика
 в школе, 1967, 3,  28-38.

\hangindent=0.3cm
 {\footnotesize Несколько статей. Одобряем проект программы!
\newline  
\it
Представленные комиссией проекты программ в основном приемлимы,
так как они полнее отвечают требования перестройки школы,
установленным решениям ЦК КПСС и Совета Министров СССР о  средней школе, чем нынедействующие программы.

}
  
  \hangindent=0.3cm \noindent 
[О-про1967-4]  \href{http://mat.univie.ac.at/~neretin/misc/reform/OP1967-04.pdf}{\it О проекте программы средней школы}. Математика
в школе, 1967, 4, 25-36.

\hangindent=0.3cm
{\footnotesize Это набор статей, обсуждающих проект программы,
преимущественно положительных или предлагающих усовершенствования.
Привожу некоторые скептические  комментарии. Надо сказать, что
и положительные отзывы содержат разные точные и неприятные замечания, которые должны были бы заставить  задуматься о реалистичности программы вообще.
\newline
Непомнящий П.Е. (Ленинград) {\it  По силам ли?},
\newline
{\it Как такое можно делать? Ссылки на какой-то опыт несолидны.
	Чудес ведь не бывает\dots 
\newline	 
Нельзя опрометчиво вводить эти программы, надо поставить обязательно опыты в обычных школах и посмотреть, что получится после X класса, если учащиеся будут обучаться по таким программам.
\newline
Естественно возникнут следующие трудности:
\newline
1)
Осилят ли учащиеся школы значительно более сложный материал, если
сравнительно более простой и легкий материал теперь в школе слабо усваивается?
\newline
2) Кто будет учить по новым программам математики учащихся
IX-X классов? Некоторые учителя средней школы не знают этого материала, а для изучения его нужно значительное время. Планируемые в крупных городах курсы для учителей едва ли справятся с этим. Что же говорить о школах небольших городов и сельской местности?
\newline
3) Зачем изучать в школе такой материал, который в вузах будет изучаться заново, а на производстве не нужен.\dots}
\newline
Г.Н.Скобелев {\it Пожелания к проекту программы.} (г.Херсон)
\newline
{\it Авторы программы не выдержали принципиального подхода к фкультативным курсам. Судя по объяснительной записке, факультатив должен помочь учащимся, желающим углубить свои математические знания. Однако вследствие нехватки времени, авторы программы порой предлагают использовать
время, отводимое на факультативные курсы, как отдушину для занятий вопросами, без которых усвоение учебной программы невозможно.}
\newline
Денисова Т.Н. {\it Замечания по проекту программы.}
 (Москва)
 \newline
 {\it Главное же пожелание Министерсту просвещения СССР -- не торопиться с введением
 новой программы. Сначала -- солидная экспериментальныя проверка новой програмы на базе новых учебников; затем -- внесение в программу результатов обсуждения и проверки; подготовка и перепедготовка учителей средней школы и только после выполнения этих условий -- введение новой программы.}

}

  \hangindent=0.3cm \noindent 
[О-про1968-1] 
\href{http://mat.univie.ac.at/~neretin/misc/reform/OP1968-01.pdf}{\it О проекте программы средней школы}. Математика
в школе, 1968,1, 16-24.

\hangindent=0.3cm
{\footnotesize
Рефераты статей, присланных в журнал и там не опубликованных.}
  
     \hangindent=0.3cm \noindent
[ПБДЛШГ1960] Пиаже Ж.,  Бет Э., Дьедонне Ж.,  Лихнерович~А., Шоке~Г., Гаттеньо~К.,
 \href{http://mat.univie.ac.at/~neretin/misc/reform/Piazhe_Prepodavanie.djvu}{\it О преподавании математики},
 Учпедгиз, 1960. Translated from 
   J.~Piaget, E.~W.~Beth, J.~Dieudonn\'e, A.~Lichnerowicz, G.~Choquet, C.~Gattegno.
 {\it L’Enseignement des Math\'ematiques.}
  Delachaux-Niestl\'e, Paris--Neuchatel,  1955.
 
    \hangindent=0.3cm \noindent
  [Пог1969] Погорелов А.В. Элементарная геометрия. 
  \href{http://mat.univie.ac.at/~neretin/misc/reform/geometr-pogorelov1-1969.djvu}{\it Планиметрия}, Наука, 1969, 
  \href{http://mat.univie.ac.at/~neretin/misc/reform/geometr-pogorelov2-1970.djvu}{\it Стереометрия}, Наука 1970,
   \href{http://mat.univie.ac.at/~neretin/misc/reform/geometr-pogorelov3-1972.djvu}{\it Элементарная геометрия}, Наука, 1972,
   Было еще  расширенное издание, Наука, 1977. 
  
      \hangindent=0.3cm \noindent
  [Пог1982] Погорелов А.В. 	  \href{http://mat.univie.ac.at/~neretin/misc/reform/geometr-06-10-1982.djvu}
 {\it Геометрия. Учебное пособие для 6-10 классов}, Просвещение, 1982 (тир 3 302 000). 
 
 \hangindent=0.3cm
{\footnotesize  В 1981 году было выпущено
  как пробный учебник (тир. 263 000)
  \newline 
 [А.В.Шевкин,
  \href{http://shevkin.ru/?action=Page&ID=134}{\it «Школьное обозрение», 2002, № 5}:
 {\it Привлекательность учебника} [планиметрии]{\it связана, видимо,
 с тем, что он является развитием хорошо продуманных учебников и задачников прошлых лет. Но самое трудное для учащихся и учителя при работе
 по этому учебнику — это отслеживание порядка вершин треугольников при обсуждении их равенства и подобия, довольно сложные для учащихся доказательства первых теорем (например, признаков равенства треугольников).
Эти трудности как раз и произрастают из желания автора все вывести из аксиом и не пользоваться, например, наложением треугольников при доказательстве признаков равенства. Обучающий и воспитательный эффект от такого способа обучения не сопоставим с теми трудностями, которые испытывают учащиеся и учителя.
\newline 
Усвоение первых тем по  учебнику} [стереометрии] {\it затрудняется тем,
что основные изучаемые геометрические объекты — точки, прямые и плоскости — «висят» в пространстве,
не имея опоры в виде знакомых с детства геометрических тел. Но опытные учителя умеют компенсировать этот недостаток, 
иллюстрируя изучение теории с помощью геометрических тел и решая с опережением на год простейшие задачи на построение сечений.}
См. также критику в  [Глад2009] и  [Вин2015]]

}

   \hangindent=0.3cm \noindent  
 [Пон1979]
Понтрягин  Л. С., \href{http://www.mathnet.ru/php/archive.phtml?wshow=paper&jrnid=mo&paperid=284&option_lang=rus}
{\it Этика и арифметика: Человек, труд, мораль.}
 Соц. индустрия. — 1979, 21 марта.
 
 \hangindent=0.3cm
{\footnotesize [эта статья, часто упоминается в связи с обсуждаемым предметом, Понтрягин там морализировал и  на кого-то
и на что-то ругался, адресаты, скорее всего
обиделись. Но о школе сказано только следующее: {\it С другой стороны, серьезные недостатки в
программе средней школы поставили определенный заслон на пути к высшему математическому образованию, создали ажиотаж вокруг репетиторов и опять-таки
искусственное выделение групп молодежи, идущей в математику. Этого быть не
должно! С этим необходимо бороться!}]

}

   \hangindent=0.3cm \noindent  
[Пон1980]
Понтрягин Л. С. \href{http://mat.univie.ac.at/~neretin/misc/reform/pontryagin-communist.html}{\it О математике и качестве её преподавания},
Коммунист, 1980. — № 14. — С. 99—112. 
На эту статью печатались отклики в  \href{http://mat.univie.ac.at/~neretin/misc/reform/Kommunist1980-18}{Коммунист. —1980. —  № 18}
и в \href{http://mat.univie.ac.at/~neretin/misc/reform/Kommunist1982-2}{Коммунист 1982. — № 2}.
  
     \hangindent=0.3cm \noindent
[Пон1998]
Понтрягин Л. С.
\newline
 \href{http://ega-math.narod.ru/LSP/book.htm}
{\it Жизнеописание
Льва Семёновича
Понтрягина,
математика,} {\it 
составленное им самим.
Рождения 1908 г.} Москва, 1998.

\hangindent=0.3cm
{\footnotesize[ Цитаты: {\it После того как в конце 1977 года до математиков, занимающихся наукой, 
наконец-то дошло, что в средней школе неблагополучно, десять академиков-математиков обратились с письмом в ЦК. 
В этом письме мы выражали тревогу по поводу происходящего в школе.
\newline 
После этого в 78-м году министр просвещения СССР М. А. Прокофьев обратился в Отделение математики АН СССР с просьбой заняться вопросами преподавания.
В результате состоялось сперва заседание Бюро Отделения математики, а затем Общее собрание Отделения математики, 
на котором присутствовали представители Министерств просвещения СССР и РСФСР. Был также и А. Н. Колмогоров.
Как на Бюро, так и на Общем собрании Отделения были решительно осуждены действующие учебники и учебные программы. 
Общее собрание Отделения продолжалось много часов и происходило в большом накале.
\newline 
Рассматривались конкретные дефекты учебников, и подавляющему большинству присутствующих было совершенно ясно, что так оставаться дальше не может.
Решительными противниками каких бы то ни было действий, направленных на исправление положения,
были академики С. Л. Соболев и Л. В. Канторович, которые говорили, что надо подождать.
Но, несмотря на их сопротивление, было принято решение, требующее вмешательства в вопросы преподавания в средней школе. 
В частности, было вынесено решение об организации комиссии по преподаванию при Отделении.
Выполнение этого решения было поручено Бюро Отделения. Следующее заседание Бюро Отделения занялось образованием комиссии по преподаванию.
И здесь возникли разногласия между математиками не по существу, а по тому, кто же будет возглавлять дело.
\newline 
Обнаружилось, что имеется два претендента — академики А. Н. Тихонов и И. М. Виноградов. И оба они были в какой-то степени поддержаны.
Поэтому было принято осложняющее всё дело решение образовать две комиссии. Одну под председательством Тихонова,
другую — под председательством Виноградова. Наличие двух комиссий указывало на раскол между математиками и затрудняло работу.
В результате длинных перипетий в Отделении, продолжавшихся около трёх лет, обе комиссии были ликвидированы и была образована одна новая комиссия,
которую возглавил Виноградов и которая называется комиссией по преподаванию математики в средней школе.
Я был единственным заместителем Виноградова.
\newline 
После смерти Виноградова председателем комиссии назначен я, а моим заместителем А. С. Мищенко — профессор мех-мата МГУ....
\newline 
Уже после того, как Отделение в 1978 году высказало своё чёткое мнение по вопросу о негодности действующих учебников и программ,
дело долго не двигалось с места....
\newline 
 С другой стороны, ещё раньше, чем Колмогоров приступил к своему изменению преподавания, учебник по геометрии начал писать хороший советский геометр,
 ныне академик, А. В. Погорелов. Погорелов рассказывал мне, что он предлагал Колмогорову использовать его учебник, но тот отказался, 
 так как учебник Погорелова не соответствовал идеям Колмогорова. Такой же отказ Погорелов получил и от Тихонова, 
 когда тот занялся подготовкой новых учебников. Здесь причина была, по-видимому, другая.
 Тихонов хотел держать всё в своих руках и не хотел вовлекать в дело столь авторитетного человека, как Погорелов, 
 поскольку Тихонов никак не мог считаться его руководителем. Комиссия Виноградова, ещё во времена существования двух комиссий, 
 рекомендовала учебник Погорелова как пригодный для школы. Это предложение было принято Прокофьевым, 
 но оно не устраивало и не устраивает Тихонова, который хочет протащить свои учебники.
 В настоящее время в семи областях Российской Федерации производится эксперимент по учебникам Тихонова, в то время 
как в остальных республиках и областях РСФСР введён официально учебник Погорелова 3. По-видимому,
Тихонов вместе с методистами Министерства просвещения РСФСР надеются продвинуть свой учебник на всю Российскую Федерацию.
\newline 
Так горько обстоит дело с геометрией. С алгеброй дело обстоит ещё хуже. Учебник, подготовленный группой Тихонова по алгебре, 
признан комиссией Виноградова не вполне удовлетворительным, а другого готового учебника у нас пока ещё нет, хотя и есть заготовки... 
Таким образом, дело, начатое ещё в 1978-м году, только к 82-му году начало немножко сдвигаться с места, однако ещё очень недостаточно.}]

}

 \hangindent=0.3cm \noindent
[Пон2008]
Понтрягина  А. И.,
\newline  \href{http://www.mathnet.ru/php/archive.phtml?wshow=paper&jrnid=mo&paperid=84&option_lang=rus}
{\it Из воспоминаний о Льве Семеновиче Понтрягине}, Матем. обр., 2008, № 3(47), 2–26 

\hangindent=0.3cm
{\footnotesize[Цитата: {\it Я не помню, как у нас в доме появился сотрудник журнала <<Коммунист>> — Леонид Вита­льевич Голованов.
Он отлично понимал значение проблемы и главный редактор этого журнала
Косолапов тоже. Они дали согласие изложить взгляды Понтрягина в своем журнале. Но отдать
дань похвалы партии и правительству первым условием было. Лев Семенович категорически
отказался от этого. Препирательство между ними шло недели две. В конце концов, Лев Семенович согласился,
но сам писать наотрез отказался. Эту часть статьи написал (как положено по
канону), Леонид Витальевич Голованов, 16 месяцев ждали появления этой статьи в печати! В
течение этого периода Леонид Витальевич временами звонил Льву Семеновичу, и в его словах
теплилась надежда, что статья все же выйдет в свет, и что он ходит на цыпочках и говорит
шепотом, как бы кого-нибудь не вспугнуть, не потревожить...}
\newline 
Ю.Н. Из публикации [Пон1980] очевидно, что эта статья многократно рецензировалась.]

}

   \hangindent=0.3cm \noindent
[Прог1959] {\it Проект программ по математике для средней школы}, «Математика  в школе», 1959,    4,      стр.       1—14

   \hangindent=0.3cm \noindent
[Прог1960]
\href{http://www.mathnet.ru/php/archive.phtml?wshow=paper&jrnid=mp&paperid=645&option_lang=rus}
{\it Проект программ по математике для средней школы}, Математика, ее преподавание, приложения и история,  5, 1960, 118–126
 
   \hangindent=0.3cm \noindent
[Прог1961]
\href{http://mat.univie.ac.at/~neretin/misc/reform/msh1961-1.pdf}
{\it Программа средней общеобразовательной политехнической школы} {\it с трудовым обучением}, Математика в школе, 1961, 1, 6-12

   \hangindent=0.3cm \noindent
[Прог1967]
\href{http://mat.univie.ac.at/~neretin/misc/reform/msh1967-1.pdf}
{\it Проект программы средней школы по математике}, Математика в школе, 1967, 1, 5-24.
 
    \hangindent=0.3cm \noindent
[Прог1967-м]
\href{http://mat.univie.ac.at/~neretin/misc/reform/msh1967-3-m.pdf}
{\it О проекте программы средней школы по математике,} Математика в школе, 1967, 3, 28
 
[Прог1967-э]
\href{http://mat.univie.ac.at/~neretin/misc/reform/msh1967-3-e.pdf}
{\it К программе курса математики}, Математика в школе, 1967, 3, 29

   \hangindent=0.3cm \noindent
[Прог1968]
\href{http://mat.univie.ac.at/~neretin/misc/reform/msh1968-2.pdf}
{\it Программа по математике для средней школы},  Математика в школе, 1968, 2, 2-20
   
      \hangindent=0.3cm \noindent
   [РезСО1980]
   \newline 
  \href{http://www.math.nsc.ru/LBRT/g2/english/ssk/resolution_council.pdf}
 {\it Резолюция Сибирского отделения АН СССР от 25 декабря 1980г},
 
 \hangindent=0.3cm
 	{\footnotesize[Она  не была опубликована, но, по-видимому, была известна.
 		Про Соболева, который, скорее всего, был главным двигателем, см. прицеп к ссылке [Соб1980].]}
 
  	        \hangindent=0.3cm \noindent 
[Рыб1931]
  Рыбкин Н. 
  \newline 
  \href{http://mat.univie.ac.at/~neretin/misc/reform/Rybkin-trig-1931.djvu}
 {\it Учебник прямолинейной тригонометрии собрание задач},
  Государственное учебно-педагогическое издательство,
  1931
 
    \hangindent=0.3cm \noindent
[Рыб1961]
  Рыбкин Н.А. {\it Сборник задач по геометрии}, \href{http://mat.univie.ac.at/~neretin/misc/reform/Rybkin-zad1.djvu}{ Часть 1},
  \href{http://mat.univie.ac.at/~neretin/misc/reform/Rybkin-zad2.djvu} {Часть 2},   Государственное учебно-педагогическое издательство, 1961
  
  \hangindent=0.3cm
  {\footnotesize[Знаменитая книга. В 1973 году было 39ое издание.]
  
}
     
        \hangindent=0.3cm \noindent
[Рыб1951]
  Рыбкин Н. \href{http://mat.univie.ac.at/~neretin/misc/reform/Rybkin-trig-1951.djvu}
  {\it  Прямолинейная тригонометрия,}
  Учпедгиз, 1951
 
     \hangindent=0.3cm \noindent
  [Соб1980]
 	С. Л. Соболев, \href{http://www.mathnet.ru/php/archive.phtml?wshow=paper&jrnid=smj&paperid=1894&option_lang=rus}
 	{\it В редакцию журнала “Коммунист}, Сиб. матем. журн., 49:5 (2008), 970–974
 	
 	\hangindent=0.3cm
 	{\footnotesize[Письмо было послано в 1980 году и опубликовано не было. Оно может звучать убедительно
 	для современного читателя. Но, цитирую,
 	 {\it Школьные программы по математике до этого времени не содержали
даже элементов математического анализа} (неправда, см. выше Раздел 2).
{\it В старших классах такими явно устаревшими оказались изучение
сложных приемов приведения тригонометрических формул к виду, удобному
для логарифмирования, изучение специальных приемов решения косоугольных
треугольников, сложные проценты и другое} (я ничего такого в предреформенных
учебниках не нашел). {\it В старых программах не было ничего сказано об отображениях}
(не  так,  осевая и центральная симметрии, преобразование подобия (?), а также проектирование были
введены в программу Глаголевым в 1938г.).
{\it В старых программах неравенства вводились в самом конце курса математики — в X классе}
(конечно, не так). Кое-что из сказанного было правдой
или отчасти правдой.  {\it В прежних программах по математике в советской средней школе не
существовало даже понятия о векторе.} Векторы в небольших дозах были
в учебнике Новоселова [Нов1962] по тригонометрии, работавшем в 1957-1966гг.;
учебник был отменен в связи с Реформой-1959, векторы должны были перейти
в курс геометрии, но учебник Болтянского и Яглома не прошел.  
 {\it В младших классах это были кустарные приемы решения арифметических 
(«текстовых») задач} (были, хотя, по-видимому, сильное смягчение было уже до реформы).]

}

       \hangindent=0.3cm \noindent
[Сов1964]
\newline
 \href{http://mat.univie.ac.at/~neretin/misc/reform/soveshchanie1964}
{\it Совещание  по  проблемам  математического 
образования в  средней  школе} {\it (Москва)}. Математика в школе,  1964, No  6, стр.  90—91. 

\hangindent=0.3cm
{\footnotesize[Цитаты: {\it А.Н.Колмогоров считает, что {\bf курс геометрии должен быть сложным,
но не следует исключать из него элементарные понятия}, 
<<Преподавание в старших классах, по мнению академика Колмогорова, должно вестись
на более строгой теоретической основе>>, <<серьезное внимание должно быть уделено элементам математического анализа>>}.
\newline 
 Были и такие слова, которые позже, кажется, уже не повторялись:
 \newline  {\it Вместе с тем докладчик отметил, что эти изменения нельзя производить
наспех на основе умозрительных заключений. Необходима постановка серьезного опережающего эксперимента.}
\newline 
Впрочем, смысл этой фразы уточнялся следующим 
\newline 
{\it Вместе с тем следует
проявить чувство меры и на ближайшие годы не вводить таких вопросов, как элементы теории вероятностей, 
математической статистики и ряду других} [напомню, что Программу-1968 они не вошли].
{\it Выступавшие единодушно поддержали идеи, высказанные А.Н.Колмогоровым и} [министра просвещения РСФСР]
 {\it Е.И.Афанасенко.}] 

}
 
    \hangindent=0.3cm \noindent
[Стен1978]
{\it Стенограмма Общего собрания отделения математики, посвященного обсуждению школьных программ и учебников по математике. 3 декабря 1978}.
В книге   \href{http://mat.univie.ac.at/~neretin/misc/bunt.djvu}{[КС2012]}.

   \hangindent=0.3cm \noindent
[Тих2009]
Тихомиров В.М. {\it Педагогические замыслы А.Н. Колмогорова о курсе геометрии в школе}.
\href{http://yspu.org/kolmogor/doc/kolmogor_08.pdf}{Труды шестых Колмогоровских чтений}. 
 Изд-во ЯГПУ, 2008,
10-17 

\hangindent=0.3cm
{\footnotesize[Цитаты:{\it
Я вижу две основные причины того, что усилия Андрея Николаевича
в итоге не привели к благоприятному исходу. Первая из них заключена в идеологии того общества, 
в котором тогда приходилось всем нам
жить и трудиться. Первоосновой всего, высшей целью жизни и деятельности каждого человека объявлялось 
тогда укрепление и развитие государства. Не личности, а именно государства. И руководящая структура
государства – партия – определяла, в частности, цели и смысл образования. Андрей Николаевич вынужден
был подчиняться этому порядку
вещей, но на самом деле он и сам искренне считал, что прогресс развития нашей страны невозможен
без широко образованной творческой
интеллигенции. И это было стимулом для его трудов и усилий.
В числе непременных аксиом того времени было требование единого, образования: каждый
должен был получить в точности то же образование, что и все остальные. Постановление партии и правительства
о реформе школьного образования было принято в 1966 году, и именно
тогда Андрею Николаевичу было поручено осуществить ту часть реформы,
которая относилась к математике. Естественно, должен был встать
вопрос: как и чему учить детей в нашей бескрайней, многонациональной, 
столь разнородной и неблагоустроенной стране с огромным числом
неблагополучных семей, детей с неполноценным умственным развитием
– читатель легко продолжит список всех наших трудностей тех и нынешних лет.
И при этом учить всех, и одинаково! Андрей Николаевич
взялся за осуществление этой реформы, не имея в виду существенно
менять исходные позиции упомянутого постановления, т. е. взялся за
неосуществимое предприятие.
\newline 
Но надо сказать и о второй причине, в силу которой планы 
по реформе образования оказались не до конца реализованными. Андрей Николаевич пребывал в мечтах о светлом будущем.
В этих мечтах он был
идеалистом...
\newline 
По мнению Андрея Николаевича (впрочем, это
можно было истолковать, как предписание в постановлении партии и
правительства от 1966 года), курс школьной математики должен быть
научным, строгим и современным. А эту цель в современном обществе
(а возможно, и в обществе будущего) осуществить невозможно.
\newline 
На двадцать второй странице пособия по геометрии для шестого
класса (изд. 1974 г.) написано: <<Систематический курс геометрии имеет
такое логическое строение:
\newline 
1. Перечисляются основные, принимаемые без определений, понятия.
\newline 
2. При их помощи даются определения других геометрических понятий.
\newline 
3. Формулируются аксиомы.
\newline 
4. На основе аксиом и определений доказываются теоремы.>>
\newline 
Уже в этом маленьком фрагменте отражается взгляд автора на своего читателя, не соответствующий возможностям последнего. Для понимания этого текста нужна высокая интеллектуальная культура, которая
формируется долго.}
]

}

   \hangindent=0.3cm \noindent 
[Тих2009]
Тихомиров В. М.  \href{http://mat.univie.ac.at/~neretin/misc/reform/Tihomirov.pdf} {\it А. И. Маркушевич (1908-1979)}.
 Историко-ма\-те\-ма\-ти\-че\-ские исследования. — М.: «Янус-К», 2009. — Т. 13 (48). — С. 128—137.

   \hangindent=0.3cm \noindent
[Фей1965] Feynmann R.  \href{http://calteches.library.caltech.edu/2362/1/feynman.pdf}
{\it New texbooks for "new" mathematics}, Engineering and science 
March 1965, 28,  6, 9-15.

    \hangindent=0.3cm \noindent
 [Фет1963]
Фетисов А.И. \href{http://mat.univie.ac.at/~neretin/misc/reform/Fetisov1963.djvu} {\it  Геометрия для 8-10 классов.}
Экспериментальный учебник.  Издательство АПН РСФСР, 1963 

\hangindent=0.3cm
{[\footnotesize Из аннотации:
{\it Все доказательства теорем и решения задач основаны на методе геометрических преобразований:
 симметpии, пepeнoca, вращения и подобия, что значительно 
упрощает в сравнении с учебником А. П. Киселева изложение и усвоение учебного материала.}
\newline 
Я этот учебник в старших классах брал в библиотеке, он мне показался очень тяжелым.]

}

   \hangindent=0.3cm \noindent
 [Фил2014]
Phillips, Christopher J. {\it The New Math: A Political History}, Chicago University Press,   2014
(можно найти в интернете)

   \hangindent=0.3cm \noindent
 [Хин1941]
 Хинчин А.Я.  \href{http://mat.univie.ac.at/~neretin/misc/reform/Hinchin1941.pdf}
{\it О понятии отношения двух чисел}, Математика в школе, 1941, 2, с.13-25 

\hangindent=0.3cm
{\footnotesize [Статья написана в связи с правкой Хинчиным учебника арифметики Киселева.
Цитата: {\it Большинство товарищей, возражающих против отождествления понятий отношения и частного,
 соглашается вместе с тем, что всякое отношение есть частное, но указывает вместе с тем, 
 что не всякое частное есть отношение...} Если кто не понял сказанного, можно обратиться к оригиналу.]

}
  
     \hangindent=0.3cm \noindent
  [Хин1961]
 	А. Я. Хинчин, 
 	\newline 
 	\href{http://www.mathnet.ru/php/archive.phtml?wshow=paper&jrnid=mp&paperid=677&option_lang=rus}
 	{\it О так называемых “задачах на соображение” в курсе арифметики}, 
 	Математика, ее преподавание, приложения и история, Матем. просв., сер. 2, 6, 1961, 29–36
 
    \hangindent=0.3cm \noindent
  [Хин1963]  Хинчин А.Я. \href{http://mat.univie.ac.at/~neretin/misc/reform/Hinchin.djvu} {\it Педагогические
  статьи.}, М., АПН СССР, 1963   

    \hangindent=0.3cm \noindent
[ЦК1977]
 ЦК КПСС и  Совет Министров СССР
 \newline 
 \href{http://mat.univie.ac.at/~neretin/misc/reform/TSK.pdf}
{<<\it О дальнейшем совершенствовании обучения, воспитания учащихся}
		 {\it общеобразовательных школ и подготовки их к труду>>}, постановление от 22 декабря 1977 г.
  
   \hangindent=0.3cm \noindent
   [Чер1988]
  Черкасов Р.С.
  \newline 
   {\it О методическом наследии А.И. Маркушевича (к 80-летию со дня рождения)}
   Математика в школе. 1988. № 2.
   
      \hangindent=0.3cm \noindent
  [Чер1993]
Черкасов Р.С.
\href{http://mat.univie.ac.at/~neretin/misc/reform/kolmogorov1.doc}
{\it Андрей Николаевич Колмогоров и школьное математическое образование.}
Колмогоров в воспоминаниях / Под ред. А.Н. Ширяева. – М.: Физматлит, 1993. С. 583–604.

\hangindent=0.3cm \noindent
[Чет1950] Четверухин Н.Ф. 
\newline
\href{http://mat.univie.ac.at/~neretin/misc/reform/Chet1950.pdf} {\it О научных принципах преподавания геометрии в советской школе 1950}, Математика в школе, 1, 5-13.

\hangindent=0.3cm \noindent
[Шаб1968] Шабат Б. В.
\newline
  \href{http://mat.univie.ac.at/~neretin/misc/reform/Shabat1968.pdf}{\it Алексей Иванович Маркушевич (к шестидесятилетию со дня рождения)}.
Математика в школе, 1968, 2, 93-94.

   \hangindent=0.3cm \noindent
[Шар2001]   Шарыгин И.Ф.
\newline
\href{http://www.mccme.ru/edu/index.php?ikey=shar_mathedu}
{\it Математическое образование: вчера, сегодня, завтра...}, 2001
 
    \hangindent=0.3cm \noindent
[Шар2002]   Шарыгин И.Ф.
\newline
	\href{http://www.strana-oz.ru/2002/2/ot-kakogo-konya-primet-smert-rossiyskaya-matematika}
	{\it От какого 'коня' примет смерть российская математика,} Отечественные записки, 2003, вып.2  	

   \hangindent=0.3cm \noindent
[Шаф1989]	
 Шафаревич И. Р.,
 \newline 
  \href{http://ega-math.narod.ru/LSP/ch8.htm#a}
 {\it Так сделайте невозможное (К 80-летию Л. С. Понтрягина)}, Советская Россия, 16.04.1989.
 \newline 
 {\footnotesize[{\it Многие математики осознавали серьёзные последствия сложившегося положения: терялась часть культуры, 
 математика теряла потенциальные таланты; но изменить его было очень нелегко. Создание новых учебников было в своем роде
 тоже проектом поворота, только не рек, а преподавания математики, и поддержка его была столь же массированной.
 Сломить её удалось Льву Семёновичу. Больше года он боролся, сочетая бурный натиск и дипломатию, пока, наконец,
 ему не удалось опубликовать свои взгляды в журнале «Коммунист». В этом ему очень помогли журналисты. В конце концов дорога
 к пересмотру учебников была открыта. В их реанимации потом участвовало много людей, но прорыв, право на её проведение завоевал
 Л.~С.~Понтрягин. Мне кажется, {\bf у нас пока нет ясных и прозрачных учебников по математике, какие были несколько десятилетий назад},
 но всё же «геркулесовы столбы» формалистического творчества остались позади, положение в какой-то степени нормализуется,
 и результат сказывается на десятках миллионов подростков, на их умах и душах.
 \newline 
 Лишь после его смерти я узнал, что одним из любимых произведений Понтрягина было жизнеописание Бенвенуто Челлини, и это сделало мне более понятным Льва Семёновича. У обоих было действительно много общего, и мне представляется, что, если бы почтенный академик жил во времена Челлини, 
 он тоже многие споры решал бы при помощи собственной шпаги и уж, во всяком случае, не обращался бы к помощи наёмных убийц.}]

}	

   \hangindent=0.3cm \noindent	
[Шевк] Шевкин А.В.  \href{http://www.shevkin.ru/novosti/nazad-k-kiselyovu} {\it Назад к Киселеву!}, 04.11.2012;
\newline
\href{http://www.shevkin.ru/novosti/eshhyo-raz-pro-vozvrashhenie-k-uchebnikam-a-p-kiseleva/}
{\it Ещё раз про возвращение к учебникам А.П. Киселева}, 07.10.2016

   \hangindent=0.3cm \noindent
[Шевч1966]
Шевченко И.Н. \href{http://mat.univie.ac.at/~neretin/misc/reform/Shevchenko.djvu}
{\it Арифметика для 5-6 классов}, 1966, Просвещение,  издание 11.

    \hangindent=0.3cm \noindent
[Шир2003] Ширяев А.Н. {\it Жизнь и творчество} [Колмогорова],  В книге
{\it Колмогоров. Юбилейное издание в трех книгах. Книга первая. Истина - благо.} Редактор-составитель А.Н.Ширяев.
Москва, Физматлит, 2003 (стр. 17-210)

   \hangindent=0.3cm \noindent
[vis2012]
\href{http://www.diary.ru/~vis1952}{\it Проект Колмогоров.} (текст неизвестного автора)

\label{bib-end}

\bigskip 

\tt 
\noindent
Yury Neretin\\
Pauli Institute/c.o Math. Dept., University of Vienna\\
\&Institute for Theoretical and Experimental Physics (Moscow); \\
\&MechMath Dept., Moscow State University;\\
\&Institute for Information Transmission Problems;\\
URL: http://mat.univie.ac.at/$\sim$neretin/	
 
\end{document}